\documentclass[11pt,a4paper]{article}
\usepackage{amsmath, amsfonts, amssymb,amsthm}
\usepackage{a4wide}
\usepackage{parskip}
\usepackage{enumitem}
\setlist[itemize,1]{label=$\bullet$}
\usepackage{etoolbox}
\usepackage{xcolor}
\usepackage[hypertexnames=false]{hyperref}
\usepackage{latexsym}
\usepackage{cite}

\usepackage[left=1in, right=1in,top=1in,bottom=1in]{geometry}

\allowdisplaybreaks

\newcommand{\be} {\begin{equation}}
	\newcommand{\ee} {\end{equation}}
\newcommand{\bea} {\begin{eqnarray}}
	\newcommand{\eea} {\end{eqnarray}}
\newcommand{\Bea} {\begin{eqnarray*}}
	\newcommand{\Eea} {\end{eqnarray*}}

\newcommand{\de} {\beta}

\def\R{{\mathbb R}}

\def\H{{H_0^1(\Omega)}}
\def\I{I_{\lambda}}
\def\J{J_{\lambda}}

\def\p{{\prime}}

\newcommand{\na} {\nabla}

\numberwithin{equation}{section}

\newtheorem{definition}{Definition}[section]
\newtheorem{theorem}{Theorem}[section]
\newtheorem{rem}{Remark}[section]
\newtheorem{lemma}{Lemma}[section]
\makeatletter
\patchcmd{\proof}{\itshape}{\normalfont}{}{}
\patchcmd{\proof}{\@addpunct{.}}{\@addpunct{.}\normalfont}{}{}
\makeatother
\newtheorem{prop}{Proposition}[section]
\newtheorem{ex}{Example}[section]
\numberwithin{equation}{section}
\def\proof{\noindent{\textbf{Proof. }}}

\begin{document}
	\date{}
	\title{A Critical Choquard Problem for a Quasilinear Schr\"odinger Operator with Singular Discontinuous Nonlinearity}
	\author{ {\bf Ankit Mishra$\,^{1,}$\footnote{e-mail: {\tt ankitmishra.rs.mat23@itbhu.ac.in, ankitmishra18295@gmail.com}}}, \bf Sarika Goyal$\,^{2,}$\footnote{e-mail: {\tt sarika@nsut.ac.in, sarika1.iitd@gmail.com}}, \bf Divya Goel$\,^{1,}$ \footnote{e-mail: {\tt divya.mat@iitbhu.ac.in},\\
     \tt Corresponding Author: Divya Goel }
     \\ $^1\,$ Department of Mathematical Sciences,\\ Indian Institute of Technology (BHU), India\\
   $^2\,$ Department of Mathematics,\\ Netaji Subhash University of Technology, India }

	\maketitle

	\begin{abstract}
In this article, we study the quasilinear Schr\"odinger equation
\begin{align*} \begin{cases}-\Delta u - \Delta \bigl(|u|^{2\alpha}\bigr)\,|u|^{2\alpha-2}u &=\lambda \left(\chi_{\{u<b\}} u^{-\beta}+\left(\int_\Omega \frac{H(y,u)}{\left|x-y\right|^\mu} dy\right)h(x,u)\right) \text{ in } \Omega,\\&u>0 \text{ in } \Omega, \\&u=0 \text{ on } \partial\Omega, \end{cases} \end{align*}
where $\Omega\subset\R^2$ is a smooth bounded domain, $b>0$, $0<\mu<2$, $\alpha>\frac12$, $0<\beta<3$, and $h$ has critical exponential growth with primitive $H$. We show that there is $\lambda^*_b\in(0,\infty)$ such that the problem has a positive solution if and only if $\lambda\in(0,\lambda^*_b]$. For $\lambda<\lambda^*_b$, this solution is $g(v_\lambda)$, where the change of variables $g$ makes the problem semilinear and $v_\lambda$ is a local minimizer of the transformed energy in $H^1_0(\Omega)$. If $h$ satisfies a lower bound at the critical scale $\exp\bigl(2\alpha\,g^{-1}(u)^2\bigr)$, we obtain a second positive solution for every $\lambda\in(0,\lambda^*_b)$. The level estimate for this solution must be done at this scale, and not at the scale $e^{u^{4\alpha}}$ used for semilinear problems. This bound excludes $h=f\,e^{u^{4\alpha}}$ with $f$ of polynomial growth. For $\alpha>1$, we obtain a second solution for $h=u^{4\alpha-1}e^{u^{4\alpha}}$, provided that the maximum of the first solution exceeds a constant depending only on $\alpha$. For $\alpha\leq1$, this question remains open.

\medskip
\noindent \textbf{Keywords:} 
	\noindent Modified quasilinear Schr\"odinger equation; Choquard equation; critical exponential growth; singular nonlinearity; discontinuous nonlinearity; multiplicity of positive solutions.
    
	\end{abstract}
	\section{Introduction}
	
We consider the following problem:
\begin{equation}
\begin{cases}
-\Delta u - \Delta \bigl(|u|^{2\alpha}\bigr)\,|u|^{2\alpha-2}u &=\lambda \left(\chi_{\{u<b\}} u^{-\beta}+\left(\int_\Omega \frac{H(y,u)}{\left|x-y\right|^\mu} dy\right)h(x,u)\right) \text{ in } \Omega,\\&u>0 \text{ in } \Omega, \\&u=0 \text{ on } \partial\Omega,
\end{cases}
\tag{$P_{\lambda}^b$}
\label{sd1}
\end{equation}
where $\Omega \subset \R^2$ is a bounded domain with smooth boundary $\partial\Omega$, $b>0, 0<\mu<2$, $\alpha>\frac{1}{2}$ and $0<\beta<3$. Here, $\chi$ is the characteristic function,
\begin{align*}
h(x,s)=f(x,s)\exp\bigl(|s|^{4\alpha}\bigr),\qquad H(x,s)=\int_0^sh(x,t)\,dt,
\end{align*}
and $f$ satisfies the hypotheses $(f_1)$--$(f_5)$ stated below.

Four features of \eqref{sd1} shape the analysis: the quasilinear operator, the singular term $u^{-\beta}$ in the full range $0<\beta<3$, the jump of the nonlinearity across the level $u=b$, and a nonlocal Choquard term whose nonlinearity $h$ has critical exponential growth in the sense of Trudinger and Moser.

Quasilinear Schr\"odinger equations containing the term $\Delta\bigl(\ell(|\psi|^2)\bigr)\ell^\p(|\psi|^2)\psi$ arise in plasma physics and condensed matter theory \cite{kurihara1981,brandi1993}. For $\ell(s)=s^{\alpha}$, this term equals $\alpha\,\Delta\bigl(|\psi|^{2\alpha}\bigr)|\psi|^{2\alpha-2}\psi$. Substituting a standing wave $\psi(t,x)=e^{-iEt}u(x)$ into the corresponding Schr\"odinger equation gives an elliptic equation for $u$ whose quasilinear part is, up to a constant factor, the operator in \eqref{sd1}.
 The existence theory for these equations began with the minimax approach of \cite{poppenberg2002} and the change of variables introduced independently in \cite{liu2003} and \cite{colin2004}. The singular, discontinuous and nonlocal terms of \eqref{sd1} are not part of this physical model. They come from the Dirichlet problems studied in \cite{dhanya2014,anthal2025,biswas2024}, and the purpose of this paper is to understand how they interact with the quasilinear operator.

The change of variables makes \eqref{sd1} accessible to variational methods. The natural energy of \eqref{sd1}, given in \eqref{sd3} below, contains the term $\int_\Omega|u|^{4\alpha-2}|\nabla u|^2dx$, which need not be finite for $u\in H^1_0(\Omega)$. Hence no functional of class $C^1$ is available on the natural space. Following \cite{liu2003,colin2004}, we write $u=g(v)$, where $g$ is the odd function defined by
\begin{align}\label{sd82}
g^\p(t)=\frac{1}{\bigl(1+2\alpha|g(t)|^{4\alpha-2}\bigr)^{1/2}}\quad\text{for }t>0,\qquad g(0)=0,\qquad g(-t)=-g(t).
\end{align}
 This transforms \eqref{sd1} into the semilinear problem \eqref{sd6} for $v$, in the sense that $u=g(v)$ is a weak solution of \eqref{sd1} if and only if $v$ is a weak solution of \eqref{sd6} (Lemma~\ref{lemsd2}). The change of variables also fixes the scale of critical growth. In dimension two, the maximal growth compatible with $H^1_0(\Omega)$ is exponential, by the Trudinger--Moser inequality \cite{trudinger1967,moser1971}. For the transformed problem, the critical scale is $\exp(2\alpha v^2)$. By Lemma~\ref{lemsd22},
\[
g(v)^{4\alpha}=2\alpha v^2-\Xi(v),\qquad \Xi(v)=o(v^2)\quad\text{as }v\to\infty,
\]
so a nonlinearity of order $\exp(|u|^{4\alpha})$ has critical growth in the new variable as well. This is the growth prescribed by $(f_2)$ below. For $\alpha=1$, it is the critical growth identified in \cite{doo2010cv}. However, $\Xi(v)\to\infty$, so $e^{|u|^{4\alpha}}$ stays below $e^{2\alpha v^2}$ by the factor $e^{-\Xi(v)}$, which tends to zero. This matters for the second solution, and our condition $(f_5)$ is a lower bound on $H$ at the scale $e^{2\alpha v^2}$, see Remark~\ref{remsd7} and Example~\ref{exsd1}.

We summarize the results closest to \eqref{sd1}. We do not attempt a complete review, and we refer to \cite{ghergu2008,moroz2017} and the references therein.

Singular problems $-\Delta u=f(x,u)$ with $f(x,u)\sim u^{-\beta}$ near $u=0$ go back to \cite{crandall1977dirichlet}, where existence was proved by the method of sub- and supersolutions. For the pure power $-\Delta u=p(x)u^{-\beta}$, Lazer and McKenna \cite{lazer1991} showed that the solution belongs to $H^1_0(\Omega)$ if and only if $\beta<3$. This is the origin of the restriction $0<\beta<3$ in \eqref{sd1}. Multiplicity for singular problems with a critical perturbation was obtained in \cite{hirano2008,giacomoni2007,giacomoni2009}. A key tool there is the comparison between local minimizers in $H^1_0(\Omega)$ and in $C^1$, which goes back to Brezis and Nirenberg \cite{brezis1993} and was extended to singular problems in \cite{giacomoni2007,hirano2008}. In dimension two, elliptic problems with critical exponential growth were treated with the Trudinger--Moser inequality in \cite{defigueiredo1995}. Singular problems with critical exponential growth were studied in \cite{adimurthi2006}.

Discontinuous nonlinearities make the energy functional merely locally Lipschitz. Their variational treatment rests on Clarke's generalized gradient and on its application to elliptic equations by Chang \cite{chang1980}, and critical growth was treated in \cite{badiale1997}. Such nonlinearities arise in free boundary problems, see \cite{chang1980} and the references therein. The combination of a singular term with the jump $\chi_{\{u<b\}}$ was treated in \cite{dhanya2014} for $-\Delta u=\lambda(\chi_{\{u<b\}}u^{-\beta}+h(u)e^{u^2})$ in a bounded domain of $\R^2$ with $0<\beta<3$, and in \cite{dhanya2016} for Sobolev-critical growth in $\R^N$. 
When the discontinuous part of the nonlinearity is monotone, as it is here, the translated energy splits as $\Lambda+\psi$, where $\Lambda$ is of class $C^1$ and $\psi$ is convex and continuous. This is the setting of Szulkin \cite{szulkin1986}, and in Section~4 we use his notion of critical point.

Choquard nonlinearities $\bigl(|x|^{-\mu}*H(u)\bigr)h(u)$ originate in the polaron model of Pekar, and they were studied by Lieb \cite{lieb1977}. The basic tool is the Hardy--Littlewood--Sobolev inequality \cite[Theorem 4.3]{lieb2001}; we refer to \cite{moroz2017} for a survey. For Choquard nonlinearities of critical exponential growth in $\R^2$, the compactness analysis is due to \cite{alves2016}. Singular Choquard problems with critical growth were treated in \cite{giacomoni2021}. The Choquard problem with a singular discontinuous nonlinearity was studied in \cite{anthal2025}, for the Laplacian in $\R^N$, $N\geq3$, with Sobolev-critical growth.

For the quasilinear operator, critical exponential growth was treated in \cite{doo2007,doo2010cv}, without singular terms. Quasilinear equations with a singular term and critical exponential growth were studied in \cite{biswas2024}, for $\alpha=1$ and a local nonlinearity. Quasilinear Schr\"odinger equations with nonlocal terms of exponential growth, without singular terms, were studied in \cite{biswas2024jga,malhotra2025}. To our knowledge, no result is available for the quasilinear operator in the presence of a singular term, a jump discontinuity and a critical Choquard nonlinearity at the same time, which is the situation of \eqref{sd1}.

The scheme of the proof is that of \cite{dhanya2014,anthal2025}, but several steps have to be redone. We list what is new with respect to \cite{dhanya2014} (semilinear, local nonlinearity), \cite{anthal2025} (Laplacian, $N\geq3$, Sobolev-critical Choquard term) and \cite{biswas2024} (quasilinear, $\alpha=1$, local nonlinearity, no jump).
\begin{enumerate}[label=\textnormal{(\arabic*)}]
\item The singular nonlinearity after the change of variables is $g(v)^{-\beta}g^\p(v)$, which is nonincreasing but not a pure power. Lemma~\ref{lemsd3} provides the boundary behavior, the comparison principle and the existence theory for this class, with constants that do not depend on the regularization parameter. Lemma~\ref{lemsd20} is a variational form of Perron's method for singular nonlinearities, for which the standard statement \cite[Chapter I, Theorem 2.4]{struwe1996} does not apply, because the nonlinearity is unbounded near $\partial\Omega$.
\item In the Brezis--Nirenberg type comparison argument (Lemma~\ref{lemsd12}), the nonlocal term cannot be estimated pointwise. Its first-order part is cancelled against the equation of the supersolution, and the remaining second-order part is absorbed by the Dirichlet energy by means of the Trudinger--Moser inequality. This is where the global ordering of the functions being compared is needed.
\item The level set $\{g(v)=b\}$ of a solution, and the set $\{w=0\}$ for the second solution, have measure zero because the nonlocal term is strictly positive (Lemmas~\ref{lemsd10}, \ref{lemsd21} and \ref{lemsd14}). The same positivity gives the existence of a solution at $\lambda=\lambda^*_b$.
\item The dichotomy of Lemma~\ref{lemsd17} avoids the separate treatment of the ``zero altitude'' case of \cite{dhanya2014,anthal2025}.
\item The growth condition $(f_5)$ for the second solution is stated for the transformed nonlinearity, at its Trudinger--Moser scale $\exp(2\alpha v^2)$. It is the analogue of condition (H3) of \cite{dhanya2014}, and it has no threshold, because the translation by the first solution produces an extra factor in the level estimate (Lemma~\ref{lemsd18}). In the original variable, $(f_5)$ requires $f$ to compensate the factor $e^{-\Xi(v)}$, where $\Xi(v)=2\alpha v^2-g(v)^{4\alpha}$ is computed in Lemma~\ref{lemsd22}. Hence $(f_5)$ excludes every $f$ of polynomial growth (Remark~\ref{remsd7}). For $\alpha>1$, $\Xi$ grows only linearly, and Proposition~\ref{propsd2} gives a second solution also when $f$ is a positive power of $u$, under a condition on the size of the first solution.
\end{enumerate}

We now state the hypotheses on $f$. Recall that $h(x,s)=f(x,s)\exp(|s|^{4\alpha})$ and $H(x,s)=\int_0^sh(x,t)\,dt$, and that $g$ is defined in \eqref{sd82}.
\begin{itemize}
\item[$(f_1)$] $f\in C^1(\bar\Omega\times\R)$, and for every $x\in\bar\Omega$, $f(x,s)>0$ for $s>0$ and $f(x,s)=0$ for $s\leq0$. In particular, $f(x,0)=0$ and $h(x,s)=H(x,s)=0$ for $s\leq0$.
\item[$(f_2)$] For any $\epsilon>0$, $\displaystyle \lim_{t\to\infty} \sup_{x\in \bar{\Omega}} f(x,t) \exp(-\epsilon \left|t\right|^{4\alpha})=0$ and  $\displaystyle \lim_{t\to\infty} \inf_{x\in \bar{\Omega}} f(x,t) \exp(\epsilon \left|t\right|^{4\alpha})=\infty.$
\item[$(f_3)$] There exists a constant $M>0$ such that $H(x,s) < M(1+h(x,s))$ for all $x\in\bar\Omega$ and $s>0$.
\item[$(f_4)$] For every $x\in\bar\Omega$, the function $s \mapsto h(x,s)s^{1-4\alpha}$ is nondecreasing on $(0,\infty)$.
\item[$(f_5)$] There exists $\theta\in[0,1)$ such that
\begin{align*}
\liminf_{s\to\infty}\,\inf_{x\in\bar\Omega}\,H(x,g(s))\exp\bigl(-2\alpha s^2+s^\theta\bigr)>0 .
\end{align*}
\end{itemize}

Condition $(f_2)$ says that $h$ has critical exponential growth $\exp(|s|^{4\alpha})$. By Remark~\ref{remsd7}(a), this is the same as Trudinger--Moser critical growth $\exp(2\alpha s^2)$ for $h(x,g(s))$. Condition $(f_3)$ is an Ambrosetti--Rabinowitz type condition adapted to exponential growth, as in \cite{defigueiredo1995}. Condition $(f_4)$ implies that $s\mapsto h(x,s)=\bigl(h(x,s)s^{1-4\alpha}\bigr)s^{4\alpha-1}$ is nondecreasing on $(0,\infty)$, and that $h(x,s)\leq h(x,1)s^{4\alpha-1}$ for $0<s\leq1$. In particular, $f(x,s)=o(s)$ as $s\to0$. Together with Lemma~\ref{lemsd1}$(g_8)$, $(f_4)$ implies that $t\mapsto h(x,g(t))g^\p(t)$ is nondecreasing on $(0,\infty)$, a fact used in all the comparison arguments.

Condition $(f_5)$ is used only for the second solution, and it is the analogue, for the transformed nonlinearity, of condition (H3) in \cite{dhanya2014}. Example~\ref{exsd1} shows that $(f_1)$--$(f_5)$ are compatible. After Theorem~\ref{thmsd1}, we write $(f_5)$ in the original variable.

\begin{ex}\label{exsd1}
(a) Let $q\geq4\alpha-1$, and let $h(x,u)=u^{q}\exp\bigl(2\alpha\,g^{-1}(u)^2\bigr)$ for $u>0$ and $h(x,u)=0$ for $u\leq0$. Thus $f(x,u)=u^{q}e^{\Xi(g^{-1}(u))}$ for $u>0$, with $\Xi$ as in Lemma~\ref{lemsd22}. By \eqref{sd84}, $f(x,u)=u^{q}\exp\bigl(c_\alpha\,\omega(u)(1+o(1))\bigr)$ as $u\to\infty$. For $\alpha=1$, this reads $h(x,u)=u^{q}\exp\bigl(u^4+u^2\log u\,(1+o(1))\bigr)$.

The function $u\mapsto\Xi(g^{-1}(u))=2\alpha g^{-1}(u)^2-u^{4\alpha}$ is of class $C^1$ on $[0,\infty)$ and vanishes at $0$. Since $q>1$, $(f_1)$ holds. Since $\Xi(s)=o(s^2)$ and $2\alpha s^2=g(s)^{4\alpha}+\Xi(s)$, we have $\Xi(g^{-1}(u))=o(u^{4\alpha})$, which gives the first limit in $(f_2)$. The second limit holds because $\Xi\geq0$. Next, $(\log h)^\p(x,u)=q/u+4\alpha g^{-1}(u)(1+2\alpha u^{4\alpha-2})^{1/2}\to\infty$ as $u\to\infty$, so $H(x,u)/h(x,u)\to0$ by l'H\^opital's rule, and $(f_3)$ follows. The function $h(x,u)u^{1-4\alpha}=u^{q+1-4\alpha}e^{2\alpha g^{-1}(u)^2}$ is nondecreasing, so $(f_4)$ holds. Finally, the substitution $\tau=g(\sigma)$ gives
\begin{align*}
H(x,g(s))=\int_0^sg(\sigma)^{q}g^\p(\sigma)e^{2\alpha\sigma^2}d\sigma .
\end{align*}
By Lemma~\ref{lemsd1}$(g_8)$ with $\gamma=q$, $g(\sigma)^qg^\p(\sigma)\geq g(1)^qg^\p(1)\sigma$ for $\sigma\geq1$. Let $s\geq2$. On $[s-1/s,s]$, we have $\sigma\geq s/2\geq1$ and $\sigma^2\geq s^2-2$, hence
\begin{align*}
H(x,g(s))\geq\int_{s-1/s}^sg(\sigma)^{q}g^\p(\sigma)e^{2\alpha\sigma^2}d\sigma\geq\frac{g(1)^qg^\p(1)}{2e^{4\alpha}}\,e^{2\alpha s^2},
\end{align*}
and $(f_5)$ holds with $\theta=0$. In this example, $h$ exceeds $u^{q}e^{u^{4\alpha}}$ by the unbounded factor $e^{\Xi(g^{-1}(u))}$.

(b) Let $f(x,s)=s^{4\alpha-1}e^{s^{p}}$ for $s>0$ and $f(x,s)=0$ for $s\leq0$, with $\max\{2,2\alpha\}<p<4\alpha$. Then $f\in C^1(\bar\Omega\times\R)$ because $4\alpha-1>1$, so $(f_1)$ holds. Condition $(f_2)$ holds because $p<4\alpha$. Condition $(f_3)$ holds because $H(x,s)/h(x,s)\to0$ as $s\to\infty$, by l'H\^opital's rule and the bound $(\log h)^\p(x,s)\geq 4\alpha s^{4\alpha-1}\to\infty$. The function $h(x,s)s^{1-4\alpha}=e^{s^p+s^{4\alpha}}$ is nondecreasing, so $(f_4)$ holds. Finally, since $p>2\alpha$ and $p>2$, $f(x,u)e^{-c\,\omega(u)}\to\infty$ for every $c>0$, and $(f_5)$ holds by Remark~\ref{remsd7}(b).
\end{ex}

\subsection{Preliminaries}\label{secsd1}
In this subsection, we recall the change of variables, introduce the functional setting and the notion of solution, and fix the notation.

\begin{definition}\label{defsd1}
Let $\phi_1$ be the first positive eigenfunction of $-\Delta$ in $\H$ with $\|\phi_1\|_{L^\infty(\Omega)}  <1.$ For $\beta>0,$ we define $\varphi_{\beta}$ as follows:
\[
\varphi_{\beta}=
\begin{cases}
\phi_1 & \text{ for } 0<\beta<1, \\ \phi_1(-\log \phi_1)^\frac{1}{2} & \text{ for } \beta=1, \\ \phi_1^{\frac{2}{\beta +1}} & \text{ for } \beta>1.
\end{cases}
\]
\end{definition}
One checks easily that $\varphi^{-\beta}_\beta \in L^\infty_{loc}(\Omega)$, and that $\varphi^{-\beta}_\beta \in L^1(\Omega)$ if and only if $\beta \in (0,1)$. We write $d(x):=\operatorname{dist}(x,\partial\Omega)$ for $x\in\bar\Omega$.

The energy functional associated with \eqref{sd1} is
\begin{equation}
{\label{sd3}} \I(u)=\frac{1}{2} \int_\Omega \left(1+2 \alpha \left|u\right|^{2(2\alpha-1)}\right) \left|\nabla u\right|^2\ dx -\lambda \int_\Omega G(u)dx-\frac{\lambda}{2} \int_\Omega \int_\Omega \frac{H(y,u(y)) H(x,u(x))}{\left|x-y\right|^\mu}dxdy,
\end{equation}
where $G$ is the primitive of $s\mapsto\chi_{\{s<b\}}s^{-\beta}$ given in \eqref{sd5} below. This functional is not well defined on $\H$, for two reasons: the term $\int_\Omega|u|^{4\alpha-2}|\nabla u|^2dx$ is not finite for every $u\in\H$, and the singular term need not be integrable. We therefore use the change of variables $u=g(v)$, with $g$ as in \eqref{sd82}. The properties of $g$ that we need are collected in the next lemma. For $\alpha=1$, $(g_1)$, $(g_2)$, $(g_4)$ and $(g_5)$ are proved in \cite[Lemmas~2.1 and~2.2]{colin2004}. We need them for every $\alpha>\frac12$, and we give the short proofs in Appendix~\ref{secsd3}.
\begin{lemma}{\label{lemsd1}}
The function $g$ satisfies the following properties:
\begin{itemize}
\item[$(g_1)$] $g$ is uniquely defined and invertible. Also, $g \in C^{1}(\R)\cap C^{\infty}(\R\setminus\{0\})$;
\item[$(g_2)$] $g(0)=0$, \quad $0<g^{\prime}(t) \leq 1, \text{ for all } t\in \R ;$
\item[$(g_3)$] $\left|g(t)\right| \leq \left|t\right|$; $\left|g(t)\right|^{2\alpha} \leq (2\alpha)^{\frac{1}{2}} \left|t\right|$, for all $t\in \R$;
\item[$(g_4)$] $ \left \{
\begin{array}{ll}
\frac{1}{2} g(t)\leq \alpha t g^{\prime}(t) \leq \alpha g(t), & t>0, \\ \frac{1}{2} g(t)\geq \alpha t g^{\prime}(t) \geq \alpha g(t), & t<0;
\end{array}
\right. $
\item[$(g_5)$] $\displaystyle \lim_ {t\to \infty}\frac{g(t)}{t^{\frac{1}{2\alpha}}}=(2\alpha)^{\frac{1}{4\alpha}}$;
\item[$(g_6)$] $\left|g(t)\right| \geq
\begin{cases}
g(1)\left|t\right|, & \text{ if } \left|t\right| \leq 1,\\ g(1)\left|t\right|^\frac{1}{2\alpha}, & \text{ if } \left|t\right| \geq 1;
\end{cases}$
\item[$(g_7)$] $\left|g(t)\right|^{2\alpha -1} g^{\prime}(t) \leq \frac{1}{(2\alpha)^{\frac{1}{2}}},$ for all $t\in \R;$
\item[$(g_{8})$] the function $\frac{(g(t))^\gamma g^{\prime}(t)}{t}$ is strictly increasing for $\gamma \geq 4\alpha -1$ and $t>0$;
\item[$(g_{9})$] the function $g(t)^{-\gamma} g^{\prime}(t)$ is decreasing for all $t>0$, where $\gamma >0.$
\end{itemize}
\end{lemma}

The next lemma measures the gap between the scales $\exp(|u|^{4\alpha})$ and $\exp(2\alpha v^2)$.
\begin{lemma}\label{lemsd22}
For $s\geq0$, let $\Xi(s):=2\alpha s^2-g(s)^{4\alpha}$. Define
\begin{align*}
\omega(u):=
\begin{cases}
u^{2\alpha}, & \alpha>1,\\ u^2\log u, & \alpha=1,\\ u^2, & \tfrac12<\alpha<1,
\end{cases}
\qquad
c_\alpha:=
\begin{cases}
2\delta_\infty, & \alpha>1,\\ 1, & \alpha=1,\\ \frac{1}{2(1-\alpha)}, & \tfrac12<\alpha<1,
\end{cases}
\end{align*}
where, for $\alpha>1$,
\begin{align}\label{sd83}
\delta_\infty:=\int_0^\infty k_\alpha(r)\,dr<\infty,\qquad k_\alpha(r):=\frac{(2\alpha)^{1/2}}{(1+2\alpha r^{4\alpha-2})^{1/2}+(2\alpha)^{1/2}r^{2\alpha-1}} .
\end{align}
Then $\Xi\geq0$ on $[0,\infty)$, $\Xi(s)=o(s^2)$ as $s\to\infty$, and
\begin{align}\label{sd84}
\lim_{s\to\infty}\frac{\Xi(s)}{\omega(g(s))}=c_\alpha .
\end{align}
\end{lemma}
\begin{proof}
For $s=0$ all the quantities below vanish, so let $s>0$ and $u:=g(s)>0$. By \eqref{sd82}, $(g^{-1})^\p(r)=(1+2\alpha r^{4\alpha-2})^{1/2}$ for $r>0$ and $g^{-1}(0)=0$, so that $s=g^{-1}(u)=\int_0^u(1+2\alpha r^{4\alpha-2})^{1/2}dr$. Consequently,
\begin{align*}
\delta(s):=(2\alpha)^{1/2}s-u^{2\alpha}=\int_0^u\Bigl((2\alpha)^{1/2}(1+2\alpha r^{4\alpha-2})^{1/2}-2\alpha r^{2\alpha-1}\Bigr)dr=\int_0^uk_\alpha(r)\,dr ,
\end{align*}
where we multiplied and divided the integrand by $(2\alpha)^{1/2}(1+2\alpha r^{4\alpha-2})^{1/2}+2\alpha r^{2\alpha-1}$. Thus $\delta(s)\geq0$ and
\begin{align*}
\Xi(s)=\bigl((2\alpha)^{1/2}s\bigr)^2-u^{4\alpha}=\delta(s)\bigl(2u^{2\alpha}+\delta(s)\bigr)\geq0 .
\end{align*}
As $r\to\infty$, $2r^{2\alpha-1}k_\alpha(r)\to1$. For $\alpha>1$, it follows that $k_\alpha$ is integrable on $(0,\infty)$. For $\alpha\leq1$, we apply l'H\^opital's rule to $\delta(s)=\int_0^uk_\alpha(r)\,dr$, regarded as a function of $u=g(s)$. Hence, as $u\to\infty$,
\begin{align*}
\delta(s)\to\delta_\infty\ \ (\alpha>1),\qquad \frac{\delta(s)}{\log u}\to\frac12\ \ (\alpha=1),\qquad \frac{\delta(s)}{u^{2-2\alpha}}\to\frac{1}{4(1-\alpha)}\ \ (\tfrac12<\alpha<1),
\end{align*}
with $\delta_\infty<\infty$ as in \eqref{sd83}. In each case $\delta(s)=o(u^{2\alpha})$, because $2-2\alpha<2\alpha$ for $\alpha>\frac12$. Therefore $\Xi(s)=2u^{2\alpha}\delta(s)(1+o(1))$, and dividing by $\omega(u)$ gives \eqref{sd84}. Finally, $\omega(u)=o(u^{4\alpha})$ and $u^{4\alpha}\leq2\alpha s^2$ by $(g_3)$, so $\Xi(s)=o(s^2)$. \qed
\end{proof}

\begin{rem}\label{remsd7}
\begin{enumerate}[label=\textnormal{(\alph*)}]
\item Since $h(x,g(s))=f(x,g(s))e^{2\alpha s^2-\Xi(s)}$ and $\Xi(s)=o(s^2)$, condition $(f_2)$ is equivalent to the following: for every $\epsilon>0$,
\begin{align*}
\lim_{s\to\infty}\sup_{x\in\bar\Omega}h(x,g(s))e^{-(2\alpha+\epsilon)s^2}=0,\qquad \lim_{s\to\infty}\inf_{x\in\bar\Omega}h(x,g(s))e^{-(2\alpha-\epsilon)s^2}=\infty .
\end{align*}
Thus $(f_2)$ says that the transformed nonlinearity has Trudinger--Moser critical growth with exponent $2\alpha$.
\item Let $\omega$ and $c_\alpha$ be as in Lemma~\ref{lemsd22}. If $\inf_{x}f(x,u)\geq c_fe^{c\,\omega(u)}$ for all large $u$, with $c_f>0$ and $c>c_\alpha$, then $(f_5)$ holds. If $\sup_xf(x,u)\leq c_fe^{c\,\omega(u)}$ for all large $u$, with $c_f>0$ and $c<c_\alpha$, then $(f_5)$ fails. The case $c=c_\alpha$ is not decided by this argument. Example~\ref{exsd1}(a) belongs to this case, and it satisfies $(f_5)$. In particular, $(f_5)$ excludes every $f$ of polynomial growth.
\end{enumerate}
\end{rem}
\begin{proof}
(a) Put $u=g(s)$. By Lemma~\ref{lemsd22}, $u^{4\alpha}=2\alpha s^2(1+o(1))$, and $h(x,g(s))e^{-(2\alpha\pm\epsilon)s^2}=f(x,u)e^{\mp\epsilon s^2-\Xi(s)}$. The equivalence follows from these two facts.

(b) Let $s>0$ and $u=g(s)$, so that $2\alpha s^2=u^{4\alpha}+\Xi(s)$. Suppose first that $f(x,\tau)\geq c_fe^{c\,\omega(\tau)}$ for all $x$ and $\tau\geq u_0$, with $c>c_\alpha$. Let $u\geq\max\{u_0+1,2\}$. In the integral defining $H(x,u)$ we keep only the interval $[u-u^{1-4\alpha},u]\subset[u-1,u]$. On this interval, $\tau^{4\alpha}\geq(u-u^{1-4\alpha})^{4\alpha}\geq u^{4\alpha}-4\alpha$ by Bernoulli's inequality, and $\omega(\tau)\geq\omega(u-1)$ because $\omega$ is increasing on $[1,\infty)$. This gives
\begin{align*}
H(x,u)\geq\int_{u-u^{1-4\alpha}}^uf(x,\tau)e^{\tau^{4\alpha}}d\tau\geq c_fe^{-4\alpha}u^{1-4\alpha}e^{u^{4\alpha}+c\,\omega(u-1)} .
\end{align*}
Hence
\begin{align*}
H(x,g(s))e^{-2\alpha s^2+s^\theta}\geq c_fe^{-4\alpha}u^{1-4\alpha}\exp\bigl(c\,\omega(u-1)-\Xi(s)+s^\theta\bigr).
\end{align*}
By \eqref{sd84} and the limit $\omega(u-1)/\omega(u)\to1$, the exponent is at least $\epsilon_0\omega(u)$ for some $\epsilon_0>0$ and all large $u$. Since $\omega(u)\geq u^{2\alpha}$ for $u\geq e$, the factor $e^{\epsilon_0\omega(u)}$ dominates $u^{1-4\alpha}$, and $(f_5)$ holds with any $\theta\in[0,1)$.

Conversely, suppose that $\sup_xf(x,u)\leq c_fe^{c\,\omega(u)}$ for large $u$, with $c<c_\alpha$. By $(f_3)$, $H(x,u)\leq M(1+f(x,u)e^{u^{4\alpha}})$, hence
\begin{align*}
H(x,g(s))e^{-2\alpha s^2+s^\theta}\leq Me^{-2\alpha s^2+s^\theta}+Mc_f\exp\bigl(c\,\omega(u)-\Xi(s)+s^\theta\bigr).
\end{align*}
By \eqref{sd84}, $c\,\omega(u)-\Xi(s)\leq-\epsilon_0\omega(u)$ for some $\epsilon_0>0$ and all large $s$. Moreover, $s\leq Cu^{2\alpha}$ by $(g_5)$, so $s^\theta\leq Cu^{2\alpha\theta}=o(\omega(u))$ for every $\theta\in[0,1)$. Hence the right-hand side tends to $0$, and $(f_5)$ fails. The last assertion follows, since every $f$ of polynomial growth satisfies the second condition with any $c\in(0,c_\alpha)$. \qed
\end{proof}

After the change of variables $u=g(v)$ in \eqref{sd3}, we obtain the functional
\begin{equation}{\label{sd4}}
\J(v)=\frac{1}{2} \int_\Omega \left|\nabla v\right|^2\ dx -\lambda \int_\Omega G(g(v))dx-\frac{\lambda}{2} \int_\Omega \int_\Omega \frac{H(y,g(v)(y)) H(x,g(v)(x))}{\left|x-y\right|^\mu}dxdy,
\end{equation}
where
\begin{equation}{\label{sd5}}
G(u)=
\begin{cases}
-\infty \ (\beta\geq1),\ 0 \ (\beta<1), & \text{ if } u\leq 0, \\ \frac{u^{1-\beta}}{1-\beta}, & \text{ if } 0<u<\frac{b}{2}, \\ \frac{1}{1-\beta}(\frac{b}{2})^{1-\beta}+ \int^u_\frac{b}{2} \chi_{\{t<b\}}(t) t^{-\beta}dt, & \text{ if } u \geq \frac{b}{2},
\end{cases}
\end{equation}
for $0<\beta<3$, $\beta \neq 1$. For $\beta=1$, the terms $(1-\beta)^{-1}t^{1-\beta}$ are replaced by $\log t$. For $\beta\geq1$, the choice $G=-\infty$ on $(-\infty,0]$ gives $J_\lambda(v)=+\infty$ whenever $|\{v\leq0\}|>0$. For $\beta<1$, the choice $G=0$ makes $G$ continuous at $0$, and $\J$ is continuous on $\H$. For $1\leq\beta<3$, $\J$ is not finite at every point of $\H$. We therefore regard $\J$ as an extended real-valued functional on $\H$.

We associate $\J$ formally with the problem
\begin{equation} {\label{sd6}}
\begin{cases}
-\Delta v=\lambda \left(\chi_{\{g(v)<b\}} g(v)^{-\de}+ \left(\int_\Omega \frac{H(y,g(v))}{\left|x-y\right|^\mu} dy\right)h(x,g(v))\right) g^\p(v) & \text{ in } \Omega,\\v>0  \text{ in } \Omega, \\v=0  \text{ on } \partial\Omega.
\end{cases}
\tag{$\mathcal{Q}_{\lambda}^b$}
\end{equation}
\begin{definition}\label{defsd2}
A function $v\in \H$ is said to be a weak solution to the problem \eqref{sd6} if $\operatorname*{ess\,inf}_K v>0$ for every compact set $K \subset \Omega$ and, for every $w \in C^\infty_0(\Omega)$,
\begin{equation}
{\label{sd7}}\int_\Omega \na v\cdot\na w dx = \lambda \int_\Omega \chi_{\{g(v)<b\}} g(v)^{-\beta} g^\p(v) w dx +\lambda \int_\Omega \int_\Omega \frac{H(y,g(v)(y)) h(x,g(v)(x))}{\left|x-y\right|^\mu} g^\p(v) w(x) dxdy.
\end{equation}
\end{definition}
\begin{rem}{\label{remsd1}}
Let $v$ and $w$ be in $\H$ such that $v \geq c \varphi_\beta$, for some $c>0.$ Let $1<\beta<3$. Since $g$ is increasing on $(0,\infty)$, Lemma~\ref{lemsd1}$(g_6)$ gives
\begin{equation*}
g(v)\geq g(c \phi^{\frac{2}{\beta+1}}_1) \geq g(c_1 d(x)^{\frac{2}{\beta+1}}) \geq
\begin{cases}
g(1) c_1 d(x)^{\frac{2}{\beta+1}}, & \text{ if } c_1 d(x)^{\frac{2}{\beta+1}}< 1, \\ g(1) c_1^{\frac{1}{2\alpha}} d(x)^{\frac{1}{\alpha(\beta+1)}}, & \text{ if } c_1 d(x)^{\frac{2}{\beta+1}} \geq 1,
\end{cases}
\end{equation*}
where $c_1$ is a positive constant. By H\"older's inequality, Hardy's inequality and $(g_2)$, we deduce
\begin{align*}
\left|\int_\Omega \chi_{\{g(v)<b\}} g(v)^{-\beta} g^\p(v) w dx\right| & \leq \int_{\Omega \cap \{x: c_1 d(x)^{\frac{2}{\beta+1}}<1\}} g(v)^{-\beta} \left|w\right| dx + \int_{\Omega \cap \{x: c_1 d(x)^{\frac{2}{\beta+1}} \geq 1\}} g(v)^{-\beta} \left|w\right| dx \\& \leq c_3 \int_\Omega \frac{1}{d(x)^{\frac{2\beta}{\beta+1}}} \left|w\right| dx + c_4 \int_\Omega \frac{1}{d(x)^{\frac{\beta}{\alpha(\beta+1)}}} \left|w\right| dx \\& \leq c_3 \int_\Omega \frac{1}{d(x)^{\frac{\beta-1}{\beta+1}}} \frac{\left|w\right|}{d(x)}dx+ c_5 \int_\Omega \left|w\right|dx \\& \leq c_3 \left(\int_\Omega \frac{1}{d(x)^{\frac{2(\beta-1)}{\beta+1}}}dx\right)^\frac{1}{2} \left(\int_\Omega \frac{\left|w\right|^2}{d(x)^2} dx\right)^\frac{1}{2}+ c_5 \|\nabla w\|_2 \leq c_6 \|\nabla w\|_2.
\end{align*}
In the last inequality, we use $\displaystyle \int_\Omega d(x)^{-\frac{2(\beta-1)}{\beta+1}}dx <\infty$ for $1<\beta<3$, and Hardy's inequality $\displaystyle \int_\Omega \frac{\left|w\right|^2}{d(x)^2} dx \leq C \|\nabla w\|_2^2$ \cite[(0.5)]{brezis1997}. For the case $0<\beta \leq 1$, the same computation with $\varphi_\beta=\phi_1$ when $\beta<1$, and with $\varphi_1=\phi_1(-\log\phi_1)^{1/2}$ when $\beta=1$, gives the same bound. For the Choquard term in \eqref{sd7}, the bound $C\|h(\cdot,g(v))\|_2\|w\|_2$ follows from \eqref{sd85} below, since $g^\p\leq1$. Hence, the right-hand side of \eqref{sd7} is finite for every $w\in\H$, so Definition~\ref{defsd2} is meaningful.
\end{rem}

Throughout the paper, we write
\begin{align*}
\rho(t)&:=\chi_{\{g(t)<b\}}\,g(t)^{-\beta}g^\p(t),\quad \Phi(t):=G(g(t)),\quad \Psi(x,t):=H(x,g(t)), \\&
\mathcal{K}[v](x):=\int_\Omega\frac{H(y,g(v(y)))}{|x-y|^\mu}dy ,
\end{align*}
so that $\Phi^\p=\rho$ on $(0,\infty)\setminus\{g^{-1}(b)\}$ and $\Psi_t(x,t)=h(x,g(t))g^\p(t)$. By $(g_9)$ and the monotonicity of $t\mapsto\chi_{\{g(t)<b\}}$, $\rho$ is nonincreasing on $(0,\infty)$, hence $\Phi$ is concave there, and on $[0,\infty)$ when $\beta<1$. Moreover, writing
\[
\Psi_t(x,t)=\frac{h(x,g(t))}{g(t)^{4\alpha-1}}\cdot \frac{g(t)^{4\alpha-1}g^\p(t)}{t}\cdot t ,
\]
we see from Lemma~\ref{lemsd1}$(g_8)$ and $(f_4)$ that $\Psi_t(x,\cdot)$ is increasing on $(0,\infty)$. We write $t_b:=g^{-1}(b)$, so that $\rho(t)=0$ for $t\geq t_b$, and $\rho$ is continuous on $(0,\infty)\setminus\{t_b\}$. Finally, we set
\[
C_\mu:=\sup_{x\in\Omega}\int_\Omega|x-y|^{-\mu}dy<\infty .
\]
For every $v\in H^1_0(\Omega)$, $(g_3)$, $(f_2)$ and $(f_3)$ give $h(x,g(v))+H(x,g(v))\leq C_\epsilon e^{2\alpha(1+\epsilon)v^2}$, and Trudinger's inequality \cite{trudinger1967} gives $e^{qv^2}\in L^1(\Omega)$ for every $q<\infty$. Together with H\"older's inequality, applied to $\mathcal{K}[v]$ with an exponent $r>2/(2-\mu)$, this yields
\begin{equation}\label{sd85}
h(\cdot,g(v)),\ H(\cdot,g(v))\in L^q(\Omega)\ \text{ for every }q<\infty,\quad \mathcal{K}[v]\in L^\infty(\Omega)\qquad\bigl(v\in H^1_0(\Omega)\bigr).
\end{equation}

We now make precise the sense in which $u=g(v)$ solves the original problem.
\begin{definition}\label{defsd3}
A function $u\in H^1_0(\Omega)\cap L^\infty_{loc}(\Omega)$ with $u>0$ a.e.\ in $\Omega$, $\operatorname*{ess\,inf}_K u>0$ for every compact $K\subset\Omega$, and $u^{2\alpha}\in H^1_0(\Omega)$ is a weak solution of \eqref{sd1} if, for every $\phi\in H^1_0(\Omega)\cap L^\infty(\Omega)$ with compact support in $\Omega$,
\begin{align}\label{sd8}
&\int_\Omega\bigl(1+2\alpha u^{4\alpha-2}\bigr)\nabla u\cdot\nabla\phi\,dx+2\alpha(2\alpha-1)\int_\Omega u^{4\alpha-3}|\nabla u|^2\phi\,dx\notag\\
&\qquad=\lambda\int_\Omega\Bigl(\chi_{\{u<b\}}u^{-\beta}+\Bigl(\int_\Omega\frac{H(y,u(y))}{|x-y|^\mu}dy\Bigr)h(x,u)\Bigr)\phi\,dx .
\end{align}
\end{definition}

\begin{lemma}\label{lemsd2}
Let $v\in H^1_0(\Omega)$ be a weak solution of \eqref{sd6} in the sense of Definition~\ref{defsd2}. Then $v\in W^{2,p}_{loc}(\Omega)$ for every $p<\infty$, and $u:=g(v)$ is a weak solution of \eqref{sd1} in the sense of Definition~\ref{defsd3}. Conversely, if $u$ is a weak solution of \eqref{sd1}, then $v:=g^{-1}(u)$ is a weak solution of \eqref{sd6}.
\end{lemma}
The proof is a direct computation, and we give it in Appendix~\ref{secsd3}.

\subsection{Main results}
Set
\begin{align}\label{sd9}
\lambda^*_b:=\sup\{\lambda>0:\ \eqref{sd6}\text{ has a weak solution}\}.
\end{align}

Let $J_\lambda$ be the functional in \eqref{sd4}. Our main result is the following.

\begin{theorem}\label{thmsd1}
Let $\Omega\subset\R^2$ be a bounded domain with smooth boundary, let $b>0$, $0<\beta<3$, $0<\mu<2$, $\alpha>\frac12$, and let $f$ satisfy $(f_1)$--$(f_4)$. Then:
\begin{enumerate}[label=\textnormal{(\roman*)}]
\item $0<\lambda^*_b<\infty$;
\item \eqref{sd1} has a weak solution if and only if $\lambda\in(0,\lambda^*_b]$. For every $\lambda\in(0,\lambda^*_b]$, there is a weak solution $u_\lambda=g(v_\lambda)$, where $v_\lambda\in H^1_0(\Omega)\cap W^{2,p}_{loc}(\Omega)$ for every $p<\infty$ and $v_\lambda\geq c_\lambda\varphi_\beta$ in $\Omega$ for some $c_\lambda>0$. If $\lambda<\lambda^*_b$, then $v_\lambda$ is a local minimizer of $J_\lambda$ in the $H^1_0(\Omega)$ topology;
\item if in addition $(f_5)$ holds, then for every $\lambda\in(0,\lambda^*_b)$, \eqref{sd1} has a second weak solution $u=g(v)$ with $v\geq v_\lambda$ and $v\not\equiv v_\lambda$.
\end{enumerate}
\end{theorem}

The lower bound in (ii) is the boundary behavior found in \cite{lazer1991}. The solution at $\lambda=\lambda^*_b$ is obtained as a limit of solutions for $\lambda<\lambda^*_b$. We do not know whether it is a local minimizer of $J_{\lambda^*_b}$, and part (iii) remains open for $\lambda=\lambda^*_b$.

In the original variable, $(f_5)$ holds if $\inf_xf(x,u)\geq c_fe^{c\,\omega(u)}$ for large $u$, with $c_f>0$ and $c>c_\alpha$, and it fails if $\sup_xf(x,u)\leq c_fe^{c\,\omega(u)}$ for large $u$, with $c<c_\alpha$ (Remark~\ref{remsd7}(b)). Here $\omega$ and $c_\alpha$ are given in Lemma~\ref{lemsd22}. In particular, part (iii) does not apply to $h(x,u)=u^{4\alpha-1}e^{u^{4\alpha}}$, nor to any $h$ whose factor $f$ has polynomial growth. For $\alpha\leq1$, the existence of a second solution for such $h$ remains open (Remark~\ref{remsd6}(b)). For $\alpha>1$, Proposition~\ref{propsd2} gives a second solution, in particular for $h(x,u)=u^{4\alpha-1}e^{u^{4\alpha}}$, provided $\|v_\lambda\|_{L^\infty(\Omega)}>\delta_\infty/(2\alpha)^{1/2}$, with $\delta_\infty$ as in \eqref{sd83}.

\subsection{Difficulties and strategy}\label{secsd2}
After the change of variables, each of the other three features of \eqref{sd1} removes a standard tool. The singular term makes the energy of \eqref{sd6} non-differentiable, and for $\beta\geq1$ the solutions are not in $C^1(\bar\Omega)$ \cite{lazer1991}. Existence, therefore, has to be built from sub- and supersolutions. Their boundary behavior, combined with Hardy's inequality, keeps them in $H^1_0(\Omega)$ for $\beta<3$. The jump across $u=b$ destroys differentiability a second time, at every point where $u=b$. We regularize the characteristic function, solve the regularized problems, and pass to the limit. We then show that the level set $\{g(v)=b\}$ of the limit has measure zero.

The Choquard term is nonlocal, so pointwise comparison arguments for local nonlinearities now require a global ordering of the functions being compared. Moreover, the translation of the problem to the first solution, which for a local nonlinearity produces a $C^1$ energy on $H^1_0(\Omega)$, here produces the sum of a $C^1$ functional and a convex continuous one. We work on the positive cone with Szulkin's notion of critical point \cite{szulkin1986}, following \cite{anthal2025}.

Finally, the critical exponential growth causes a loss of compactness. Palais--Smale sequences can be handled only below the level
\begin{align*}
c^*=\frac{(4-\mu)\pi}{4\alpha},
\end{align*}
which combines the Trudinger--Moser threshold with the exponent $\frac{4}{4-\mu}$ of the Hardy--Littlewood--Sobolev inequality; see \cite{alves2016} for the semilinear counterpart. The level estimate below $c^*$ is obtained with Moser functions in the variable $v$ (Lemma~\ref{lemsd18}). On the ball where the Moser function concentrates, the first solution satisfies $v_\lambda\geq\nu>0$. The translation by $v_\lambda$ therefore produces an extra factor $e^{8\alpha\nu tm_n}$ in the estimate, where $m_n$ is the height of the Moser function. As in \cite{dhanya2014}, this factor allows a growth condition without any threshold, namely $(f_5)$. In the original variable, $(f_5)$ requires $f$ to compensate the factor $e^{-\Xi(v)}$ of Lemma~\ref{lemsd22}, which measures the gap between the scales $\exp(|u|^{4\alpha})$ and $\exp(2\alpha v^2)$. This is why $(f_5)$ excludes $f$ of polynomial growth.

In summary, we use sub- and supersolutions and the variational form of Perron's method (Lemma~\ref{lemsd20}) for the regularized problems, and a comparison argument of Brezis--Nirenberg type \cite{brezis1993,hirano2008} for the first solution. For the second solution, we use a mountain-pass argument on the positive cone with Szulkin's notion of critical point \cite{szulkin1986}. The Palais--Smale sequence at the mountain-pass level is constructed directly (Lemma~\ref{lemsd19}), because our Palais--Smale sequences are only known to have weak limits.

The paper is organized as follows. Section~2 treats the purely singular discontinuous problem and provides the subsolution $z_\lambda$, the boundary behavior and the comparison principle. Section~3 proves parts (i) and (ii) of Theorem~\ref{thmsd1}: Perron's method, the threshold, the regularized problems, the passage to the limit, the local minimality, and the existence at $\lambda=\lambda^*_b$. Section~4 proves part (iii) and Proposition~\ref{propsd2}, which covers $h(x,u)=u^{4\alpha-1}e^{u^{4\alpha}}$ for $\alpha>1$. Appendix~\ref{secsd3} contains the proofs of Lemmas~\ref{lemsd1}, \ref{lemsd2}, \ref{lemsd3}(c), \ref{lemsd20} and~\ref{lemsd19}, which are elementary or follow standard schemes.
	
	\section{The purely singular discontinuous problem}
In this section, we study the purely singular problem obtained by removing the Choquard term from \eqref{sd6}. Lemma~\ref{lemsd3} proves a comparison principle and a lower bound near $\partial\Omega$ for singular problems with a nonincreasing nonlinearity, which are used throughout the paper. At the end of this section, we show that the purely singular problem has a unique solution $z_\lambda$.

For fixed $b>0$ and $0<\beta<3$, we consider the singular discontinuous problem
\begin{align} {\label{sd10}}
\begin{cases}
-\Delta v=\lambda \chi_{\{g(v)<b\}} g(v)^{-\de} g^\p(v) & \text{ in } \Omega,\\v>0  \text{ in } \Omega, \\v=0  \text{ on } \partial\Omega,
\end{cases}
\tag{$\mathcal{R}_{\lambda}^b$}
\end{align}
and the approximating problem
\begin{align} {\label{sd11}}
\begin{cases}
-\Delta v=\lambda \chi_{\eta}(g(v)-b) g(v)^{-\de}g^\p(v) & \text{ in } \Omega,\\v>0  \text{ in } \Omega, \\v=0  \text{ on } \partial\Omega,
\end{cases}
\tag{$\mathcal{R}_{\lambda}^{b,\eta}$}
\end{align}
where
\begin{align*}
\chi_\eta(s)=\chi_{(-\infty, -\eta)}(s)-s \eta^{-1} \chi_{[-\eta,0)}(s), \quad s\in \mathbb{R}.
\end{align*}

\begin{rem}{\label{remsd2}}
If $t_1<t_2$ and $\eta < \eta_0$, then
\begin{align*}
\chi_{\eta_0}(g(t_2)-b) \leq \chi_{\eta}(g(t_1)-b).
\end{align*}
\end{rem}

We first collect what we need about singular problems with a nonincreasing nonlinearity. The nonlinearities in this paper, such as $g(t)^{-\beta}g^\p(t)$, are not pure powers, so the results of \cite{crandall1977dirichlet,lazer1991} do not apply directly. Moreover, we need the boundary estimate for every weak solution, and not only for the solution obtained by approximation. The result of \cite{crandall1977dirichlet} is used only once, for the supersolution $\zeta^{(1)}$ in the proof of part~(c), given in Appendix~\ref{secsd3}.

\begin{lemma}\label{lemsd3}
Let $0<\beta<3$, $\lambda>0$, and let $\rho:(0,\infty)\to[0,\infty)$ be nonincreasing and satisfy, for constants $C_\rho,c_\rho,t_0>0$,
\begin{align}\label{sd12}
\rho(t)\leq C_\rho\,(t^{-\beta}+1)\ \text{ for all } t>0,\qquad \rho(t)\geq c_\rho\,t^{-\beta}\ \text{ for } 0<t\leq t_0 .
\end{align}
Consider the problem
\begin{align}\label{sd13}
-\Delta z=\lambda\rho(z)\ \text{ in }\Omega,\qquad z>0\ \text{ in }\Omega,\qquad z=0\ \text{ on }\partial\Omega,
\end{align}
understood as in Definition~\ref{defsd2}: $z\in H^1_0(\Omega)$, $\operatorname*{ess\,inf}_K z>0$ for every compact $K\subset\Omega$, and
\begin{align}{\label{sd14}}
\int_\Omega\nabla z\cdot\nabla\phi\,dx=\lambda\int_\Omega\rho(z)\phi\,dx \quad\text{for all }\phi\in C^\infty_0(\Omega).
\end{align}
\begin{enumerate}[label=\textnormal{(\alph*)}]
\item There is $\kappa>0$, depending only on $\lambda,c_\rho,t_0,\beta$ and $\Omega$, such that every weak solution $z$ of \eqref{sd13} satisfies $z\geq\kappa\varphi_\beta$ in $\Omega$.
\item \eqref{sd13} has at most one weak solution. More generally, if $\rho_1\leq\rho_2$ both satisfy \eqref{sd12} and $z_i$ is a weak solution of \eqref{sd13} with $\rho=\rho_i$ ($i=1,2$), then $z_1\leq z_2$. The same conclusion holds if $z_1$ is only a weak subsolution, or $z_2$ only a weak supersolution, provided this sub- or supersolution is bounded below by a positive multiple of $\varphi_\beta$.
\item If $\rho$ is continuous, \eqref{sd13} has a weak solution $z$, and $z\in C(\bar\Omega)\cap W^{2,p}_{loc}(\Omega)$ for every $p<\infty$, with $\kappa\varphi_\beta\leq z\leq\zeta_\lambda$, where $\zeta_\lambda\in H^1_0(\Omega)\cap C(\bar\Omega)$ depends only on $\lambda,C_\rho,\beta$ and $\Omega$.
\end{enumerate}
\end{lemma}
\begin{proof}
Put $\gamma:=\frac{2}{\beta+1}$ and $L:=-\log\phi_1$. As $0<\phi_1<1$, we have $L\geq L_0:=-\log\|\phi_1\|_{L^\infty(\Omega)}>0$. From $-\Delta\phi_1=\lambda_1\phi_1$, where $\lambda_1$ is the first eigenvalue of $-\Delta$ in $H^1_0(\Omega)$, a direct computation gives
\begin{align*}
-\Delta\phi_1^{\gamma}&=\gamma\lambda_1\phi_1^{\gamma}+\gamma(1-\gamma)|\nabla\phi_1|^2\phi_1^{-\gamma\beta}\qquad(\beta>1),\\
-\Delta\bigl(\phi_1L^{1/2}\bigr)&=\lambda_1\phi_1\bigl(L^{1/2}-\tfrac12L^{-1/2}\bigr)+\phi_1^{-1}|\nabla\phi_1|^2\bigl(\tfrac12L^{-1/2}+\tfrac14L^{-3/2}\bigr)\qquad(\beta=1),\\
-\Delta\phi_1&=\lambda_1\phi_1\qquad(\beta<1).
\end{align*}
We claim that there is $A=A(\Omega,\beta)>0$ such that
\begin{align}\label{sd81}
-\Delta\varphi_\beta\leq A\,\varphi_\beta^{-\beta}\qquad\text{in }\Omega .
\end{align}
We consider three cases. If $\beta<1$, then $\varphi_\beta=\phi_1$, and $\phi_1\leq\phi_1^{-\beta}$ because $\phi_1<1$.
If $\beta>1$, then $\varphi_\beta^{-\beta}=\phi_1^{-\gamma\beta}$. The second term is a bounded multiple of $\phi_1^{-\gamma\beta}$, as $0<\gamma<1$ and $|\nabla\phi_1|$ is bounded. For the first term, we use $\phi_1^{\gamma}\leq\phi_1^{-\gamma\beta}$.
Finally, if $\beta=1$, then $\varphi_1^{-1}=\phi_1^{-1}L^{-1/2}$. The term $-\frac12\lambda_1\phi_1L^{-1/2}$ is negative, and we drop it. The function $s\mapsto s^2\log(1/s)$ is bounded on $(0,1)$, so
\begin{align*}
\lambda_1\phi_1L^{1/2}=\lambda_1(\phi_1^2L)\,\phi_1^{-1}L^{-1/2}\leq C\varphi_1^{-1}.
\end{align*}
As $L\geq L_0$, we have $\frac12L^{-1/2}+\frac14L^{-3/2}\leq(\frac12+\frac14L_0^{-1})L^{-1/2}$, so the last term is at most $\|\nabla\phi_1\|_{L^\infty(\Omega)}^2(\frac12+\frac14L_0^{-1})\varphi_1^{-1}$. Thus \eqref{sd81} holds in all three cases.

Choose $\kappa>0$ so small that
\begin{align*}
\kappa^{1+\beta}A\leq\lambda c_\rho,\qquad \kappa\|\varphi_\beta\|_{L^\infty(\Omega)}\leq t_0 .
\end{align*}
Then, by \eqref{sd81} and \eqref{sd12},
\begin{align}\label{sd15}
-\Delta(\kappa\varphi_\beta)\leq\kappa^{1+\beta}A\,(\kappa\varphi_\beta)^{-\beta}\leq\lambda c_\rho(\kappa\varphi_\beta)^{-\beta}\leq\lambda\rho(\kappa\varphi_\beta)\quad\text{in }\Omega .
\end{align}

(a) Let $z$ be a weak solution and $\epsilon>0$. The function $\psi:=(\kappa\varphi_\beta-z-\epsilon)^+$ belongs to $H^1_0(\Omega)$ and is supported in $K_\epsilon:=\{\kappa\varphi_\beta\geq\epsilon\}$. This is a compact subset of $\Omega$, and $\rho(z)\leq\rho(\operatorname*{ess\,inf}_{K_\epsilon}z)$ is bounded there. Mollifying $\psi$, we see that the weak identity for $z$ holds with $\phi=\psi$. Integrating \eqref{sd15} against $\psi\geq0$ and subtracting the weak identity for $z$, we obtain
\begin{align*}
\int_\Omega|\nabla\psi|^2dx=\int_\Omega\nabla(\kappa\varphi_\beta-z)\cdot\nabla\psi\,dx\leq\lambda\int_\Omega\bigl(\rho(\kappa\varphi_\beta)-\rho(z)\bigr)\psi\,dx\leq0 ,
\end{align*}
because $\kappa\varphi_\beta>z$ on $\{\psi>0\}$ and $\rho$ is nonincreasing. Hence $\psi=0$, that is, $z\geq\kappa\varphi_\beta-\epsilon$. Letting $\epsilon\to0$, we obtain (a).

(b) Let $z_1,z_2$ be as stated. By (a), the monotonicity of $\rho_i$ and \eqref{sd12}, $\rho_i(z_i)\leq C_\rho((\kappa\varphi_\beta)^{-\beta}+1)$. The proof of the estimate in Remark~\ref{remsd1} applies with $g(v)$ replaced by $\kappa\varphi_\beta$ and shows that
\begin{align*}
\int_\Omega\rho_i(z_i)|w|dx\leq C\|\nabla w\|_2\quad\text{for all }w\in H^1_0(\Omega).
\end{align*}
Hence both sides of \eqref{sd14} are continuous on $H^1_0(\Omega)$, so \eqref{sd14} holds for every $\phi\in H^1_0(\Omega)$. Testing the difference of the two identities with $(z_1-z_2)^+$, we obtain
\begin{align*}
\int_\Omega|\nabla(z_1-z_2)^+|^2dx&\leq\lambda\int_\Omega\bigl(\rho_1(z_1)-\rho_2(z_2)\bigr)(z_1-z_2)^+dx\\&\leq\lambda\int_\Omega\bigl(\rho_2(z_1)-\rho_2(z_2)\bigr)(z_1-z_2)^+dx\leq0 ,
\end{align*}
so $z_1\leq z_2$. Uniqueness follows by taking $\rho_1=\rho_2$.

(c) The proof is standard and is given in Appendix~\ref{secsd3}. \qed
\end{proof}

\begin{lemma}{\label{lemsd6}}
For any $\lambda>0$, the problem \eqref{sd10} admits a unique weak solution $z_\lambda$ in $H_0^1(\Omega)$. In addition, $d_\lambda \varphi_\beta \leq z_\lambda$ in $\Omega$ for some $d_\lambda>0$, and $z_\lambda\leq z_{\lambda^\p}$ for $\lambda\leq\lambda^\p$.
\end{lemma}

\begin{proof}
We use the notation
\begin{align*}
\rho_\eta(t):=\chi_\eta(g(t)-b)\,g(t)^{-\beta}g^\p(t),\quad \hat\rho(t):=g(t)^{-\beta}g^\p(t),\quad t>0 .
\end{align*}
The map $t\mapsto\chi_\eta(g(t)-b)$ is nonincreasing and nonnegative, and $\hat\rho$ is positive and nonincreasing by Lemma~\ref{lemsd1}$(g_9)$. Hence $\rho_\eta$ is continuous, nonnegative, and nonincreasing. By $(g_3)$ and $(g_6)$, it satisfies \eqref{sd12} with constants independent of $\eta\in(0,b/2)$. Lemma~\ref{lemsd3}(b),(c) provides a unique solution $z_\eta\in H^1_0(\Omega)\cap C(\bar\Omega)$ of \eqref{sd11}, and $z_\eta\geq\kappa\varphi_\beta$ with $\kappa$ independent of $\eta$.

By Remark~\ref{remsd2}, $\rho_{\eta_0}\leq\rho_\eta$ for $\eta<\eta_0$, so Lemma~\ref{lemsd3}(b) yields $z_\eta\geq z_{\eta_0}$ whenever $\eta<\eta_0$. Similarly, Lemma~\ref{lemsd3}(b),(c), applied to $\hat\rho$, provides a unique solution $w_\lambda\in H^1_0(\Omega)\cap C(\bar\Omega)$ of
\begin{align*}
-\Delta w=\lambda\hat\rho(w)\ \text{ in }\Omega,\quad w>0\ \text{ in }\Omega,\quad w=0\ \text{ on }\partial\Omega,
\end{align*}
and Lemma~\ref{lemsd3}(b) gives $z_\eta\leq w_\lambda$, because $\rho_\eta\leq\hat\rho$. As $z_\eta\geq\kappa\varphi_\beta$ and $\rho_\eta(t)\leq g(t)^{-\beta}$, the estimate in Remark~\ref{remsd1} applies and shows that
\begin{align*}
\int_\Omega|\nabla z_\eta|^2dx=\lambda\int_\Omega\rho_\eta(z_\eta)z_\eta\,dx\leq C\|\nabla z_\eta\|_2 ,
\end{align*}
with $C$ independent of $\eta$. Hence $\{z_\eta\}$ is bounded in $H^1_0(\Omega)$. Up to a subsequence, $z_\eta\rightharpoonup z_\lambda$ in $H^1_0(\Omega)$ and $z_\eta\to z_\lambda$ a.e.\ in $\Omega$, with $\kappa\varphi_\beta\leq z_\lambda\leq w_\lambda$.

On each compact $K\subset\Omega$, we have $\rho_\eta(z_\eta)\leq\hat\rho(\kappa\min_K\varphi_\beta)$. Since $z_\eta$ satisfies
\begin{align}{\label{sd16}}
-\Delta z_\eta=\lambda \chi_{\eta}(g(z_\eta)-b) g(z_\eta)^{-\de}g^\p(z_\eta),
\end{align}
elliptic regularity shows that $\{z_\eta\}$ is bounded in $W^{2,p}_{loc}(\Omega)$ for every $p<\infty$, hence in $C^\theta_{loc}(\Omega)$ for some $\theta\in(0,1)$. Up to a subsequence, $z_\eta\to z_\lambda$ uniformly on compact subsets of $\Omega$, and $z_\lambda\in W^{2,p}_{loc}(\Omega)$.

Let $\Omega_b:=\{x\in\Omega:g(z_\lambda(x))=b\}$, which is a closed subset of $\Omega$. On $\Omega\setminus\Omega_b$, we have $\chi_\eta(g(z_\eta)-b)\to\chi_{\{g(z_\lambda)<b\}}$, and hence $\rho_\eta(z_\eta)\to\rho(z_\lambda)$ pointwise. Passing to the limit in \eqref{sd16}, we obtain
\begin{align}{\label{sd17}}
-\Delta z_\lambda=\lambda \chi_{\{g(z_\lambda)<b\}} g(z_\lambda)^{-\de}g^\p(z_\lambda) \text{ a.e.\ in } \Omega \setminus \Omega_b.
\end{align}
On $\Omega_b$, which is a level set of $z_\lambda\in W^{2,p}_{loc}(\Omega)$, we have $\nabla z_\lambda=0$ and $\Delta z_\lambda=0$ a.e., while $\chi_{\{g(z_\lambda)<b\}}=0$. Hence \eqref{sd17} holds a.e.\ in $\Omega$. Since $\rho(z_\lambda)\leq g(\kappa\varphi_\beta)^{-\beta}$, the estimate in Remark~\ref{remsd1} shows that \eqref{sd17} also holds in the weak sense, so $z_\lambda$ is a weak solution of \eqref{sd10}.

By $(g_3)$ and $(g_6)$, the function $\rho$ satisfies \eqref{sd12}. Parts (a) and (b) of Lemma~\ref{lemsd3} give the uniqueness of the weak solution and the bound $z_\lambda\geq\kappa\varphi_\beta$ in $\Omega$. Finally, if $\lambda\leq\lambda^\p$, Lemma~\ref{lemsd3}(b), applied with $1$ in place of $\lambda$ and with $\rho_1=\lambda\rho$, $\rho_2=\lambda^\p\rho$, gives $z_\lambda\leq z_{\lambda^\p}$. \qed
\end{proof}
	\section{Existence of the first solution}
In this section, we prove parts (i) and (ii) of Theorem~\ref{thmsd1}. We begin with some auxiliary lemmas.

\begin{lemma}\label{lemsd7}
Let $\lambda>0$ and let $\bar v$ be a weak solution of \eqref{sd6}. Then $\bar v\geq z_\lambda$ in $\Omega$, and $\bar v\in L^\infty(\Omega)$.
\end{lemma}

\begin{proof}
Since the Choquard term is nonnegative,
\begin{align*}
-\Delta \bar v\geq \lambda\rho(\bar v),\qquad -\Delta z_\lambda=\lambda\rho(z_\lambda).
\end{align*}
For $\epsilon>0$, the function $(z_\lambda-\bar v-\epsilon)^+$ is supported in $\{z_\lambda>\epsilon\}\Subset\Omega$, where every term is bounded, so it is an admissible test function. Since $\rho$ is nonincreasing,
\begin{align*}
\int_\Omega|\nabla(z_\lambda-\bar v-\epsilon)^+|^2dx
\leq \lambda\int_\Omega\bigl(\rho(z_\lambda)-\rho(\bar v)\bigr)(z_\lambda-\bar v-\epsilon)^+dx\leq0,
\end{align*}
hence $\bar v\geq z_\lambda-\epsilon$. Letting $\epsilon\to0$ gives $\bar v\geq z_\lambda$.

Put $F:=\mathcal{K}[\bar v]\,h(\cdot,g(\bar v))g^\p(\bar v)$. By \eqref{sd85} and $(g_2)$, $F\in L^q(\Omega)$ for every $q<\infty$.  Fix $q>1$, and let $\zeta\in W^{2,q}(\Omega)\cap H^1_0(\Omega)$ be the solution of $-\Delta\zeta=\lambda F$. Then $\zeta\geq0$ and $\zeta\in L^\infty(\Omega)$. The function $\bar v-\zeta$ satisfies
\begin{align*}
-\Delta(\bar v-\zeta)=\lambda\chi_{\{g(\bar v)<b\}}g(\bar v)^{-\beta}g^\p(\bar v),
\end{align*}
and the right-hand side vanishes on $\{\bar v\geq t_b\}\supset\operatorname{supp}(\bar v-\zeta-t_b)^+$. Since $\bar v\geq z_\lambda\geq\kappa\varphi_\beta$, the estimate in Remark~\ref{remsd1} shows that this identity holds for test functions in $H^1_0(\Omega)$. Testing with $(\bar v-\zeta-t_b)^+$ gives $\int_\Omega|\nabla(\bar v-\zeta-t_b)^+|^2dx=0$, so $\bar v \leq t_b + \zeta$, and hence $\|\bar v\|_{L^\infty(\Omega)}\leq t_b+\|\zeta\|_{L^\infty(\Omega)}$.\qed
\end{proof}

Next, we state the variational form of Perron's method that we use. The standard version \cite[Chapter I, Theorem 2.4]{struwe1996} assumes that the nonlinearity is bounded between the subsolution and the supersolution. Here, the singular term is unbounded near $\partial\Omega$, and the Hardy-type estimate in Remark~\ref{remsd1} takes the place of this assumption.

\begin{lemma}\label{lemsd20}
Let $\lambda>0$ and let $\rho_*:(0,\infty)\to[0,\infty)$ be continuous and nonincreasing, with $\rho_*(t)\leq C_\rho(t^{-\beta}+1)$ for $t>0$. Put $P_*(t):=\int_1^t\rho_*(s)\,ds$ for $t>0$. Let $\underline v,\bar v\in H^1_0(\Omega)$ satisfy $\underline v\leq\bar v$ in $\Omega$, $\underline v\geq c\varphi_\beta$ for some $c>0$, $\bar v\in L^\infty(\Omega)$, and
\begin{align}\label{sd70}
\begin{split}
\int_\Omega\nabla\underline v\cdot\nabla\phi\,dx&\leq\lambda\int_\Omega\bigl(\rho_*(\underline v)+\mathcal{K}[\underline v]\Psi_t(x,\underline v)\bigr)\phi\,dx,\\
\int_\Omega\nabla\bar v\cdot\nabla\phi\,dx&\geq\lambda\int_\Omega\bigl(\rho_*(\bar v)+\mathcal{K}[\bar v]\Psi_t(x,\bar v)\bigr)\phi\,dx
\end{split}
\end{align}
for all $0\leq\phi\in H^1_0(\Omega)$. Let $\mathcal{N}_*:=\{v\in H^1_0(\Omega):\underline v\leq v\leq\bar v\}$ and
\begin{align*}
\mathcal{E}(v):=\frac12\int_\Omega|\nabla v|^2dx-\lambda\int_\Omega P_*(v)\,dx-\frac\lambda2\int_\Omega\int_\Omega\frac{\Psi(y,v(y))\Psi(x,v(x))}{|x-y|^\mu}dxdy,\qquad v\in\mathcal{N}_* .
\end{align*}
Then $\mathcal{E}$ is finite on $\mathcal{N}_*$ and attains its minimum there. Moreover, every minimizer $u$ of $\mathcal{E}$ on $\mathcal{N}_*$ is a weak solution of
\begin{align}\label{sd71}
-\Delta u=\lambda\rho_*(u)+\lambda\mathcal{K}[u]\Psi_t(x,u)\ \text{ in }\Omega,\qquad u=0\ \text{ on }\partial\Omega,
\end{align}
that is, $\operatorname*{ess\,inf}_Ku>0$ for every compact $K\subset\Omega$, and
\begin{align*}
\int_\Omega\nabla u\cdot\nabla\phi\,dx=\lambda\int_\Omega\bigl(\rho_*(u)+\mathcal{K}[u]\Psi_t(x,u)\bigr)\phi\,dx\qquad\text{for every }\phi\in H^1_0(\Omega).
\end{align*}
\end{lemma}

In the setting of Lemma~\ref{lemsd20}, every $v\in\mathcal{N}_*$ satisfies $c\varphi_\beta\leq v\leq\|\bar v\|_{L^\infty(\Omega)}$. Hence $0\leq\rho_*(v)\leq C_\rho((c\varphi_\beta)^{-\beta}+1)$, and the estimate in Remark~\ref{remsd1}, with $g(v)$ replaced by $c\varphi_\beta$, gives
\begin{align}\label{sd72}
\int_\Omega\rho_*(v)|\phi|\,dx\leq C\|\nabla\phi\|_2\qquad\text{for all }\phi\in H^1_0(\Omega),\ v\in\mathcal{N}_* .
\end{align}
The proof of Lemma~\ref{lemsd20} follows the usual scheme of Perron's method, and we give it in Appendix~\ref{secsd3}.

\begin{lemma}{\label{lemsd8}}
$0< \lambda_b^*< \infty.$
\end{lemma}

\begin{proof}
Let $\lambda>0$ and let $v$ be a weak solution of \eqref{sd6}. By Lemmas~\ref{lemsd7} and~\ref{lemsd3}(a), $v\geq z_\lambda\geq\kappa\varphi_\beta$. Hence, by Remark~\ref{remsd1}, the right-hand side of \eqref{sd7} is continuous on $H^1_0(\Omega)$. For the nonlocal term, we use the bound $\|\mathcal{K}[v]\|_\infty\|h(\cdot,g(v))\|_{2}\|w\|_2$, which is finite by \eqref{sd85}. Therefore \eqref{sd7} holds for every $w\in H^1_0(\Omega)$. Testing \eqref{sd7} with the first eigenfunction $\phi_1$ of $-\Delta$, we obtain
\begin{align}\label{sd18}
\begin{split}
\lambda_1 \int_\Omega v \phi_1 dx&= \lambda \int_\Omega \chi_{\{g(v)<b\}} g(v)^{-\beta} g^\p(v) \phi_1(x) dx \\
& \qquad +\lambda \int_\Omega \int_\Omega \frac{H(y,g(v)(y)) h(x,g(v)(x))}{\left|x-y\right|^\mu} g^\p(v) \phi_1(x) dxdy \\
& \geq \lambda \int_\Omega \left(\chi_{\{g(v)<b\}} g(v)^{-\beta} + \left(\int_\Omega \frac{H(y,g(v)(y))}{d^\mu _\Omega}dy\right) h(x,g(v)(x))\right) g^\p(v)\phi_1(x) dx,
\end{split}
\end{align}
where $d_\Omega$ is the diameter of $\Omega.$ Now, let $\lambda\geq1$. By Lemma~\ref{lemsd7}, $v\geq z_\lambda$, and by Lemma~\ref{lemsd6}, $z_\lambda\geq z_1$. Since $H(y,\cdot)$ is nondecreasing,
\begin{align*}
\int_\Omega \frac{H(y,g(v))}{d^\mu _\Omega}dy\geq d_\Omega^{-\mu}\int_\Omega H(y,g(z_1(y)))dy=:c_0>0 \quad\text{for all }x\in\Omega,
\end{align*}
where $c_0$ does not depend on $\lambda$ or on $v$. Set
\begin{align*}
\Pi(t):=\inf_{x\in\bar\Omega}\bigl(\rho(t)+c_0h(x,g(t))g^\p(t)\bigr)\quad\text{for }t>0.
\end{align*}
Then $\Pi>0$ on $(0,\infty)$, $\Pi(t)\to\infty$ as $t\to0^+$, and $\Pi(t)/t\to\infty$ as $t\to\infty$. For the last claim, note that $(f_2)$ yields $\displaystyle \inf_{x\in \bar{\Omega}} h(x,s)\geq \exp{(\frac{s^{4\alpha}}{2})}$ for sufficiently large $s$, while $(g_4)$ and $(g_6)$ give
\begin{align*}
g(t)^{4\alpha}\geq g(1)^{4\alpha}t^{2},\qquad g^\p(t)\geq g(t)/(2\alpha t)\geq c\,t^{-(2\alpha-1)/(2\alpha)}\quad\text{for }t\geq1.
\end{align*}
Hence $\Pi_0:=\displaystyle\inf_{t>0}\Pi(t)/t>0$. The constant $\Pi_0$ depends only on $b$, $\beta$, $c_0$, $f$ and $\alpha$, and not on $\lambda$ or on $v$.  Using this in \eqref{sd18}, we deduce
\begin{align*}
\lambda_1\int_\Omega v\phi_1dx\geq\lambda\int_\Omega \Pi(v)\phi_1dx\geq\lambda \Pi_0\int_\Omega v\phi_1dx ,
\end{align*}
and, dividing by $\int_\Omega v\phi_1dx>0$, we obtain $\lambda\leq\lambda_1/\Pi_0$. Thus, $\lambda^*_b\leq\max\{1,\lambda_1/\Pi_0\}<\infty$.
Next, we show that $0<\lambda^*_b$. We first consider the following singular problem without the jump discontinuity
\begin{align} {\label{sd19}}
\begin{cases}
-\Delta v=\lambda \left(g(v)^{-\de}+ \left(\int_\Omega \frac{H(y,g(v))}{\left|x-y\right|^\mu} dy\right)h(x,g(v))\right) g^\p(v) & \text{ in } \Omega,\\v>0  \text{ in } \Omega, \\v=0  \text{ on } \partial\Omega,
\end{cases}
\tag{$\mathcal{Q}_{\lambda}$}
\end{align}
and then show that its solution solves \eqref{sd6} when $\lambda$ is small.

Let $w_\lambda$ be the solution of $-\Delta w=\lambda\hat\rho(w)$ from the proof of Lemma~\ref{lemsd6}. Since the Choquard term is nonnegative,
\begin{align*}
-\Delta w_\lambda= \lambda g(w_\lambda)^{-\beta}g^\p(w_\lambda) \leq \lambda \left(g(w_\lambda)^{-\de}+ \left(\int_\Omega \frac{H(y,g(w_\lambda))}{\left|x-y\right|^\mu} dy\right)h(x,g(w_\lambda))\right) g^\p(w_\lambda) ,
\end{align*}
so $w_\lambda$ is a subsolution of \eqref{sd19}.
We claim that $\|w_\lambda\|_{L^\infty(\Omega)}\to0$ as $\lambda\to0^+$. By $(g_6)$,
\[
g(t)^{-\beta}g^\p(t)\leq C_\rho(t^{-\beta}+1) \text{ for } t>0,
\]
with $C_\rho:=g(1)^{-\beta}$. For $\lambda\in(0,1]$, let $\xi_\lambda$ be the solution of
\begin{align*}
-\Delta\xi=C_\rho(\xi^{-\beta}+\lambda^{\beta/(1+\beta)})\ \text{ in }\Omega,\qquad \xi>0\ \text{ in }\Omega,\qquad \xi=0\ \text{ on }\partial\Omega,
\end{align*}
given by Lemma~\ref{lemsd3}(c). Then by Lemma~\ref{lemsd3}(b), $\xi_\lambda\leq\xi_1$. A direct computation shows that $\zeta_\lambda:=\lambda^{1/(1+\beta)}\xi_\lambda$ solves $-\Delta\zeta=\lambda C_\rho(\zeta^{-\beta}+1)$. Since $w_\lambda$ is a subsolution of this problem, Lemma~\ref{lemsd3}(b) gives
\begin{align*}
w_\lambda\leq\zeta_\lambda\leq\lambda^{1/(1+\beta)}\|\xi_1\|_{L^\infty(\Omega)},
\end{align*}
which proves the claim.

Next, let $z$ solve $-\Delta z=1, z>0$ in $\Omega$, $z=0$ on $\partial\Omega$, and let
\begin{align*}
S_0&:=\|\xi_1\|_{L^\infty(\Omega)}+\|z\|_{L^\infty(\Omega)},\qquad
A_0:=\sup\{h(x,s):x\in\bar\Omega,\ 0\leq s\leq g(S_0)\},\\
S_1&:=M(1+A_0)C_\mu A_0,
\end{align*}
where $M$ is the constant in $(f_3)$ and $C_\mu$ is defined in Subsection~\ref{secsd1}. For $0<\lambda\leq\lambda_0:=\min\{1,S_1^{-1}\}$, set $\tilde z_\lambda:=w_\lambda+\lambda S_1 z$. Then $0<\tilde z_\lambda\leq S_0$, so by $(f_3)$, we get
\begin{align*}
\Bigl(\int_\Omega \frac{H(y,g(\tilde{z}_\lambda))}{\left|x-y\right|^\mu} dy\Bigr)h(x,g(\tilde{z}_\lambda))g^\p(\tilde{z}_\lambda)
\leq M(1+A_0)C_\mu A_0=S_1 .
\end{align*}
Since $g^{-\beta}g^\p$ is nonincreasing and $\tilde z_\lambda\geq w_\lambda$,
\begin{align*}
-\Delta\tilde z_\lambda=\lambda g(w_\lambda)^{-\beta}g^\p(w_\lambda)+\lambda S_1
\geq\lambda\Bigl(g(\tilde z_\lambda)^{-\beta}+\Bigl(\int_\Omega \frac{H(y,g(\tilde{z}_\lambda))}{\left|x-y\right|^\mu} dy\Bigr)h(x,g(\tilde{z}_\lambda))\Bigr)g^\p(\tilde z_\lambda).
\end{align*}
Thus $\tilde z_\lambda$ is a supersolution of \eqref{sd19}, and
\begin{align*}
\|\tilde z_\lambda\|_{L^\infty(\Omega)}\leq\lambda^{1/(1+\beta)}\|\xi_1\|_{L^\infty(\Omega)}+\lambda S_1\|z\|_{L^\infty(\Omega)}\to0\quad\text{as }\lambda\to0^+.
\end{align*}

We check \eqref{sd70} with $\rho_*:=\hat\rho$, $\underline v:=w_\lambda$ and $\bar v:=\tilde z_\lambda$. All the integrals below are finite by \eqref{sd72}. By the proof of Lemma~\ref{lemsd3}(b), the identity $\int_\Omega\nabla w_\lambda\cdot\nabla\phi\,dx=\lambda\int_\Omega\hat\rho(w_\lambda)\phi\,dx$ holds for every $\phi\in H^1_0(\Omega)$. Since the Choquard term is nonnegative, this gives the first inequality in \eqref{sd70} for $0\leq\phi\in H^1_0(\Omega)$. Similarly, $\int_\Omega\nabla\tilde z_\lambda\cdot\nabla\phi\,dx=\lambda\int_\Omega\bigl(\hat\rho(w_\lambda)+S_1\bigr)\phi\,dx$ for every $\phi\in H^1_0(\Omega)$, and the pointwise inequality between the right-hand sides, proved above, yields the second inequality in \eqref{sd70} for $0\leq\phi\in H^1_0(\Omega)$.

Moreover, $w_\lambda\geq\kappa\varphi_\beta$ by Lemma~\ref{lemsd3}(a), $\tilde z_\lambda\in L^\infty(\Omega)$, and $w_\lambda\leq\tilde z_\lambda$. Lemma~\ref{lemsd20} therefore provides a solution $\tilde v_\lambda$ of \eqref{sd71}, that is, of \eqref{sd19}, with $w_\lambda\leq\tilde v_\lambda\leq\tilde z_\lambda$.  By $(g_3)$, $g(\tilde v_\lambda)\leq\tilde v_\lambda\leq\|\tilde z_\lambda\|_{L^\infty(\Omega)}<b$ for small $\lambda$, so $\chi_{\{g(\tilde v_\lambda)<b\}}\equiv1$ and $\tilde v_\lambda$ solves \eqref{sd6}. It is a weak solution in the sense of Definition~\ref{defsd2} because $\tilde v_\lambda\geq w_\lambda$, which is positive and continuous in $\Omega$. Hence $\lambda_b^*>0$. \qed
\end{proof}

We now consider the regularized problem associated with \eqref{sd6},
\begin{align} {\label{sd20}}
\begin{cases}
-\Delta v=\lambda \left(\chi_{\eta}(g(v)-b) g(v)^{-\de} +\left(\int_\Omega \frac{H(y,g(v))}{\left|x-y\right|^\mu} dy\right)h(x,g(v))\right) g^\p(v)& \text{ in } \Omega,\\v>0  \text{ in } \Omega, \\v=0  \text{ on } \partial\Omega,
\end{cases}
\tag{$\mathcal{S}_{\lambda}^{b,\eta}$}
\end{align}
where $\chi_\eta$ is as in Section~2. Its energy functional on $H^1_0(\Omega)$ is
\begin{align}{\label{sd21}}
\J^\eta (v)=\frac{1}{2} \int_\Omega \left|\nabla v\right|^2\ dx -\lambda \int_\Omega G_\eta(g(v))dx-\frac{\lambda}{2} \int_\Omega \int_\Omega \frac{H(y,g(v)(y)) H(x,g(v)(x))}{\left|x-y\right|^\mu}dxdy,
\end{align}
where
\begin{align}{\label{sd22}}
G_\eta(u)=
\begin{cases}
-\infty \ (\beta\geq1),\ 0 \ (\beta<1), & \text{ if } u\leq 0, \\ \frac{u^{1-\beta}}{1-\beta}, & \text{ if } 0<u<\frac{b}{2}, \\ \frac{1}{1-\beta}(\frac{b}{2})^{1-\beta}+ \int^u_\frac{b}{2} \chi_\eta (t-b) t^{-\beta}dt, & \text{ if } u \geq \frac{b}{2},
\end{cases}
\end{align}
for $0<\beta<3$, $\beta \neq 1$. For $\beta=1$, the terms $(1-\beta)^{-1}t^{1-\beta}$ are replaced by $\log t$.

Let $\lambda\in(0,\lambda^*_b)$. By the proof of Lemma~\ref{lemsd6}, $\kappa\varphi_\beta\leq z_\eta\leq z_\lambda\leq w_\lambda$ for $\eta\in(0,b/2)$, where $\kappa$ does not depend on $\eta$. We fix $\sigma\in(0,1)$ with $\sigma\|z_\lambda\|_{L^\infty(\Omega)}<b/2$. Then $\sigma z_\eta\leq\sigma z_\lambda<b/2$ for every $\eta\in(0,b/2)$. Moreover, since $\lambda<\lambda^*_b$, by the definition \eqref{sd9} of $\lambda^*_b$ as a supremum, there is $\bar\lambda\in(\lambda,\lambda^*_b]$ such that $(\mathcal{Q}^b_{\bar\lambda})$ has a weak solution. We fix such a $\bar\lambda$ and a solution $\bar v_\lambda$ of $(\mathcal{Q}^b_{\bar\lambda})$. By Lemma~\ref{lemsd7}, $\bar v_\lambda\geq z_{\bar\lambda}\geq\kappa\varphi_\beta$ and $\bar v_\lambda\in L^\infty(\Omega)$.

The number $\sigma$, the parameter $\bar\lambda$ and the function $\bar v_\lambda$ depend on $\lambda$. For each $\lambda$, they are kept fixed in Lemmas~\ref{lemsd9}, \ref{lemsd10} and~\ref{lemsd12}, and in Section~4. Lemmas~\ref{lemsd9} and~\ref{lemsd10} hold for every $\sigma\in(0,1)$ with $\sigma\|z_\lambda\|_{L^\infty(\Omega)}<b/2$, and Lemma~\ref{lemsd21} uses one such $\sigma$ for a sequence $\lambda_n$.

\begin{lemma}{\label{lemsd9}}
Let $\eta\in(0,b/2)$ and $\lambda \in (0, \lambda_b^*)$ be fixed, and let $\sigma$ and $\bar v_\lambda$ be as above. Then the approximating problem \eqref{sd20} has a solution $v_\eta.$
\end{lemma}

\begin{proof}
Since the map $t \mapsto \chi_\eta(g(t)-b)g(t)^{- \beta} g^\p(t) $ is nonincreasing for $t>0$, and $\sigma<1$,
\begin{align*}
-\Delta (\sigma z_\eta) &= \sigma \lambda \chi_\eta(g(z_\eta)-b) g(z_\eta)^{-\beta}g^\p(z_\eta) \\
& < \lambda \chi_\eta(g(\sigma z_\eta)-b)g(\sigma z_\eta)^{- \beta} g^\p(\sigma z_\eta) \\
& < \lambda \chi_\eta(g(\sigma z_\eta)-b)g(\sigma z_\eta)^{- \beta} g^\p(\sigma z_\eta) +\lambda \left(\int_\Omega \frac{H(y,g(\sigma z_\eta))}{\left|x-y\right|^\mu} dy\right)h(x,g(\sigma z_\eta)) g^\p(\sigma z_\eta).
\end{align*}
This shows that $\sigma z_\eta$ is a subsolution of \eqref{sd20}.

Let $\bar\lambda$ and $\bar v_\lambda$ be as fixed before the lemma.
Recall that $\rho_\eta$, defined in the proof of Lemma~\ref{lemsd6}, is nonincreasing on $(0,\infty)$. Since $\chi_\eta(g(t)-b)\leq\chi_{\{g(t)<b\}}$, $\lambda<\bar\lambda$, and the Choquard term is nonnegative, we have
\begin{align*}
-\Delta\bar v_\lambda\geq\lambda\rho_\eta(\bar v_\lambda)+\lambda\mathcal{K}[\bar v_\lambda]\Psi_t(x,\bar v_\lambda),\quad
-\Delta(\sigma z_\eta)\leq\lambda\rho_\eta(\sigma z_\eta).
\end{align*}

The two inequalities give $-\Delta (\sigma z_\eta-\bar v_\lambda) \leq \lambda( \rho_\eta(\sigma z_\eta)- \rho_\eta(\bar v_\lambda))$. We test with $(\sigma z_\eta-\bar v_\lambda)^+\in H^1_0(\Omega)$, which is admissible by the estimate in Remark~\ref{remsd1}, because $\sigma z_\eta\geq\sigma\kappa\varphi_\beta$ and $\bar v_\lambda\geq\kappa\varphi_\beta$. This gives
\begin{align*}
\int_\Omega|\nabla(\sigma z_\eta-\bar v_\lambda)^+|^2dx
\leq\lambda\int_\Omega\bigl(\rho_\eta(\sigma z_\eta)-\rho_\eta(\bar v_\lambda)\bigr)(\sigma z_\eta-\bar v_\lambda)^+dx\leq0 ,
\end{align*}
which shows that $\sigma z_\eta\leq\bar v_\lambda$, and hence
\begin{align*}
\mathcal{N}_\eta:=\{v\in H^1_0(\Omega):\sigma z_\eta\leq v\leq\bar v_\lambda\}\neq\emptyset.
\end{align*}

We check \eqref{sd70} with $\rho_*=\rho_\eta$, $\underline v=\sigma z_\eta$ and $\bar v=\bar v_\lambda$. By the proof of Lemma~\ref{lemsd3}(b),
\[
\int_\Omega\nabla z_\eta\cdot\nabla\phi\,dx=\lambda\int_\Omega\rho_\eta(z_\eta)\phi\,dx,
\]
for every $\phi\in H^1_0(\Omega)$. Multiplying by $\sigma$ and using the pointwise inequality $\sigma\rho_\eta(z_\eta)\leq\rho_\eta(\sigma z_\eta)$, proved above, together with the nonnegativity of the Choquard term, we get the first inequality in \eqref{sd70} for $0\leq\phi\in H^1_0(\Omega)$. For $\bar v_\lambda$, the weak identity of Definition~\ref{defsd2} extends to all $\phi\in H^1_0(\Omega)$ as in the proof of Lemma~\ref{lemsd8}. Since $\bar\lambda>\lambda$ and $\chi_\eta(g(t)-b)\leq\chi_{\{g(t)<b\}}$, this identity yields the second inequality in \eqref{sd70}. The function $\rho_\eta$ is continuous and nonincreasing, and $\rho_\eta(t)\leq g(1)^{-\beta}(t^{-\beta}+1)$ by $(g_6)$. Hence Lemma~\ref{lemsd20} applies.

Since $\frac{d}{dt}G_\eta(g(t))=\rho_\eta(t)$ for $t>0$, the functional $\mathcal{E}$ in Lemma~\ref{lemsd20} and the functional $J^\eta_\lambda$ in \eqref{sd21}, with $G_\eta$ as in \eqref{sd22}, differ on $\mathcal{N}_\eta$ by the constant $\lambda|\Omega|G_\eta(g(1))$. Hence $J^\eta_\lambda$ attains its minimum on $\mathcal{N}_\eta$ at some $v_\eta$, and $v_\eta$ solves \eqref{sd20}.

We also record bounds, uniform in $\eta$, that are used in the proof of Lemma~\ref{lemsd10}. Every $v\in\mathcal{N}_\eta$ satisfies $\sigma\kappa\varphi_\beta\leq v\leq\bar v_\lambda$. The functions $G_\eta$ and $G$ are nondecreasing on $(0,\infty)$ and satisfy $G_\eta\leq G$ there, because $\chi_\eta(t-b)\leq\chi_{\{t<b\}}$. Also, $G_\eta=G$ on $(0,b/2]$. Since $\sigma\kappa\varphi_\beta\leq\sigma z_\lambda<b/2$ and $g(t)\leq t$, we obtain
\begin{align}\label{sd76}
G(g(\sigma\kappa\varphi_\beta))=G_\eta(g(\sigma\kappa\varphi_\beta))\leq G_\eta(g(v))\leq G(g(\bar v_\lambda)),\quad 0\leq H(\cdot,g(v))\leq H(\cdot,g(\bar v_\lambda)).
\end{align}
The outer bounds do not depend on $\eta$, and they are integrable over $\Omega$. Indeed, $G(g(\bar v_\lambda))\leq G(g(\|\bar v_\lambda\|_\infty))$ because $\bar v_\lambda$ is bounded. Moreover, $G(g(\sigma\kappa\varphi_\beta))\in L^1(\Omega)$, as shown in the proof of Lemma~\ref{lemsd20} under the assumption $\beta<3.$ More precisely, this function is bounded when $\beta<1$, while for $\beta\geq1$ it is negative and unbounded near $\partial\Omega$, but remains integrable. Finally, the boundedness of $\bar v_\lambda$ and the assumptions on $h$ give
\[
D(\bar v_\lambda):=\int_\Omega\int_\Omega\frac{H(y,g(\bar v_\lambda))H(x,g(\bar v_\lambda))}{|x-y|^\mu}dxdy<\infty.
\]

Hence, as in \eqref{sd73}, for all $v\in\mathcal{N}_\eta$,
\begin{align}\label{sd77}
J^\eta_\lambda(v) \geq \frac12\|\nabla v\|_2^2-C_{\mathcal N},\quad J^\eta_\lambda(v) \leq \frac12\|\nabla v\|_2^2-\lambda\int_\Omega G(g(\sigma\kappa\varphi_\beta))dx,
\end{align}
where $C_{\mathcal N}:=\lambda\int_\Omega G(g(\bar v_\lambda))dx+\frac\lambda2D(\bar v_\lambda)$, independent of $\eta$ and of $v$. \qed
\end{proof}

Define
\begin{align*}
\mathcal{N}_0=\{v \in H^1_0(\Omega): \sigma z_\lambda \leq v \leq \bar{v}_\lambda\},\qquad
\mathcal{N}=\{v \in H^1_0(\Omega): \sigma\kappa \varphi_\beta \leq v \leq \bar{v}_\lambda\}.
\end{align*}
\begin{lemma}{\label{lemsd10}}
Let $\lambda \in (0, \lambda^*_b)$, and let $\sigma$ and $\bar v_\lambda$ be as fixed before Lemma~\ref{lemsd9}. Then
\begin{enumerate}[label=\textnormal{(\roman*)}]
\item  the sequence of approximating solutions $v_\eta$ of \eqref{sd20} converges weakly to $v_\lambda \in \mathcal{N}_0$. Moreover, $v_\lambda$ solves \eqref{sd6}, $v_\lambda\in W^{2,p}_{loc}(\Omega)$ for every $p<\infty$, and $|\{g(v_\lambda)=b\}|=0$;
\item $J^\eta_\lambda (v) \to J_\lambda (v)$ uniformly in $\mathcal{N}$, as $\eta \to 0^+$;
\item $J_\lambda(v_\lambda)=\min_{v \in \mathcal{N}_0} J_\lambda(v).$
\end{enumerate}
\end{lemma}

\begin{proof}
(i) Since $\sigma z_\lambda\in\mathcal{N}_\eta$ and $v_\eta$ minimizes $J^\eta_\lambda$ on $\mathcal{N}_\eta$, the bounds \eqref{sd77} give
\begin{align*}
\tfrac12\|\nabla v_\eta\|_2^2-C\leq J^\eta_\lambda(v_\eta)\leq J^\eta_\lambda(\sigma z_\lambda)\leq C,
\end{align*}
with $C$ independent of $\eta$. Hence, along a sequence $\eta_k\to0$, $v_{\eta_k}\rightharpoonup v_\lambda$ in $H^1_0(\Omega)$ and a.e.\ in $\Omega$. Since $z_{\eta_k}\uparrow z_\lambda$ by the proof of Lemma~\ref{lemsd6}, we get $\sigma z_\lambda\leq v_\lambda\leq\bar v_\lambda$, i.e.\ $v_\lambda\in\mathcal{N}_0$.

Fix $\eta_0\in(0,b/2)$ and a compact set $K\subset\Omega$. For $\eta<\eta_0$, we have $v_\eta\geq\sigma z_{\eta_0}\geq c_K>0$ on $K$. Let
\begin{align*}
F_\eta:=\lambda\bigl(\chi_\eta(g(v_\eta)-b)g(v_\eta)^{-\beta}+\mathcal{K}[v_\eta]\,h(x,g(v_\eta))\bigr)g^\p(v_\eta)
\end{align*}
be the right-hand side of \eqref{sd20}. Since $\Psi_t(x,\cdot)$ is nondecreasing and $\bar v_\lambda\in L^\infty(\Omega)$ (Lemma~\ref{lemsd7}), there is $C_K$ such that
\begin{align*}
0<\lambda\,d_\Omega^{-\mu}\int_\Omega H(y,g(\sigma z_{\eta_0}))dy\;\inf_{x\in K}\Psi_t(x,c_K)
\leq F_\eta\leq C_K\quad\text{on }K .
\end{align*}
By elliptic regularity, $\{v_\eta\}$ is bounded in $W^{2,p}(K)$ for every $p<\infty$. Hence $v_{\eta_k}\to v_\lambda$ in $C^1(K)$, $-\Delta v_{\eta_k}\rightharpoonup-\Delta v_\lambda$ weakly in $L^p(K)$, and $v_\lambda\in W^{2,p}_{loc}(\Omega)$.

Passing to the limit in the lower bound, $-\Delta v_\lambda\geq c>0$ a.e.\ on $K$. On the other hand, $\Delta v_\lambda=0$ a.e.\ on the level set $\Omega_b:=\{g(v_\lambda)=b\}$. Hence $|\Omega_b|=0$.

On $\Omega\setminus\Omega_b$, the convergence $v_{\eta_k}\to v_\lambda$ gives $\chi_{\eta_k}(g(v_{\eta_k})-b)\to\chi_{\{g(v_\lambda)<b\}}$ pointwise, exactly as in the proof of Lemma~\ref{lemsd6}. Moreover, $\mathcal{K}[v_{\eta_k}]\to\mathcal{K}[v_\lambda]$ pointwise by dominated convergence, since $H(y,g(v_{\eta_k}))\leq H(y,g(\bar v_\lambda))$. Hence
\begin{align*}
F_{\eta_k}\to F:=\lambda\bigl(\rho(v_\lambda)+\mathcal{K}[v_\lambda]\,h(x,g(v_\lambda))\bigr)g^\p(v_\lambda)
\end{align*}
a.e.\ in $\Omega$ and, being uniformly bounded on $K$, in $L^p(K)$. Thus $-\Delta v_\lambda=F$ a.e.\ in $\Omega$, and $v_\lambda$ is a weak solution of \eqref{sd6} in the sense of Definition~\ref{defsd2}. The positivity requirement holds because $v_\lambda\geq\sigma z_\lambda$.

(ii) Let $v\in \mathcal{N}.$ We write
\begin{align}\label{sd23}
J^\eta_\lambda(v)-J_\lambda(v)=\lambda\int_\Omega \left(G(g(v))-G_\eta(g(v))\right)dx= \lambda\sum_{i=1}^3 \int_{A_i^\eta} \left(G(g(v))-G_\eta(g(v))\right)dx,
\end{align}
where $A_1^\eta=\{x\in \Omega:g(v(x))\leq b-\eta\}$, $A_2^\eta= \{x\in \Omega: b-\eta <g(v(x))\leq b\}$, and $A_3^\eta=\{x\in \Omega: g(v(x))>b\}$. For $\eta<\frac{b}{2}$,
\begin{align*}
G_\eta(g(v)) -G_\eta(\tfrac{b}{2})=\int_{\frac{b}{2}}^{g(v)} \chi_\eta(t-b)t^{-\beta}dt,\qquad G(g(v))-G(\tfrac{b}{2})=\int_{\frac{b}{2}}^{g(v)} \chi_{\{t<b\}}t^{-\beta}dt,
\end{align*}
and $\chi_\eta(t-b)=\chi_{(-\infty,-\eta)}(t-b)-\frac{t-b}{\eta} \chi_{[-\eta,0)}(t-b)$. Hence the integral over $A_1^\eta$ in \eqref{sd23} vanishes. On $A_2^\eta$, $0\leq 1-\frac{b-t}{\eta}\leq 1$ and $t^{-\beta} \leq C$ for some $C>0$, so
\begin{align*}
\left|\int_{A_2^\eta} \left(G(g(v))-G_\eta(g(v))\right)dx\right| &\leq \int_{A^\eta_2} \int_{\frac{b}{2}}^{g(v)} \left|\chi_\eta(t-b)-\chi_{\{t<b\}}(t)\right| t^{-\beta} dtdx \\
&= \int_{A^\eta_2} \int_{b-\eta}^{g(v)} \left|1-\frac{b-t}{\eta}\right| t^{-\beta} dtdx \\
& \leq C \int_{A_2^\eta}\left|g(v)-b+\eta\right|dx \leq C\left|\Omega\right|\eta.
\end{align*}
Similarly, $\left|\int_{A_3^\eta} \left(G(g(v))-G_\eta(g(v))\right)dx\right| \leq C\left|\Omega\right|\eta.$ These bounds together with \eqref{sd23} complete the proof of (ii).

(iii) By (i), $v_\eta \rightharpoonup v_\lambda$ in $H^1_0(\Omega)$ and $v_\eta \to v_\lambda$ locally uniformly in $\Omega$. The weak lower semicontinuity of the Dirichlet norm, (ii) and the dominated convergence theorem give
\begin{equation}{\label{sd24}}
J_\lambda(v_\lambda)
\leq
\liminf_{\eta\to 0} J_\lambda^\eta(v_\eta)
=
\liminf_{\eta\to 0}
\inf_{v\in \mathcal{N}_\eta} J_\lambda^\eta(v).
\end{equation}
Let $r>0$. By (ii), there is $\eta_r>0$ such that, for $\eta \leq \eta_r$,
\begin{align*}
J_\lambda(v)-r < J_\lambda^\eta(v) < J_\lambda(v)+r,
\text{ for all } v\in\mathcal{N}.
\end{align*}
Since $\mathcal{N}_0\subset \mathcal{N}_\eta\subset\mathcal{N}$, taking the infimum over $v\in \mathcal{N}_\eta$, we obtain
\begin{align*}
\inf_{v\in\mathcal{N}_\eta} J_\lambda(v)-r
\leq
\inf_{v\in\mathcal{N}_\eta} J_\lambda^\eta(v)
\leq
\inf_{v\in\mathcal{N}_\eta} J_\lambda(v)+r
\leq
\inf_{v\in\mathcal{N}_0} J_\lambda(v)+r.
\end{align*}
Hence
\begin{equation}{\label{sd86}}
\limsup_{\eta\to 0}
\inf_{v\in\mathcal{N}_\eta} J_\lambda^\eta(v)
\leq
\inf_{v\in\mathcal{N}_0} J_\lambda(v)+r.
\end{equation}
Letting $r\to0$ in \eqref{sd86} and combining the result with \eqref{sd24}, we get $J_\lambda(v_\lambda)\leq \displaystyle\inf_{v\in\mathcal{N}_0} J_\lambda(v)$. Since $v_\lambda\in\mathcal{N}_0$, this proves (iii). \qed
\end{proof}

\begin{lemma}{\label{lemsd12}}
Let $\lambda \in(0,\lambda^*_b).$ Then $v_\lambda$ is a local minimizer of $J_\lambda$ in the $H^1_0(\Omega)$ topology.
\end{lemma}

\begin{proof}
Suppose, by contradiction, that $v_\lambda$ is not a local minimizer of $J_\lambda$ in $H^1_0(\Omega)$. Then, there exists a sequence $\{v_n\} \subset H^1_0(\Omega)$ such that $v_n \to v_\lambda$ in $H^1_0(\Omega)$ and $J_\lambda(v_n)<J_\lambda(v_\lambda)<\infty.$ We may assume that $v_n\geq0$, and even $v_n>0$ a.e.\ in $\Omega$ when $\beta\geq1$. Indeed, for $\beta\geq1$, $J_\lambda(v_n)<\infty$ and the definition \eqref{sd5} of $G$ on $(-\infty,0]$ imply that $v_n>0$ a.e.\ in $\Omega$. Now let $\beta<1$. Then $G\geq0$ on $[0,\infty)$ and $G=0$ on $(-\infty,0]$, and $H(x,g(|v_n|))\geq H(x,g(v_n))$ by $(f_1)$. Since $|\nabla|v_n||=|\nabla v_n|$, this gives $J_\lambda(|v_n|)\leq J_\lambda(v_n)$. Moreover, $v_n \to v_\lambda$ in $H^1_0(\Omega)$ and $v_\lambda \geq 0$, so $|v_n|\to|v_\lambda|=v_\lambda$ in $H^1_0(\Omega)$, and we may replace $v_n$ by $|v_n|$.

Let $\sigma$ be as above and define $\underline{v}=\sigma z_\lambda$. Let $\bar\lambda\in(\lambda,\lambda^*_b]$ and $\bar v:=\bar v_\lambda$, a weak solution of $(\mathcal{Q}^b_{\bar\lambda})$, be as fixed before Lemma~\ref{lemsd9}. We shall use the following facts.
\begin{enumerate}[label=\textnormal{(\alph*)}]
\item $\underline v\leq v_\lambda\leq\bar v$, since $v_\lambda\in\mathcal{N}_0$. The three functions belong to $W^{2,p}_{loc}(\Omega)$ for every $p<\infty$ (proof of Lemma~\ref{lemsd6}, Lemma~\ref{lemsd2}). Moreover, $\bar v\in L^\infty(\Omega)$ by Lemma~\ref{lemsd7}, so that $\mathcal{K}[\underline v]\leq\mathcal{K}[\bar v]\leq C$ in $\Omega$. We set $t^*:=\|\bar v\|_{L^\infty(\Omega)}+1$.
\item The function $t\mapsto h(x,g(t))g^\p(t)$ is nondecreasing on $(0,\infty)$ (see Subsection~\ref{secsd1}) and Lipschitz on $(0,t^*]$, uniformly in $x\in\bar\Omega$, with a constant $L$. Indeed, its derivative is $\partial_sh(x,g(t))g^\p(t)^2+h(x,g(t))g^{\p\p}(t)$, where $g^{\p\p}=-\alpha(4\alpha-2)g^{4\alpha-3}(g^\p)^4$. The first term is bounded on $\bar\Omega\times(0,t^*]$ because $f\in C^1$ by $(f_1)$, and the second because $(f_4)$ gives $h(x,s)\leq h(x,1)s^{4\alpha-1}$ for $0<s\leq1$, so that $h(x,g)g^{4\alpha-3}\leq Cg^{8\alpha-4}$ is bounded.
\end{enumerate}
Now, we set
$w_n=\max \{\underline{v}, \min \{v_n,\bar{v}\}\}$, $\bar{z}_n=(v_n-\bar{v})^+$ and $\underline{z}_n=(v_n-\underline{v})^-$, and we put
$\overline{R}_n=\{v_n>\bar v\}$ and $\underline{R}_n=\{v_n<\underline v\}$.
Then $w_n \in \mathcal{N}_0$ and $v_n=w_n-\underline{z}_n+ \bar{z}_n$. Moreover,
\begin{align*}
v_n=\bar v+\bar z_n,\ w_n=\bar v\ \ \text{on }\overline R_n,\quad v_n=\underline v-\underline z_n,\ w_n=\underline v\ \ \text{on }\underline R_n,\quad v_n=w_n\ \ \text{elsewhere}.
\end{align*}
Using the symmetry of the kernel and the notation $\Phi$, $\Psi$, $\mathcal{K}$ of Subsection~\ref{secsd1}, we obtain
\begin{align}\label{sd26}
\begin{split}
J_\lambda(v_n)={}&J_\lambda(w_n)+\frac12\|\nabla\bar z_n\|_2^2+\int_{\overline R_n}\nabla\bar v\cdot\nabla\bar z_n\,dx+\frac12\|\nabla\underline z_n\|_2^2-\int_{\underline R_n}\nabla\underline v\cdot\nabla\underline z_n\,dx\\
&-\lambda\int_\Omega\bigl(\Phi(v_n)-\Phi(w_n)\bigr)dx-\frac\lambda2\int_\Omega\bigl(\Psi(x,v_n)-\Psi(x,w_n)\bigr)\bigl(\mathcal{K}[v_n]+\mathcal{K}[w_n]\bigr)dx .
\end{split}
\end{align}
We test the inequality satisfied by the supersolution $\bar v$ with $\bar {z}_n$. This is admissible by Remark~\ref{remsd1}, because $\bar z_n\geq0$, $\bar v\geq\sigma\kappa\varphi_\beta$ and $\bar v\in L^\infty(\Omega)$. In the same way, we test the inequality satisfied by the subsolution $\underline v$ with $\underline {z}_n $, and we get
\begin{align*}
\int_{\overline R_n}\nabla\bar v\cdot\nabla\bar z_n\,dx&\geq\lambda\int_{\overline R_n}\bigl(\rho(\bar v)+\mathcal{K}[\bar v]\Psi_t(x,\bar v)\bigr)\bar z_n\,dx,\\
\int_{\underline R_n}\nabla\underline v\cdot\nabla\underline z_n\,dx&\leq\lambda\int_{\underline R_n}\bigl(\rho(\underline v)+\mathcal{K}[\underline v]\Psi_t(x,\underline v)\bigr)\underline z_n\,dx .
\end{align*}
Inserting these two inequalities in \eqref{sd26}, we deduce
\begin{align}{\label{sd27}}
J_\lambda(v_n) \geq J_\lambda(w_n)+\overline C_n + \underline C_n,
\end{align}
where
\begin{align*}
\overline C_n:={}&\frac12\|\nabla\bar z_n\|_2^2+\lambda\int_{\overline R_n}\bigl[\rho(\bar v)\bar z_n-\Phi(\bar v+\bar z_n)+\Phi(\bar v)\bigr]dx\\
&+\frac\lambda2\int_{\overline R_n}\Bigl[2\mathcal{K}[\bar v]\Psi_t(x,\bar v)\bar z_n-\bigl(\Psi(x,\bar v+\bar z_n)-\Psi(x,\bar v)\bigr)\bigl(\mathcal{K}[v_n]+\mathcal{K}[w_n]\bigr)\Bigr]dx,\\
\underline C_n:={}&\frac12\|\nabla\underline z_n\|_2^2-\lambda\int_{\underline R_n}\bigl[\rho(\underline v)\underline z_n+\Phi(\underline v-\underline z_n)-\Phi(\underline v)\bigr]dx\\
&-\frac\lambda2\int_{\underline R_n}\Bigl[2\mathcal{K}[\underline v]\Psi_t(x,\underline v)\underline z_n-\bigl(\Psi(x,\underline v)-\Psi(x,\underline v-\underline z_n)\bigr)\bigl(\mathcal{K}[v_n]+\mathcal{K}[w_n]\bigr)\Bigr]dx .
\end{align*}

By Lemma~\ref{lemsd10}(iii) and $w_n\in\mathcal{N}_0$, we have $J_\lambda(w_n)\geq J_\lambda(v_\lambda)$. Therefore, by \eqref{sd27}, we get
\begin{align}\label{sd28}
J_\lambda(v_\lambda)>J_\lambda(v_n)\geq J_\lambda(v_\lambda)+\overline C_n+\underline C_n \quad\text{for all }n .
\end{align}
We shall prove that $\overline C_n\geq0$ and $\underline C_n\geq0$ for $n$ large, which contradicts \eqref{sd28}.

\textbf{Step 1.} We show that the sets $S:=\{v_\lambda=\bar v\}$ and $S^\p:=\{v_\lambda=\sigma z_\lambda\}$ have measure zero. All three functions $\underline v$, $v_\lambda$ and $\bar v$ lie in $W^{2,p}_{loc}(\Omega)$, so $\Delta(\bar v-v_\lambda)=0$ a.e.\ on $S$ and $\Delta(v_\lambda-\underline v)=0$ a.e.\ on $S^\p$. We have
\begin{align*}
-\Delta\bar v=\bar\lambda\bigl(\rho(\bar v)+\mathcal{K}[\bar v]\Psi_t(x,\bar v)\bigr),\quad -\Delta v_\lambda=\lambda\bigl(\rho(v_\lambda)+\mathcal{K}[v_\lambda]\Psi_t(x,v_\lambda)\bigr)\qquad\text{a.e.\ in }\Omega .
\end{align*}
On $S$, the arguments of $\rho$ and $\Psi_t$ coincide, and $\mathcal{K}[v_\lambda]\leq\mathcal{K}[\bar v]$ since $v_\lambda\leq\bar v$ in $\Omega$. Hence, a.e.\ on $S$,
\begin{align*}
-\Delta(\bar v-v_\lambda)\geq(\bar\lambda-\lambda)\Bigl(\chi_{\{g(v_\lambda)<b\}}g(v_\lambda)^{-\beta}+\mathcal{K}[v_\lambda]\,h(x,g(v_\lambda))\Bigr)g^\p(v_\lambda)>0 ,
\end{align*}
so $|S|=0$. On $S^\p$, we have $g(\sigma z_\lambda)\leq\sigma z_\lambda<b/2$. Since $g^{-\beta}g^\p$ is nonincreasing and $\sigma<1$, we get, a.e.\ on $S^\p$,
\begin{align*}
-\Delta(v_\lambda-\sigma z_\lambda)\geq\lambda g(\sigma z_\lambda)^{-\beta}g^\p(\sigma z_\lambda)-\sigma\lambda g(z_\lambda)^{-\beta}g^\p(z_\lambda)\geq(1-\sigma)\lambda\, g(\sigma z_\lambda)^{-\beta}g^\p(\sigma z_\lambda)>0 ,
\end{align*}
so $|S^\p|=0$.

\textbf{Step 2.} We prove that $|\overline R_n|,|\underline R_n|\to0$ and $\|\nabla\bar z_n\|_2,\|\nabla\underline z_n\|_2\to0$ as $n \to \infty$. Since $u\mapsto u^{\pm}$ is continuous on $H^1_0(\Omega)$ and $v_\lambda\in\mathcal{N}_0$, as $n \to \infty$,
\begin{align*}
\bar z_n=(v_n-\bar v)^+\to(v_\lambda-\bar v)^+=0,\quad \underline z_n=(v_n-\underline v)^-\to(v_\lambda-\underline v)^-=0\quad\text{in }H^1_0(\Omega).
\end{align*}
For $\delta>0$,
\begin{align*}
\overline R_n\subset\{v_\lambda>\bar v-\delta\}\cup\{|v_n-v_\lambda|\geq\delta\}.
\end{align*}
Here $|\{|v_n-v_\lambda|\geq\delta\}|\leq\delta^{-2}\|v_n-v_\lambda\|_2^2\to0$ as $n\to \infty$, while $|\{v_\lambda>\bar v-\delta\}|\downarrow|\{v_\lambda\geq\bar v\}|=|S|=0$ as $\delta\to0$ by Step~1. Hence $|\overline R_n|\to0$, and likewise $|\underline R_n|\to0$.

\textbf{Step 3.} We show that $\overline C_n\geq0$ for large $n$. On $\overline R_n$, we have $v_n=\bar v+\bar z_n$ and $w_n=\bar v$.
Since $\Phi$ is concave on $(0,\infty)$,
\begin{align*}
\Phi(\bar v+\bar z_n)-\Phi(\bar v)\leq\rho(\bar v)\bar z_n\quad\text{on }\overline R_n ,
\end{align*}
so the second term of $\overline C_n$ is nonnegative.
We need an upper bound for $H(y,g(v_n))+H(y,g(w_n))$. Recall that $H(y,\cdot)$ and $g$ are nondecreasing, $w_n\leq\bar v$ in $\Omega$, $v_n\leq\bar v$ in $\Omega\setminus\overline R_n$ and $v_n=\bar v+\bar z_n$ in $\overline R_n$. Hence, for a.e.\ $y\in\Omega$,
\begin{align}\label{sd29}
H(y,g(v_n))+H(y,g(w_n))\leq2H(y,g(\bar v))+\chi_{\overline R_n}(y)\bigl\{H(y,g(\bar v+\bar z_n))-H(y,g(\bar v))\bigr\}.
\end{align}
We integrate \eqref{sd29} against $|x-y|^{-\mu}$, multiply by the nonnegative factor $\{H(x,g(\bar v+\bar z_n))-H(x,g(\bar v))\}$, $x\in\overline R_n$, and insert the result in $\overline C_n$. The two terms containing $\mathcal{K}[\bar v]$ combine into a Taylor remainder, and we obtain
\begin{align}\label{sd30}
\begin{split}
\overline C_n\geq{}&\frac12\|\nabla\bar z_n\|_2^2-\lambda\int_{\overline R_n}\mathcal{K}[\bar v](x)\Bigl[H(x,g(\bar v+\bar z_n))-H(x,g(\bar v))-h(x,g(\bar v))g^\p(\bar v)\bar z_n\Bigr]dx\\
&-\frac\lambda2\int_{\overline R_n}\int_{\overline R_n}\frac{\bigl\{H(x,g(\bar v+\bar z_n))-H(x,g(\bar v))\bigr\}\bigl\{H(y,g(\bar v+\bar z_n))-H(y,g(\bar v))\bigr\}}{|x-y|^\mu}dxdy .
\end{split}
\end{align}
Both integrals in \eqref{sd30} are of second order in $\bar z_n$.
We estimate the two terms involving $\mathcal{K}[\bar v]$ together, since taking them separately would lead to first-order terms that cannot be absorbed by the quadratic term $\frac12\|\nabla\bar z_n\|_2^2$.

On $\overline R_n$, we write
\begin{align*}
H(x,g(\bar v+\bar z_n))-H(x,g(\bar v))=\int_{\bar v}^{\bar v+\bar z_n}h(x,g(s))g^\p(s)\,ds .
\end{align*}
The integrand is nondecreasing in $s$, and it is Lipschitz on $(0,t^*]$ by (b). Moreover, since $g^\p\leq1$, we obtain, on $\overline R_n$,
\begin{align}\label{sd31}
\begin{split}
0\leq H(x,g(\bar v+\bar z_n))-H(x,g(\bar v))&\leq h(x,g(v_n))\,\bar z_n,\\
0\leq H(x,g(\bar v+\bar z_n))-H(x,g(\bar v))-h(x,g(\bar v))g^\p(\bar v)\bar z_n&\leq
\begin{cases}
\frac L2\,\bar z_n^2, & \text{if }\bar z_n\leq1,\\ h(x,g(v_n))\,\bar z_n^2, & \text{if }\bar z_n>1.
\end{cases}
\end{split}
\end{align}
For $\bar z_n\leq1$, we used $\bar v+\bar z_n\leq t^*$, and for $\bar z_n>1$, we used the first inequality and $\bar z_n\leq\bar z_n^2$. Also,  using $(f_2)$, $(g_3)$, and $\bar v\leq t^*$, we deduce, on $\overline R_n$,
\begin{align*}
h(x,g(v_n))\leq C_\epsilon e^{(1+\epsilon)g(v_n)^{4\alpha}}\leq C_\epsilon e^{2\alpha(1+\epsilon)(\bar v+\bar z_n)^2}\leq Ce^{4\alpha(1+\epsilon)\bar z_n^2}.
\end{align*}
Since $\|\nabla\bar z_n\|_2\to0$ by Step~2, the Trudinger--Moser inequality \cite{moser1971} gives, for every $q<\infty$, some $n_q$ such that
\begin{align}\label{sd32}
\int_{\overline R_n}h(x,g(v_n))^q\,dx\leq C_q, \quad\text{for }n\geq n_q .
\end{align}

By H\"older's inequality, the Sobolev embedding $H^1_0(\Omega)\hookrightarrow L^4(\Omega)$, (a), \eqref{sd31} and \eqref{sd32} with $q=4$, the first integral in \eqref{sd30} is at most
\begin{align*}
C\int_{\overline R_n}\bar z_n^2dx+C\int_{\overline R_n}h(x,g(v_n))\bar z_n^2dx
&\leq C|\overline R_n|^{\frac12}\|\bar z_n\|_4^2+C|\overline R_n|^{\frac14}\Bigl(\int_{\overline R_n}h(x,g(v_n))^4dx\Bigr)^{\frac14}\|\bar z_n\|_4^2\\
&\leq C|\overline R_n|^{\frac14}\|\nabla\bar z_n\|_2^2\quad(n\geq n_4).
\end{align*}
For the second integral, let $p:=\frac{4}{4-\mu}$. We use the Hardy--Littlewood--Sobolev inequality \cite[Theorem 4.3]{lieb2001}, the first inequality in \eqref{sd31}, H\"older's inequality with exponents $2$, $4$, $4$ applied to $h(x,g(v_n))^p\cdot\bar z_n^p\cdot\chi_{\overline R_n}$, and \eqref{sd32} with $q=2p$. Then the second integral is at most
\begin{align*}
C\Bigl(\int_{\overline R_n}\bigl(h(x,g(v_n))\bar z_n\bigr)^pdx\Bigr)^{\frac2p}
&\leq C\Bigl(\int_{\overline R_n}h(x,g(v_n))^{2p}dx\Bigr)^{\frac1p}\|\bar z_n\|_{4p}^2\,|\overline R_n|^{\frac1{2p}}\\
&\leq C|\overline R_n|^{\frac1{2p}}\|\nabla\bar z_n\|_2^2\quad(n\geq n_{2p}).
\end{align*}
Inserting the two estimates in \eqref{sd30}, we get
\begin{align*}
\overline C_n\geq\Bigl(\frac12-C\lambda|\overline R_n|^{\frac14}-C\lambda|\overline R_n|^{\frac1{2p}}\Bigr)\|\nabla\bar z_n\|_2^2\geq0,
\end{align*}
for large $n$, by Step~2.

\textbf{Step 4.} We prove that $\underline C_n\geq0$ for large $n$. On $\underline R_n$, we have $v_n=\underline v-\underline z_n$ and $w_n=\underline v$, with $0\leq v_n<\underline v$, and $v_n>0$ a.e.\ when $\beta\geq1$.
Since $\Phi$ is concave on $(0,\infty)$, and on $[0,\infty)$ when $\beta<1$,
\begin{align*}
\Phi(\underline v-\underline z_n)-\Phi(\underline v)\leq-\rho(\underline v)\underline z_n\quad\text{on }\underline R_n ,
\end{align*}
so the second term of $\underline C_n$ is nonnegative.
Also, since $w_n\geq\underline v$ in $\Omega$, $v_n\geq\underline v$ in $\Omega\setminus\underline R_n$ and $v_n=\underline v-\underline z_n$ in $\underline R_n$, we have, for a.e.\ $y\in\Omega$,
\begin{align}\label{sd33}
H(y,g(v_n))+H(y,g(w_n))\geq2H(y,g(\underline v))-\chi_{\underline R_n}(y)\bigl\{H(y,g(\underline v))-H(y,g(\underline v-\underline z_n))\bigr\}.
\end{align}
We integrate \eqref{sd33} against $|x-y|^{-\mu}$, multiply by the nonnegative factor $\{H(x,g(\underline v))-H(x,g(\underline v-\underline z_n))\}$, $x\in\underline R_n$, and insert the result in $\underline C_n$. We obtain
\begin{align}\label{sd34}
\begin{split}
\underline C_n\geq{}&\frac12\|\nabla\underline z_n\|_2^2-\lambda\int_{\underline R_n}\mathcal{K}[\underline v](x)\Bigl[h(x,g(\underline v))g^\p(\underline v)\underline z_n-H(x,g(\underline v))+H(x,g(\underline v-\underline z_n))\Bigr]dx\\
&-\frac\lambda2\int_{\underline R_n}\int_{\underline R_n}\frac{\bigl\{H(x,g(\underline v))-H(x,g(\underline v-\underline z_n))\bigr\}\bigl\{H(y,g(\underline v))-H(y,g(\underline v-\underline z_n))\bigr\}}{|x-y|^\mu}dxdy .
\end{split}
\end{align}
On $\underline R_n$, all the arguments lie in $[0,t^*]$, so (b) gives
\begin{align*}
&0\leq H(x,g(\underline v))-H(x,g(\underline v-\underline z_n))=\int_{\underline v-\underline z_n}^{\underline v}h(x,g(s))g^\p(s)\,ds\leq C\,\underline z_n,\\
&0\leq h(x,g(\underline v))g^\p(\underline v)\underline z_n-H(x,g(\underline v))+H(x,g(\underline v-\underline z_n))\\
&\qquad=\int_{\underline v-\underline z_n}^{\underline v}\bigl(h(x,g(\underline v))g^\p(\underline v)-h(x,g(s))g^\p(s)\bigr)ds\leq\frac L2\,\underline z_n^2 .
\end{align*}
By \eqref{sd34}, (a), the Hardy--Littlewood--Sobolev inequality, H\"older's inequality (note that $p<4$) and the Sobolev embedding, we deduce
\begin{align*}
\underline C_n&\geq\frac12\|\nabla\underline z_n\|_2^2-C\int_{\underline R_n}\underline z_n^2\,dx-C\Bigl(\int_{\underline R_n}\underline z_n^p\,dx\Bigr)^{\frac2p}\\
&\geq\frac12\|\nabla\underline z_n\|_2^2-C|\underline R_n|^{\frac12}\|\underline z_n\|_4^2-C|\underline R_n|^{\frac2p-\frac12}\|\underline z_n\|_4^2\\
&\geq\Bigl(\frac12-C|\underline R_n|^{\frac12}-C|\underline R_n|^{\frac2p-\frac12}\Bigr)\|\nabla\underline z_n\|_2^2\geq0,
\end{align*}
for large $n$, by Step~2. No Trudinger--Moser argument is needed here, because on $\underline R_n$ all the arguments of $h(x,g(\cdot))$ stay in $[0,t^*]$.

By Steps~3 and~4, we have $\overline C_n\geq0$ and $\underline C_n\geq0$ for large $n$, and \eqref{sd28} gives $J_\lambda(v_\lambda)>J_\lambda(v_n)\geq J_\lambda(v_\lambda)$, a contradiction. Hence $v_\lambda$ is a local minimizer of $J_\lambda$ in the $H^1_0(\Omega)$ topology. \qed
\end{proof}

We now treat the endpoint $\lambda=\lambda^*_b$.

\begin{lemma}\label{lemsd21}
Problem \eqref{sd6} has a weak solution for $\lambda=\lambda^*_b$. More precisely, let $\lambda_n\uparrow\lambda^*_b$ with $\lambda^*_b/2\leq\lambda_n<\lambda^*_b$, fix $\sigma\in(0,1)$ with $\sigma\|z_{\lambda^*_b}\|_\infty<b/2$, and let $v_n$ be the solution of $(\mathcal{Q}^b_{\lambda_n})$ given by Lemma~\ref{lemsd10} with this $\sigma$. Then $\{v_n\}$ is bounded in $H^1_0(\Omega)$ and, along a subsequence, $v_n\rightharpoonup v$ in $H^1_0(\Omega)$. The limit satisfies $v\in W^{2,p}_{loc}(\Omega)$ for every $p<\infty$, $v\geq\sigma\underline z\geq\sigma\kappa\varphi_\beta$, where $\underline z:=z_{\lambda^*_b/2}$, and $|\{g(v)=b\}|=0$. Finally, $v$ is a weak solution of $(\mathcal{Q}^b_{\lambda^*_b})$.
\end{lemma}
\begin{proof}
Since $z_{\lambda_n}\leq z_{\lambda^*_b}$ by Lemma~\ref{lemsd6}, this $\sigma$ can be used in Lemma~\ref{lemsd9} for every $n$. Hence
\begin{align*}
v_n\geq\sigma z_{\lambda_n}\geq\sigma\underline z\geq\sigma\kappa\varphi_\beta ,
\end{align*}
where $\kappa$ is the constant of Lemma~\ref{lemsd3}(a) for $\underline z$. By Lemma~\ref{lemsd10}(iii), $v_n$ minimizes $J_{\lambda_n}$ on a set that contains $\sigma z_{\lambda_n}$. Moreover, $G(g(\sigma z_{\lambda_n}))\geq G(g(\sigma\kappa\varphi_\beta))\in L^1(\Omega)$ and, since $z_{\lambda_n}\leq z_{\lambda^*_b}$, we obtain
\begin{align*}
\|\nabla z_{\lambda_n}\|_2^2=\lambda_n\int_\Omega\rho(z_{\lambda_n})z_{\lambda_n}dx\leq\lambda^*_b\int_\Omega g(\kappa\varphi_\beta)^{-\beta}z_{\lambda^*_b}\,dx\leq C .
\end{align*}
Hence
\begin{align*}
J_{\lambda_n}(v_n)\leq J_{\lambda_n}(\sigma z_{\lambda_n})\leq\frac12\|\nabla\sigma z_{\lambda_n}\|_2^2-\lambda_n\int_\Omega G(g(\sigma\kappa\varphi_\beta))dx\leq C .
\end{align*}

We test the equation of $v_n$ with $v_n$, which is admissible by Remark~\ref{remsd1} and Lemma~\ref{lemsd7}. Substituting the resulting identity into $J_{\lambda_n}(v_n)$, we get
\begin{align*}
J_{\lambda_n}(v_n)&=\frac{\lambda_n}{2}\int_\Omega\rho(v_n)v_n\,dx-\lambda_n\int_\Omega\Phi(v_n)dx\\
&\qquad+\frac{\lambda_n}{2}\int_\Omega\int_\Omega\frac{H(y,g(v_n))\bigl[h(x,g(v_n))g^\p(v_n)v_n-H(x,g(v_n))\bigr]}{|x-y|^\mu}dxdy .
\end{align*}
By $(g_4)$, we have $g^\p(t)t\geq g(t)/(2\alpha)$, while $(f_3)$ yields $H(x,s)\leq M(1+h(x,s))$. Moreover, for $s\geq4\alpha M$, $Mh(x,s)\leq\frac{s}{4\alpha}h(x,s)$, whereas for $s\leq4\alpha M$, the continuity of $h$ gives $h(x,s)\leq C$, uniformly in $x\in \bar{\Omega}.$ Hence
\begin{align*}
h(x,s)g^\p(t)t-H(x,s)\geq\frac{s}{4\alpha}h(x,s)-C_h\quad\text{for }s=g(t),\ t>0 .
\end{align*}
Let
\begin{align*}
Q_n:=\int_\Omega\int_\Omega\frac{H(y,g(v_n))h(x,g(v_n))g(v_n(x))}{|x-y|^\mu}dxdy,\quad c_H:=\int_\Omega H(y,g(\sigma\underline z))dy>0 .
\end{align*}
Then
\[
\int_\Omega H(y,g(v_n))dy\geq c_H, \quad Q_n\geq d_\Omega^{-\mu}c_H\int_\Omega h(x,g(v_n))g(v_n)dx.
\]
Moreover, $\Phi(t)=G(g(t))\leq \displaystyle \max_{s\in \mathbb{R}} G(s)<\infty$ if $\beta\geq1$ and $\Phi(t)\leq Ct^{1-\beta}$ if $\beta<1$, and $\int_\Omega\rho(v_n)v_n\,dx\geq0$. Therefore, with $C_\mu$ as in Subsection~\ref{secsd1}, we deduce
\begin{align*}
\frac{\lambda_n}{8\alpha}Q_n\leq C+\frac{\lambda_nC_hC_\mu}{2}\int_\Omega H(y,g(v_n))dy+C\bigl(1+\|\nabla v_n\|_2^{(1-\beta)^+}\bigr).
\end{align*}
For $R>0$, $(f_3)$ gives
\begin{align*}
\int_\Omega H(y,g(v_n))dy\leq C_R+\frac{M}{R}\int_\Omega h(y,g(v_n))g(v_n)dy\leq C_R+\frac{Md_\Omega^\mu}{Rc_H}Q_n .
\end{align*}
Choosing $R$ large and using $\lambda_n\geq\lambda^*_b/2$, we get $Q_n\leq C(1+\|\nabla v_n\|_2^{(1-\beta)^+})$. Testing the equation again with $v_n$, and using $g^\p(t)t\leq g(t)$ together with the estimate
\[
\int_\Omega\rho(v_n)v_n\,dx\leq\int_\Omega g(\sigma\kappa\varphi_\beta)^{-\beta}v_n\,dx\leq C\|\nabla v_n\|_2
\]
by Remark~\ref{remsd1}, we get
\[
\|\nabla v_n\|_2^2\leq\lambda_nC\|\nabla v_n\|_2+\lambda_nQ_n\leq C\bigl(1+\|\nabla v_n\|_2\bigr).
\]
Hence $\{v_n\}$ is bounded in $H^1_0(\Omega)$, and $Q_n\leq C$. In particular,
\begin{align}\label{sd78}
\int_\Omega h(x,g(v_n))g(v_n)\,dx\leq C,\qquad \int_\Omega\mathcal{K}[v_n]\,h(x,g(v_n))g(v_n)\,dx=Q_n\leq C .
\end{align}

Along a subsequence, $v_n\rightharpoonup v$ in $H^1_0(\Omega)$ and $v_n\to v$ a.e. Hence $v\geq\sigma\underline z\geq\sigma\kappa\varphi_\beta$, and $g(v_n)\to g(v)$ a.e. By the first bound in \eqref{sd78} and the lemma of de Figueiredo, Miyagaki and Ruf \cite[Lemma 2.1]{defigueiredo1995}, $h(\cdot,g(v_n))\to h(\cdot,g(v))$ in $L^1(\Omega)$. By $(f_3)$ and the generalized dominated convergence theorem, $H(\cdot,g(v_n))\to H(\cdot,g(v))$ in $L^1(\Omega)$. Since $|x|^{-\mu}\in L^1_{loc}(\R^2)$, $\mathcal{K}[v_n]\to\mathcal{K}[v]$ in $L^1(\Omega)$ and, along a further subsequence, a.e.

The second bound in \eqref{sd78} gives, for every $R>0$,
\begin{align*}
\int_{\{v_n\geq R\}}\mathcal{K}[v_n]h(x,g(v_n))\,dx\leq\frac{C}{g(R)} .
\end{align*}
On $\{v_n<R\}$, the integrand is at most $\displaystyle\max_{\bar\Omega}h(\cdot,g(R))\,\mathcal{K}[v_n]$, and $\{\mathcal{K}[v_n]\}$ is uniformly integrable because it converges in $L^1(\Omega)$. Hence $\{\mathcal{K}[v_n]h(\cdot,g(v_n))\}$ is uniformly integrable over $\Omega$. Since $g^\p\leq1$ and $g^\p(v_n)\to g^\p(v)$ a.e., Vitali's theorem gives
\begin{align}\label{sd79}
\mathcal{K}[v_n]\,h(\cdot,g(v_n))\,g^\p(v_n)\to\mathcal{K}[v]\,h(\cdot,g(v))\,g^\p(v)\quad\text{in }L^1(\Omega).
\end{align}

Let $\phi\in C^\infty_0(\Omega)$ with support $K$. On $K$,
\begin{align*}
0\leq\rho(v_n)\leq g(\min_K\sigma\kappa\varphi_\beta)^{-\beta}=:C_K .
\end{align*}
Along a further subsequence, $\rho(v_n)\rightharpoonup\tilde\rho$ weakly-$*$ in $L^\infty(K)$, with $0\leq\tilde\rho\leq C_K$. By a diagonal argument, we may assume that this holds for every compact $K\subset\Omega$, with one function $\tilde\rho\in L^\infty_{loc}(\Omega)$. Passing to the limit in the weak identity for $v_n$, and using \eqref{sd79}, we get
\begin{align*}
\int_\Omega\nabla v\cdot\nabla\phi\,dx=\lambda^*_b\int_\Omega\bigl(\tilde\rho+\mathcal{K}[v]\,h(x,g(v))g^\p(v)\bigr)\phi\,dx\quad\text{for all }\phi\in C^\infty_0(\Omega).
\end{align*}
By \eqref{sd85} and $g^\p\leq1$, $\mathcal{K}[v]h(\cdot,g(v))g^\p(v)\in L^p(\Omega)$ for every $p<\infty$. Hence $-\Delta v\in L^p_{loc}(\Omega)$, and $v\in W^{2,p}_{loc}(\Omega)$ by elliptic regularity.

On the level set $\{v=t_b\}$, we have $\Delta v=0$ a.e. On the other hand, $\tilde\rho\geq0$ and $\mathcal{K}[v]\geq d_\Omega^{-\mu}c_H>0$, so, a.e.\ on $\{v=t_b\}$,
\begin{align*}
-\Delta v=\lambda^*_b\bigl(\tilde\rho+\mathcal{K}[v]h(x,g(v))g^\p(v)\bigr)\geq\lambda^*_b\mathcal{K}[v]h(x,b)g^\p(t_b)>0 .
\end{align*}
Hence $|\{g(v)=b\}|=|\{v=t_b\}|=0$. Since $\rho$ is continuous on $(0,\infty)\setminus\{t_b\}$, we get $\rho(v_n)\to\rho(v)$ a.e.\ in $\Omega$. With the local bound $C_K$, dominated convergence identifies the weak-$*$ limit as $\tilde\rho=\rho(v)$ a.e. Thus $v$ satisfies \eqref{sd7} for every $\phi\in C^\infty_0(\Omega)$ with $\lambda=\lambda^*_b$. The positivity requirement of Definition~\ref{defsd2} holds because $v\geq\sigma\underline z$. \qed
\end{proof}

\noindent\textbf{Proof of Theorem~\ref{thmsd1}(i)--(ii).}
Part (i) is precisely Lemma~\ref{lemsd8}.
For part (ii), first suppose that \eqref{sd1} has a weak solution $u$. By Lemma~\ref{lemsd2}, $v=g^{-1}(u)$ is a weak solution of \eqref{sd6}, so $\lambda\leq\lambda^*_b$ by \eqref{sd9}.  Conversely, let $\lambda\in(0,\lambda^*_b)$. Lemma~\ref{lemsd10}(i) gives a weak solution $v_\lambda\in\mathcal{N}_0$ of \eqref{sd6} with $v_\lambda\in W^{2,p}_{loc}(\Omega)$ and $v_\lambda\geq\sigma z_\lambda\geq\sigma\kappa\varphi_\beta$. Lemma~\ref{lemsd12} shows that $v_\lambda$ is a local minimizer of $J_\lambda$ in $H^1_0(\Omega)$, and Lemma~\ref{lemsd2} shows that $u_\lambda=g(v_\lambda)$ is a weak solution of \eqref{sd1}. For $\lambda=\lambda^*_b$, Lemma~\ref{lemsd21} gives a weak solution $v$ of \eqref{sd6} with $v\in W^{2,p}_{loc}(\Omega)$ and $v\geq\sigma\kappa\varphi_\beta$, and Lemma~\ref{lemsd2} applies again. \qed
	\section{Existence of a second solution}
	In this section, we prove Theorem~\ref{thmsd1}(iii). Throughout, $\lambda\in(0,\lambda^*_b)$ is fixed, and $v_\lambda$ is the solution given by Lemma~\ref{lemsd10}, which is a local minimizer of $J_\lambda$ by Lemma~\ref{lemsd12}. We look for a second solution of the form $v_\lambda+w$ with $w\geq0$, $w\not\equiv0$.

If $v$ solves the problem below, then $v_\lambda+v$ solves \eqref{sd6}, and it is a second solution if $v\not\equiv0$:
\begin{align*}{\label{sd35}}
\begin{cases}
-\Delta v&= \lambda \left(\chi_{\{g(v+v_\lambda)<b\}} g(v+v_\lambda)^{-\beta} g^\p(v+v_\lambda)-\chi_{\{g(v_\lambda)<b\}} g(v_\lambda)^{-\beta} g^\p(v_\lambda)\right)\\& \qquad+ \lambda\left(
\begin{array}{l}
\left(\int_\Omega \frac{H(y,g(v+v_\lambda))}{\left|x-y\right|^\mu}dy\right)h(x, g(v+v_\lambda)) g^\p(v+v_\lambda)\\[2mm] \qquad -\left(\int_\Omega \frac{H(y,g(v_\lambda))}{\left|x-y\right|^\mu}dy\right)h(x, g(v_\lambda)) g^\p(v_\lambda)
\end{array}
\right)  \text{ in } \Omega,\\& v>0 \text{ in } \Omega, \quad v=0 \text{ on } \partial\Omega.
\end{cases}
\tag{$\widetilde{\mathcal{Q}}_{\lambda}^b$}
\end{align*}

For $x\in \Omega$ and $s\in\R$, let
\begin{align*}
\tilde{f}(x,s)&=\left(\chi_{\{g(s+v_\lambda)<b\}} g(s+v_\lambda)^{-\beta} g^\p(s+v_\lambda)-\chi_{\{g(v_\lambda)<b\}} g(v_\lambda)^{-\beta} g^\p(v_\lambda)\right)\,\chi_{\mathbb{R}^+}(s),\\ \widetilde{F}(x,t)&=\int_0^t \tilde{f}(x,s)ds .
\end{align*}
Thus $\tilde f(x,s)=\rho(v_\lambda+s)-\rho(v_\lambda)$ for $s>0$, and $\tilde f(x,s)=0$ for $s\leq0$. Since $\rho$ is nonincreasing, $\tilde f(x,\cdot)$ is nonincreasing and nonpositive on $\R$. Hence $\widetilde F(x,\cdot)$ is nonpositive and concave on $\R$, and $\widetilde F(x,s)=0$ for $s\leq0$.

The energy functional $\mathcal{J}_\lambda:H^1_0(\Omega)\to\mathbb{R}$ associated with \eqref{sd35} is
\begin{align*}
\mathcal{J}_\lambda(v)&=\frac12\int_\Omega|\nabla v|^2dx-\lambda\int_\Omega\widetilde{F}(x,v)dx
+\lambda\int_\Omega\mathcal{K}[v_\lambda]\,h(x,g(v_\lambda))g^\p(v_\lambda)\,v\,dx\\
&\qquad-\frac{\lambda}{2}\int_\Omega\int_\Omega\frac{H(y,g(v_\lambda+v))H(x,g(v_\lambda+v))-H(y,g(v_\lambda))H(x,g(v_\lambda))}{|x-y|^\mu}dxdy .
\end{align*}

\begin{prop}\label{propsd1}
Let $\psi(v):=-\lambda\int_\Omega\widetilde F(x,v)\,dx$ and $\Lambda:=\mathcal{J}_\lambda-\psi$. Then $\psi$ is convex, nonnegative and Lipschitz continuous on $H^1_0(\Omega)$, and $\Lambda$ is of class $C^1$ on $H^1_0(\Omega)$, with
\begin{align}\label{sd36}
\langle\Lambda^\p(v),\phi\rangle=\int_\Omega\nabla v\cdot\nabla\phi\,dx-\lambda\int_\Omega\bigl(\mathcal{K}[v_\lambda+v]\,\Psi_t(x,v_\lambda+v)-\mathcal{K}[v_\lambda]\,\Psi_t(x,v_\lambda)\bigr)\phi\,dx .
\end{align}
In particular, $\mathcal{J}_\lambda=\Lambda+\psi$ is locally Lipschitz on $H^1_0(\Omega)$.
\end{prop}

\begin{proof}
By $(f_1)$, $\Psi(x,t)=0$ for $t\leq0$, and $\Psi_t(x,t)=h(x,g(t))g^\p(t)\to h(x,0)g^\p(0)=0$ as $t\downarrow0$ because $f(x,0)=0$. Hence $\Psi(x,\cdot)\in C^1(\mathbb{R})$. Fix $u\in H^1_0(\Omega)$. By Trudinger's inequality \cite{trudinger1967}, $e^{qu^2}\in L^1(\Omega)$ for every $q<\infty$. Since $e^{q(u+w)^2}\leq e^{2qu^2}e^{2qw^2}$, Moser's inequality \cite{moser1971}, together with H\"older's inequality, bounds $\int_\Omega e^{q(u+w)^2}dx$ uniformly for $\|\nabla w\|_2\leq\delta(q)$. With $(g_3)$, $(f_2)$ and $(f_3)$, this gives $L^q(\Omega)$ bounds for both $\Psi(\cdot,v_\lambda+u+w)$ and $\Psi_t(\cdot,v_\lambda+u+w)$, for every $q<\infty$, uniformly for $\|\nabla w\|_2\leq\delta(q)$. By the Hardy--Littlewood--Sobolev inequality \cite[Theorem 4.3]{lieb2001}, the double integral in $\mathcal{J}_\lambda$ is then of class $C^1$. Its derivative is obtained by differentiation under the integral sign, and the symmetry of the kernel produces the factor $2$ that cancels the factor $\frac12$. The third term of $\mathcal{J}_\lambda$ is linear in $v$, and \eqref{sd36} follows.

By the properties of $\tilde f$ above, $|\tilde f(x,s)|\leq\rho(v_\lambda)$. Hence
\begin{align*}
\Bigl|\int_\Omega\widetilde{F}(x,v+w)dx-\int_\Omega\widetilde{F}(x,v)dx\Bigr|
\leq\int_\Omega\Bigl|\int_v^{v+w}|\tilde f(x,s)|\,ds\Bigr|\,dx\leq\int_\Omega\rho(v_\lambda)|w|dx\leq C\|\nabla w\|_2 .
\end{align*}
The last inequality is the estimate of Remark~\ref{remsd1}, which applies because $v_\lambda\geq\sigma z_\lambda\geq\sigma\kappa\varphi_\beta$ by Lemmas~\ref{lemsd10}(i) and~\ref{lemsd3}(a). The convexity of $\psi$ follows from the concavity of $\widetilde F(x,\cdot)$, and $\psi\geq0$ from $\widetilde F\leq0$. \qed
\end{proof}

We follow the scheme of \cite[Section 5]{anthal2025}, where the same difficulty arises for a Choquard term. For $v\in H^1_0(\Omega)$, we write
\begin{align*}
V:=v_\lambda+v,\quad D(V):=\int_\Omega\int_\Omega\frac{H(y,g(V))H(x,g(V))}{|x-y|^\mu}dxdy,\quad c^*:=\frac{(4-\mu)\pi}{4\alpha},
\end{align*}
and $X^+:=\{v\in H^1_0(\Omega):v\geq0\text{ a.e.\ in }\Omega\}$, a closed convex cone, hence weakly closed. For $s\geq0$, we have
\[
\widetilde F(x,s)=\Phi(v_\lambda+s)-\Phi(v_\lambda)-\rho(v_\lambda)s.
\]
Consequently, testing the equation of $v_\lambda$ with $v\in X^+$, we obtain
\begin{align}\label{sd37}
\mathcal{J}_\lambda(v)=J_\lambda(v_\lambda+v)-J_\lambda(v_\lambda)\quad\text{for all }v\in X^+ .
\end{align}
By Lemma~\ref{lemsd12}, $0$ is a local minimizer of $\mathcal{J}_\lambda$ on $X^+$, and $\mathcal{J}_\lambda(0)=0$. The identity \eqref{sd37} fails off $X^+$, which is why the second solution is sought in $X^+$.

For $x\in\Omega$ and $s\in\R$, let $\underline{\tilde f}(x,s):=\tilde f(x,s+0)$ and $\overline{\tilde f}(x,s):=\tilde f(x,s-0)$ denote the right and left limits of the monotone function $\tilde f(x,\cdot)$ at $s$. Then $\underline{\tilde f}\leq\tilde f\leq\overline{\tilde f}$.  Since $\rho$ is nonincreasing and $\tilde f\equiv0$ on $(-\infty,0]$,
\begin{align}\label{sd38}
-\rho(v_\lambda)\leq\underline{\tilde f}\leq\overline{\tilde f}\leq0\quad\text{on }[0,\infty),\quad \underline{\tilde f}=\overline{\tilde f}=0\quad\text{on }(-\infty,0).
\end{align}
Moreover, $\underline{\tilde f}(x,0)=\overline{\tilde f}(x,0)=0$ for a.e.\ $x\in\Omega$, because $|\{v_\lambda=t_b\}|=0$ by Lemma~\ref{lemsd10}(i) and $\rho$ is continuous at every $t\neq t_b$.

By Proposition~\ref{propsd1}, $\mathcal{J}_\lambda=\Lambda+\psi$, with $\Lambda\in C^1(H^1_0(\Omega))$ and $\psi$ convex and continuous. Adding to $\psi$ the indicator function of the closed convex set $X^+$, which is $0$ on $X^+$ and $+\infty$ elsewhere, we are in the setting of Szulkin \cite{szulkin1986}. We use the corresponding notion of critical point.

We say that $w\in X^+$ is a \emph{critical point of $\mathcal{J}_\lambda$ on $X^+$} if
\begin{align}\label{sd39}
\langle\Lambda^\p(w),\ell-w\rangle+\psi(\ell)-\psi(w)\geq0\quad\text{for all }\ell\in X^+ .
\end{align}
A local minimizer $w$ of $\mathcal{J}_\lambda$ on $X^+$ is a critical point on $X^+$. Indeed, for $\ell\in X^+$ and small $t>0$, we have $w+t(\ell-w)\in X^+$, and by the convexity of $\psi$,
\begin{align*}
0\leq\mathcal{J}_\lambda(w+t(\ell-w))-\mathcal{J}_\lambda(w)&\leq\Lambda(w+t(\ell-w))-\Lambda(w)+t\bigl(\psi(\ell)-\psi(w)\bigr)\\
&=t\bigl(\langle\Lambda^\p(w),\ell-w\rangle+\psi(\ell)-\psi(w)\bigr)+o(t).
\end{align*}
Dividing by $t$ gives \eqref{sd39}.

We next compute the directional derivatives of $\psi$. For $x\in\Omega$, $s,r\in\R$, let
\begin{align}\label{sd40}
\Theta(x,s,r):=
\begin{cases}
-\underline{\tilde f}(x,s)\,r, & r\geq0,\\ -\overline{\tilde f}(x,s)\,r, & r\leq0 .
\end{cases}
\end{align}
By \eqref{sd38}, $|\Theta(x,s,r)|\leq\rho(v_\lambda)|r|$. Let $w,\phi\in H^1_0(\Omega)$. Since $\tilde f(x,\cdot)$ is monotone, as $\epsilon\downarrow0$,
\begin{align*}
\frac1\epsilon\int_{w}^{w+\epsilon\phi}\tilde f(x,\tau)\,d\tau\to
\begin{cases}
\tilde f(x,w+0)\,\phi & \text{on }\{\phi>0\},\\ \tilde f(x,w-0)\,\phi & \text{on }\{\phi<0\},
\end{cases}
\end{align*}
and these quotients are bounded by $\rho(v_\lambda)|\phi|\in L^1(\Omega)$ (Remark~\ref{remsd1}). By dominated convergence and the convexity of $\psi$,
\begin{align}\label{sd41}
\psi^\p(w;\phi):=\lim_{\epsilon\downarrow0}\frac{\psi(w+\epsilon\phi)-\psi(w)}{\epsilon}=\lambda\int_\Omega\Theta(x,w,\phi)\,dx,\quad \psi(w+\phi)-\psi(w)\geq\psi^\p(w;\phi).
\end{align}

Since $\widetilde F(x,s)=0$ for $s\leq0$, we have $\widetilde F(x,\max\{s,0\})=\widetilde F(x,s)$. Hence
\begin{align}\label{sd80}
\psi(u^+)=\psi(u)\quad\text{for all }u\in H^1_0(\Omega).
\end{align}

\begin{lemma}\label{lemsd13}
Let $w_k\rightharpoonup w$ in $H^1_0(\Omega)$. Then $\int_\Omega\rho(v_\lambda)|w_k-w|dx\to0$. In particular, $\psi$ is weakly continuous on $H^1_0(\Omega)$.
\end{lemma}
\begin{proof}
By Lemmas~\ref{lemsd10}(i) and~\ref{lemsd3}(a), $v_\lambda\geq\sigma z_\lambda\geq\sigma\kappa\varphi_\beta$. Hence $\rho(v_\lambda)\leq g(\sigma\kappa\varphi_\beta)^{-\beta}$ by $(g_2)$ and the monotonicity of $g$, and
\begin{align*}
\rho(v_\lambda)\,d\leq d\,g(\sigma\kappa\varphi_\beta)^{-\beta}\in L^2(\Omega).
\end{align*}
The membership holds precisely because $\beta<3$. Indeed, near $\partial\Omega$, the right-hand side behaves like $d^{\frac{1-\beta}{1+\beta}}$ for $\beta>1$, like $(\log(1/d))^{-1/2}$ for $\beta=1$, and like $d^{1-\beta}$ for $\beta<1$, and $\frac{2(\beta-1)}{\beta+1}<1$ if and only if $\beta<3$. By Hardy's inequality \cite[(0.5)]{brezis1997}, $\{d^{-1}|w_k-w|\}$ is bounded in $L^2(\Omega)$. Every subsequence has a further subsequence converging to $0$ a.e., hence weakly in $L^2(\Omega)$, so the whole sequence converges weakly to $0$. Pairing the two factors in
\begin{align*}
\int_\Omega\rho(v_\lambda)|w_k-w|dx=\int_\Omega\bigl(\rho(v_\lambda)d\bigr)\frac{|w_k-w|}{d}dx
\end{align*}
gives the first assertion. The second follows from $|\widetilde F(x,s)-\widetilde F(x,s^\p)|\leq\rho(v_\lambda)|s-s^\p|$. \qed
\end{proof}

\begin{lemma}\label{lemsd14}
Let $w\in X^+$ be a critical point of $\mathcal{J}_\lambda$ on $X^+$, that is, suppose that $w$ satisfies \eqref{sd39}. Then $V=v_\lambda+w\in W^{2,p}_{loc}(\Omega)$ for every $p<\infty$, $|\{g(V)=b\}|=0$, and $V$ is a weak solution of \eqref{sd6}. Moreover, if $w\not\equiv0$, then $|\{w=0\}|=0$, and $u=g(V)$ is a second weak solution of \eqref{sd1}.
\end{lemma}

\begin{proof}
Let $\phi\in C^\infty_0(\Omega)$. For $\epsilon>0$, define
\[
\vartheta_\epsilon:=(w+\epsilon\phi)^-, \quad  \Omega_\epsilon:=\{w+\epsilon\phi<0\}.
\]
Then $\ell:=(w+\epsilon\phi)^+=w+\epsilon\phi+\vartheta_\epsilon\in X^+$, $\ell-w=\epsilon\phi+\vartheta_\epsilon$, and $\psi(\ell)=\psi(w+\epsilon\phi)$ by \eqref{sd80}. Taking this $\ell$ in \eqref{sd39} and using the linearity of $\Lambda^\p(w)$, we obtain
\begin{align*}
0\leq\epsilon\langle\Lambda^\p(w),\phi\rangle+\langle\Lambda^\p(w),\vartheta_\epsilon\rangle+\psi(w+\epsilon\phi)-\psi(w) .
\end{align*}
We show that the second term is at most $o(\epsilon)$. First, since $\nabla\vartheta_\epsilon=-\nabla(w+\epsilon\phi)\chi_{\Omega_\epsilon}$,
\begin{align*}
\int_\Omega\nabla w\cdot\nabla\vartheta_\epsilon\,dx=-\int_{\Omega_\epsilon}\bigl(|\nabla w|^2+\epsilon\nabla w\cdot\nabla\phi\bigr)dx\leq-\epsilon\int_{\Omega_\epsilon}\nabla w\cdot\nabla\phi\,dx=o(\epsilon),
\end{align*}
because $\chi_{\Omega_\epsilon}\to\chi_{\{w=0\}\cap\{\phi<0\}}$ a.e.\ and $\nabla w=0$ a.e.\ on $\{w=0\}$. Second, consider the bracket
\begin{align*}
\mathcal{B}(w):=\mathcal{K}[V]\,\Psi_t(x,V)-\mathcal{K}[v_\lambda]\,\Psi_t(x,v_\lambda)
\end{align*}
in \eqref{sd36}. It is nonnegative, because $V\geq v_\lambda$, $H(y,\cdot)$ is nondecreasing and $\Psi_t(x,\cdot)$ is nondecreasing on $[0,\infty)$. Hence $-\lambda\int_\Omega \mathcal{B}(w)\vartheta_\epsilon\, dx\leq0$. Dividing by $\epsilon$, passing to the limit as $\epsilon\to0$, and using \eqref{sd41}, we obtain
\begin{align}\label{sd43}
\langle\Lambda^\p(w),\phi\rangle+\lambda\int_\Omega\Theta(x,w,\phi)dx\geq0\quad\text{for all }\phi\in C^\infty_0(\Omega).
\end{align}

Let $0\leq\phi\in C^\infty_0(\Omega)$. By \eqref{sd40}, $\Theta(x,w,\phi)=-\underline{\tilde f}(x,w)\phi$ and $\Theta(x,w,-\phi)=\overline{\tilde f}(x,w)\phi$. Applying \eqref{sd43} to $\phi$ and to $-\phi$, we get
\begin{align}\label{sd44}
\lambda\int_\Omega\underline{\tilde f}(x,w)\phi\,dx\leq\langle\Lambda^\p(w),\phi\rangle\leq\lambda\int_\Omega\overline{\tilde f}(x,w)\phi\,dx .
\end{align}
Since $v_\lambda$ is a weak solution of \eqref{sd6},
\begin{align*}
\int_\Omega\nabla v_\lambda\cdot\nabla\phi\,dx=\lambda\int_\Omega\rho(v_\lambda)\phi\,dx+\lambda\int_\Omega\mathcal{K}[v_\lambda]\Psi_t(x,v_\lambda)\phi\,dx .
\end{align*}
On adding this identity to \eqref{sd36}, we get
\begin{align*}
\langle\Lambda^\p(w),\phi\rangle+\lambda\int_\Omega\rho(v_\lambda)\phi\,dx=\int_\Omega\nabla V\cdot\nabla\phi\,dx-\lambda\int_\Omega\mathcal{K}[V]\Psi_t(x,V)\phi\,dx .
\end{align*}
Adding $\lambda\int_\Omega\rho(v_\lambda)\phi\,dx$ to each of the three terms of \eqref{sd44}, we obtain
\begin{align*}
\lambda\bigl(\rho(v_\lambda)+\underline{\tilde f}(x,w)\bigr)\leq-\Delta V-\lambda\mathcal{K}[V]\,\Psi_t(x,V)\leq\lambda\bigl(\rho(v_\lambda)+\overline{\tilde f}(x,w)\bigr)\quad \text{in }\mathcal{D}^\p(\Omega).
\end{align*}

The outer terms lie in $L^\infty_{loc}(\Omega)$. Hence, the middle term, a distribution squeezed between them, is a function in $L^\infty_{loc}(\Omega)$. Indeed, its difference from the left-hand side is a nonnegative Radon measure bounded by an $L^\infty_{loc}$-density. By \eqref{sd85}, applied to $V\in H^1_0(\Omega)$, and $g^\p\leq1$, we have $\mathcal{K}[V]\Psi_t(\cdot,V)\in L^p_{loc}(\Omega)$. It follows that $-\Delta V\in L^p_{loc}(\Omega)$ and $V\in W^{2,p}_{loc}(\Omega)$ for every $p<\infty$, and the two inequalities hold a.e.\ in $\Omega$.

Recall that $\rho$ is continuous at every $t>0$, $t\neq t_b$. Wherever $\tilde f(x,\cdot)$ is continuous at $s$, its one-sided limits at $s$ coincide with $\tilde f(x,s)$. We consider three cases.
First, on the set where $w>0$ and $V\neq t_b$, we have $\underline{\tilde f}(x,w)=\overline{\tilde f}(x,w)=\rho(V)-\rho(v_\lambda)$, and both bounds equal $\lambda\rho(V)$. Next, on $\{w=0\}$, the sentence following \eqref{sd38} shows that both bounds equal $\lambda\rho(v_\lambda)=\lambda\rho(V)$ a.e. Finally, on $\{V=t_b\}$, we have $D^2V=0$ a.e. Indeed, the gradient of a function in $W^{1,1}_{loc}(\Omega)$ vanishes a.e.\ on each of its level sets, and we apply this to $V$ and to each $\partial_iV$. In particular, $\Delta V=0$ a.e.\ on $\{V=t_b\}$. On the other hand, $\rho(v_\lambda)+\underline{\tilde f}(x,w)\geq0$ by \eqref{sd38}, and the left inequality yields
\begin{align*}
-\Delta V\geq\lambda\mathcal{K}[V]\Psi_t(x,t_b)=\lambda\mathcal{K}[V]h(x,b)g^\p(t_b)>0 \text{ a.e.\ on }\{V=t_b\}.
\end{align*}
Together with $\Delta V=0$ a.e.\ on $\{V=t_b\}$, this forces $|\{g(V)=b\}|=|\{V=t_b\}|=0$. By the first two cases, $-\Delta V=\lambda\rho(V)+\lambda\mathcal{K}[V]\,\Psi_t(x,V) \text{ a.e.\ in }\Omega.$
Thus, $V$ satisfies \eqref{sd7}, and the positivity in Definition~\ref{defsd2} holds because $V\geq v_\lambda$.

Since $w=V-v_\lambda\in W^{2,p}_{loc}(\Omega)$, $\Delta w=0$ a.e.\ on $\{w=0\}$. Subtracting the equations satisfied by $V$ and $v_\lambda$, we get, a.e.\ on $\{w=0\}$, where $V=v_\lambda$,
\begin{align*}
-\Delta w=\lambda\bigl(\mathcal{K}[V]-\mathcal{K}[v_\lambda]\bigr)\Psi_t(x,v_\lambda).
\end{align*}
As $w\geq0$, $w\not\equiv0$ and $H(y,\cdot)$ is strictly increasing on $(0,\infty)$ by $(f_1)$, we have $\mathcal{K}[V]-\mathcal{K}[v_\lambda]>0$ in $\Omega$. Moreover, $\Psi_t(x,v_\lambda)=h(x,g(v_\lambda))g^\p(v_\lambda)>0$. As the left-hand side vanishes a.e.\ on $\{w=0\}$, we conclude that $|\{w=0\}|=0$. By Lemma~\ref{lemsd2}, $u=g(V)$ is a weak solution of \eqref{sd1}, and $u\neq g(v_\lambda)$ by the injectivity of $g$. \qed
\end{proof}

\begin{lemma}\label{lemsd15}
Let $\{v_n\}\subset X^+$ and $\epsilon_n\to0$ satisfy $\mathcal{J}_\lambda(v_n)\to c$ and
\begin{align}\label{sd45}
\langle\Lambda^\p(v_n),\ell-v_n\rangle+\psi(\ell)-\psi(v_n)\geq-\epsilon_n\max\{1,\|\nabla(\ell-v_n)\|_2\}\quad\text{for all }\ell\in X^+ .
\end{align}
Then $\{v_n\}$ is bounded in $H^1_0(\Omega)$ and, with $V_n:=v_\lambda+v_n$,
\begin{align}\label{sd46}
\int_\Omega h(x,g(V_n))g(V_n)dx+\int_\Omega\mathcal{K}[V_n]\,h(x,g(V_n))g(V_n)dx+D(V_n)\leq C .
\end{align}
\end{lemma}

\begin{proof}
Set $X_n:=\int_\Omega\mathcal{K}[V_n]\Psi_t(x,V_n)v_n\,dx$. Take $\ell=2v_n$ in \eqref{sd45}. Since $-\tilde f\leq\rho(v_\lambda)$,
\begin{align*}
\psi(2v_n)-\psi(v_n)=-\lambda\int_\Omega\int_{v_n}^{2v_n}\tilde f(x,\tau)\,d\tau\,dx\leq\lambda\int_\Omega\rho(v_\lambda)v_n\,dx\leq C\|\nabla v_n\|_2 ,
\end{align*}
by Remark~\ref{remsd1}. Hence $\langle\Lambda^\p(v_n),v_n\rangle\geq-C(1+\|\nabla v_n\|_2)$. By the definition of $\mathcal{J}_\lambda$ and \eqref{sd36},
\begin{align*}
\tfrac\lambda2\bigl[X_n-D(V_n)\bigr]=\mathcal{J}_\lambda(v_n)-\tfrac12\langle\Lambda^\p(v_n),v_n\rangle+\lambda\int_\Omega\widetilde F(x,v_n)dx-\tfrac\lambda2D(v_\lambda)-\tfrac\lambda2\int_\Omega\mathcal{K}[v_\lambda]\Psi_t(x,v_\lambda)v_n\,dx .
\end{align*}
The last three terms are nonpositive, since $\widetilde F\leq0$ and $v_n\geq0$. Hence
\begin{align}\label{sd47}
\tfrac\lambda2\bigl[X_n-D(V_n)\bigr]\leq C\bigl(1+\|\nabla v_n\|_2\bigr).
\end{align}

Fix $R_0:=2\|v_\lambda\|_{L^\infty(\Omega)}+1$ and let $m_0:=\displaystyle \min_{x\in\bar\Omega}h(x,g(R_0))$, which is positive by $(f_1)$. Let $R\geq R_0$. On $\{V_n\geq R\}$, one has $v_\lambda\leq R/2\leq V_n/2$, hence $v_n\geq V_n/2$, and $(g_4)$ gives
\begin{align*}
g^\p(V_n)v_n\geq\tfrac12g^\p(V_n)V_n\geq\frac{g(V_n)}{4\alpha}\geq\frac{g(R)}{4\alpha}.
\end{align*}
Therefore
\begin{align}\label{sd48}
h(x,g(V_n))\leq\frac{4\alpha}{g(R)}\,\Psi_t(x,V_n)v_n\quad\text{on }\{V_n\geq R\}.
\end{align}
On this set, we also have $h(x,g(V_n))\geq h(x,g(R_0))\geq m_0$ by $(f_4)$. Hence $(f_3)$ and \eqref{sd48} give
\begin{align}\label{sd49}
H(x,g(V_n))\leq M\Bigl(1+\frac1{m_0}\Bigr)h(x,g(V_n))\leq\frac{C_0}{g(R)}\,\Psi_t(x,V_n)v_n\quad\text{on }\{V_n\geq R\},
\end{align}
where $C_0:=4\alpha M(1+1/m_0)$. On $\{V_n<R\}$, we use instead $0\leq g(V_n)\leq g(R)$ and the monotonicity of $h$ and $H$:
\begin{align*}
h(x,g(V_n))\leq C_R:=\max_{x\in\bar\Omega}h(x,g(R)),\qquad H(x,g(V_n))\leq\bar H_R:=\max_{x\in\bar\Omega}H(x,g(R)).
\end{align*}
Write $H(\cdot,g(V_n))=H_n^<+H_n^{\geq}$, where $H_n^<:=H(\cdot,g(V_n))\chi_{\{V_n<R\}}$ and $H_n^{\geq}:=H(\cdot,g(V_n))\chi_{\{V_n\geq R\}}$. Since $0\leq H_n^<\leq H(\cdot,g(V_n))$, we have $\int_\Omega H_n^<(x)|x-y|^{-\mu}dx\leq\mathcal{K}[V_n](y)$. Using the symmetry of the kernel and \eqref{sd49},
\begin{align*}
D(V_n)&=\int_\Omega\int_\Omega\frac{H_n^<(x)H_n^<(y)}{|x-y|^\mu}dxdy+\int_\Omega H_n^{\geq}(y)\Bigl(\int_\Omega\frac{H_n^<(x)}{|x-y|^\mu}dx\Bigr)dy\\
&\qquad+\int_\Omega H_n^{\geq}(x)\,\mathcal{K}[V_n](x)\,dx\\
&\leq\bar H_R^2C_\mu|\Omega|+2\int_\Omega H_n^{\geq}(x)\,\mathcal{K}[V_n](x)\,dx\leq\bar H_R^2C_\mu|\Omega|+\frac{2C_0}{g(R)}X_n .
\end{align*}

The constant $C_0$ does not depend on $R$. We now choose $R\geq R_0$ with $2C_0/g(R)\leq\frac12$ and keep it fixed. From now on, constants may depend on $R$. Inserting the last bound into \eqref{sd47} gives $X_n\leq C(1+\|\nabla v_n\|_2)$, hence $D(V_n)\leq C(1+\|\nabla v_n\|_2)$. Since $\mathcal{J}_\lambda(v_n)\to c$, we deduce $\frac12\|\nabla v_n\|_2^2\leq C(1+\|\nabla v_n\|_2)$. Thus $\{v_n\}$ is bounded in $H^1_0(\Omega)$, and $X_n$, $D(V_n)$ are bounded.

We now prove \eqref{sd46}. Let $c_H:=d_\Omega^{-\mu}\int_\Omega H(y,g(v_\lambda))dy$. Since $V_n\geq v_\lambda$, we have $\mathcal{K}[V_n]\geq c_H>0$, and \eqref{sd49} gives
\begin{align*}
\int_\Omega H(x,g(V_n))dx\leq\bar H_R|\Omega|+\frac{C_0}{g(R)}\int_\Omega\Psi_t(x,V_n)v_n\,dx\leq\bar H_R|\Omega|+\frac{C_0}{g(R)c_H}X_n\leq C .
\end{align*}
On $\{V_n\geq R\}$ we have $h(x,g(V_n))g(V_n)\leq4\alpha\,\Psi_t(x,V_n)v_n$, by the inequality before \eqref{sd48}. On $\{V_n<R\}$, we have $h(x,g(V_n))g(V_n)\leq C_Rg(R)$. Since $\mathcal{K}[V_n]\geq c_H$,
\begin{align*}
\int_\Omega h(x,g(V_n))g(V_n)\,dx\leq\frac{4\alpha}{c_H}X_n+C_Rg(R)|\Omega|.
\end{align*}
Finally, $\int_\Omega\mathcal{K}[V_n]\,dx\leq C_\mu\int_\Omega H(y,g(V_n))\,dy$ by Tonelli's theorem, with $C_\mu$ as in Subsection~\ref{secsd1}, so
\begin{align*}
\int_\Omega\mathcal{K}[V_n]\,h(x,g(V_n))g(V_n)\,dx\leq4\alpha X_n+C_Rg(R)\,C_\mu\int_\Omega H(y,g(V_n))\,dy .
\end{align*}
Both right-hand sides are bounded, which proves \eqref{sd46}. \qed
\end{proof}

\begin{lemma}\label{lemsd16}
Let $\{v_n\}$ be as in Lemma~\ref{lemsd15}. Then, along a subsequence, $v_n\rightharpoonup v$ in $H^1_0(\Omega)$ with $v\in X^+$, and the conclusions of Lemma~\ref{lemsd14} hold for $w=v$. Furthermore, if $0<c<c^*$, then $v\not\equiv0$.
\end{lemma}

\begin{proof}
By Lemma~\ref{lemsd15}, we may assume $v_n\rightharpoonup v$, $v_n\to v$ a.e.\ and in $L^q(\Omega)$ for all $q<\infty$. Moreover, $v\in X^+$ because $X^+$ is weakly closed, and $V_n\to V=v_\lambda+v$ a.e. The bound \eqref{sd46} contains the two bounds in \eqref{sd78}. Hence the argument that follows \eqref{sd78} in the proof of Lemma~\ref{lemsd21}, from the application of \cite[Lemma 2.1]{defigueiredo1995} onwards, applies with $v_n$ replaced by $V_n$. It shows that $\mathcal{K}[V_n]\to\mathcal{K}[V]$ in $L^1(\Omega)$ and, along a further subsequence, a.e., and that $\{\mathcal{K}[V_n]h(\cdot,g(V_n))\}$ is uniformly integrable. By $(f_3)$, $\{\mathcal{K}[V_n]H(\cdot,g(V_n))\}$ is uniformly integrable as well, and Vitali's theorem gives
\begin{align}\label{sd50}
\mathcal{K}[V_n]\,h(\cdot,g(V_n))\to\mathcal{K}[V]\,h(\cdot,g(V))\ \text{ in }L^1(\Omega),\quad D(V_n)\to D(V) \text{ as } n\to\infty.
\end{align}

Fix $\phi\in C^\infty_0(\Omega)$ and $\epsilon>0$. Let $\vartheta_n:=(v_n+\epsilon\phi)^-$ and $\upsilon_n:=\epsilon\phi+\vartheta_n$, and take $\ell=(v_n+\epsilon\phi)^+=v_n+\upsilon_n$ in \eqref{sd45}. Then $0\leq\vartheta_n\leq\epsilon|\phi|$, $\{\vartheta_n\}$ is bounded in $H^1_0(\Omega)$, and $\vartheta_n\rightharpoonup\vartheta:=(v+\epsilon\phi)^-$; put $\upsilon:=\epsilon\phi+\vartheta$. We pass to the limit in the three terms of the left-hand side of \eqref{sd45}.

For the first term, $\nabla v_n\cdot\nabla\vartheta_n=-|\nabla\vartheta_n|^2-\epsilon\nabla\phi\cdot\nabla\vartheta_n$ a.e., and weak lower semicontinuity of the norm gives
\begin{align*}
\limsup_{n\to\infty}\int_\Omega\nabla v_n\cdot\nabla\upsilon_n\,dx\leq\epsilon\int_\Omega\nabla v\cdot\nabla\phi\,dx-\|\nabla\vartheta\|_2^2-\epsilon\int_\Omega\nabla\phi\cdot\nabla\vartheta\,dx=\int_\Omega\nabla v\cdot\nabla\upsilon\,dx .
\end{align*}
For the second term, we use $|\upsilon_n|\leq2\epsilon|\phi|$, $\upsilon_n\to\upsilon$ a.e., \eqref{sd50}, $g^\p\leq1$ and Vitali's theorem to obtain
\[
\int_\Omega \mathcal{B}(v_n)\upsilon_n\,dx\to\int_\Omega \mathcal{B}(v)\upsilon\,dx,
\]
with $\mathcal{B}$ as in the proof of Lemma~\ref{lemsd14}.
For the third term, \eqref{sd80} shows that $\psi(\ell)-\psi(v_n)=\psi(v_n+\epsilon\phi)-\psi(v_n)$, and this tends to $\psi(v+\epsilon\phi)-\psi(v)$, because $\psi$ is weakly continuous by Lemma~\ref{lemsd13}.

Since the right-hand side of \eqref{sd45} tends to $0$, and $\psi(v+\epsilon\phi)=\psi((v+\epsilon\phi)^+)$ by \eqref{sd80},
\begin{align*}
\langle\Lambda^\p(v),\epsilon\phi+\vartheta\rangle+\psi\bigl((v+\epsilon\phi)^+\bigr)-\psi(v)\geq0 .
\end{align*}
This is \eqref{sd39} for $w=v$ and $\ell=(v+\epsilon\phi)^+$, the only choice of $\ell$ used in the proof of Lemma~\ref{lemsd14}. Hence the conclusions of Lemma~\ref{lemsd14} hold for $w=v$.

Suppose now $0<c<c^*$ and $v\equiv0$. Then $V=v_\lambda$, and $D(V_n)\to D(v_\lambda)$ by \eqref{sd50}. Moreover, $\int_\Omega|\widetilde F(x,v_n)|dx\leq\int_\Omega\rho(v_\lambda)|v_n|dx\to0$ by Lemma~\ref{lemsd13}, and the linear term tends to $0$. Hence $\frac12\|\nabla v_n\|_2^2\to c$, so for $n$ large
\begin{align*}
\|\nabla v_n\|_2^2\leq\frac{(4-\mu)\pi}{2\alpha}(1-\delta)\quad\text{for some }\delta>0 .
\end{align*}
Let $p:=\frac{4}{4-\mu}$ and choose $\epsilon>0$ and $s>1$ with $s(1+\epsilon)^2(1-\delta)<1$. By $(f_2)$, $(g_3)$ and $V_n\leq\|v_\lambda\|_\infty+v_n$,
\begin{align*}
h(x,g(V_n))\leq C_\epsilon e^{2\alpha(1+\epsilon)V_n^2}\leq Ce^{2\alpha(1+\epsilon)^2v_n^2},\qquad 2\alpha ps(1+\epsilon)^2\|\nabla v_n\|_2^2\leq4\pi s(1+\epsilon)^2(1-\delta)<4\pi .
\end{align*}
Hence Moser's inequality \cite{moser1971} shows that $h(\cdot,g(V_n))$, and by $(f_3)$ also $H(\cdot,g(V_n))$, are bounded in $L^{ps}(\Omega)$. By the Hardy--Littlewood--Sobolev inequality \cite[Theorem 4.3]{lieb2001} and H\"older's inequality,
\begin{align*}
X_n=\int_\Omega\mathcal{K}[V_n]\,\Psi_t(x,V_n)v_n\,dx\leq C\|H(\cdot,g(V_n))\|_{p}\,\|h(\cdot,g(V_n))\|_{ps}\,\|v_n\|_{ps^\p}\to0 .
\end{align*}
Now take $\ell=0$ in \eqref{sd45}. Since $\psi(0)-\psi(v_n)=\lambda\int_\Omega\widetilde F(x,v_n)\,dx\leq0$, we get $\langle\Lambda^\p(v_n),v_n\rangle\leq C\epsilon_n$, that is,
\begin{align*}
\|\nabla v_n\|_2^2\leq C\epsilon_n+\lambda X_n-\lambda\int_\Omega\mathcal{K}[v_\lambda]\Psi_t(x,v_\lambda)v_n\,dx\leq C\epsilon_n+\lambda X_n\to0 .
\end{align*}
Hence $c=\lim\mathcal{J}_\lambda(v_n)=0$, a contradiction. \qed
\end{proof}

\begin{lemma}\label{lemsd17}
There is $r_1>0$ with the following property. Either $\mathcal{J}_\lambda$ has a critical point $w\not\equiv0$ on $X^+$, or, for every $r_0\in(0,r_1)$, there is $a>0$ such that $\mathcal{J}_\lambda(v)\geq a$ for all $v\in X^+$ with $\|\nabla v\|_2=r_0$.
\end{lemma}

\begin{proof}
By \eqref{sd37} and Lemma~\ref{lemsd12}, there is $r_1\in(0,(\pi/(4\alpha))^{1/2})$ such that $\mathcal{J}_\lambda\geq0=\mathcal{J}_\lambda(0)$ on $X^+\cap B_{r_1}$, where $B_{r_1}:=\{\|\nabla v\|_2<r_1\}$. If $\mathcal{J}_\lambda(w)=0$ for some $w\in(X^+\cap B_{r_1})\setminus\{0\}$, then $w$ is a local minimizer of $\mathcal{J}_\lambda$ on $X^+$, hence a critical point on $X^+$ (see the paragraph following \eqref{sd39}).

Otherwise $\mathcal{J}_\lambda>0$ on $(X^+\cap B_{r_1})\setminus\{0\}$. Fix $r_0\in(0,r_1)$. Then $\inf\{\mathcal{J}_\lambda(v):v\in X^+,\ \|\nabla v\|_2=r_0\}\geq0$, and we show that this infimum is positive. Suppose that it is $0$, and choose $q>2$ with $2\alpha qr_0^2<\pi$. Let $\{v_k\}\subset X^+$ satisfy $\|\nabla v_k\|_2=r_0$ and $\mathcal{J}_\lambda(v_k)\to0$, and, passing to a subsequence, assume that $v_k\rightharpoonup v$. Then $v\in X^+$.

We claim that the non-quadratic part of $\mathcal{J}_\lambda$ is weakly continuous on $\bar B_{r_0}$. For $\psi$ this is Lemma~\ref{lemsd13}. For the double integral, Moser's inequality, as in the proof of Lemma~\ref{lemsd16}, shows that $h(\cdot,g(V_k))$ and $H(\cdot,g(V_k))$ remain bounded in $L^q(\Omega)$. Vitali's theorem then gives $H(\cdot,g(V_k))\to H(\cdot,g(V))$ in $L^{p}(\Omega)$, $p=\frac4{4-\mu}$, and the Hardy--Littlewood--Sobolev inequality yields $D(V_k)\to D(V)$. The linear term is clearly continuous.

Hence $\mathcal{J}_\lambda(v)\leq \displaystyle \liminf_{k \to\infty}\mathcal{J}_\lambda(v_k)=0$, so $\mathcal{J}_\lambda(v)=0$ and $v=0$. But then the non-quadratic part of $\mathcal{J}_\lambda(v_k)$ tends to $0$, and $\frac{r_0^2}{2}=\displaystyle\lim_{k \to\infty}\mathcal{J}_\lambda(v_k)=0$, a contradiction. \qed
\end{proof}
\begin{rem}\label{remsd5}
In \cite{dhanya2014,anthal2025}, the corresponding step is split into the two cases
\begin{align*}
&\text{(ZA)}\quad \inf_{v\in X^+,\,\|\nabla v\|_2=r_0}\mathcal{J}_\lambda=\mathcal{J}_\lambda(0)\ \text{ for all small }r_0,\\
&\text{(MP)}\quad \inf_{v\in X^+,\,\|\nabla v\|_2=r_0}\mathcal{J}_\lambda>\mathcal{J}_\lambda(0)\ \text{ for some }r_0.
\end{align*}
In those papers, (ZA) is treated with Ekeland's variational principle on the annulus $\{v\in X^+:r_0-\delta\leq\|\nabla v\|_2\leq r_0+\delta\}$, and the nontriviality of the limit requires the strong convergence of the resulting sequence. The dichotomy in Lemma~\ref{lemsd17} is finer. In the degenerate case, the critical point is a local minimizer on $X^+$, whereas in the nondegenerate case, (ZA) is excluded by the weak continuity established above. In the degenerate case, no level estimate is needed.
\end{rem}

We now show that the mountain-pass level lies below $c^*$. Here we use $(f_5)$.

\begin{lemma}\label{lemsd18}
Assume $(f_5)$. For every $r_0>0$ there is a continuous path $\gamma:[0,1]\to X^+$ with $\gamma(0)=0$, $\|\nabla\gamma(1)\|_2>r_0$, $\mathcal{J}_\lambda(\gamma(1))<0$ and $\max_{[0,1]}\mathcal{J}_\lambda\circ\gamma<c^*$.
\end{lemma}
\begin{proof}
Since $v_\lambda\in C(\Omega)$ and $|\{v_\lambda=t_b\}|=0$ (Lemma~\ref{lemsd10}(i)), the open set $\{v_\lambda\neq t_b\}$ is nonempty. Hence there are $x_0\in\Omega$ and $r>0$ with $B_r:=B_r(x_0)\Subset\Omega$ such that either $\sup_{B_r}v_\lambda<t_b$ or $\inf_{B_r}v_\lambda>t_b$. Put $\nu:=\inf_{B_r}v_\lambda>0$. For $n\geq2$ let $M_n$ be the Moser function
\begin{align*}
M_n(x):=\frac{1}{\sqrt{2\pi}}
\begin{cases}
\sqrt{\log n}, & |x-x_0|\leq r/n,\\[1mm] \dfrac{\log(r/|x-x_0|)}{\sqrt{\log n}}, & r/n<|x-x_0|\leq r,\\[2mm] 0, & |x-x_0|>r .
\end{cases}
\end{align*}
Then $M_n\in X^+$, $\|\nabla M_n\|_2=1$, and $M_n=m_n:=\sqrt{\log n/(2\pi)}$ on $B_{r/n}$. Since $\int_0^r\log^2(r/s)s\,ds=r^2/4$,
\begin{align}\label{sd51}
\int_\Omega M_n^2\,dx=\frac{r^2\log n}{2n^2}+\frac{1}{\log n}\int_{r/n}^{r}\log^2\Bigl(\frac rs\Bigr)s\,ds\leq\frac{C_r}{\log n}.
\end{align}
For $t\geq0$, we set $w_{n,t}:=tM_n$. Then $\|\nabla w_{n,t}\|_2=t$, and $t\mapsto w_{n,t}$ is continuous in $H^1_0(\Omega)$. Throughout the proof, $w=w_{n,t}$, $V=v_\lambda+w$, and $\xi$ is defined by $t^2=2c^*(1+\xi)$ so that
\begin{align}\label{sd52}
\frac12\|\nabla w\|_2^2=c^*(1+\xi),\quad 4\alpha t^2m_n^2=(4-\mu)(1+\xi)\log n .
\end{align}

We claim that $|\widetilde F(x,s)|\leq C_rs^2$ on $B_r$, for $s\geq0$. If $\inf_{B_r}v_\lambda>t_b$, then $\rho(v_\lambda+\tau)=\rho(v_\lambda)=0$ for $\tau\geq0$, and $\widetilde F\equiv0$ on $B_r$. Otherwise $\sup_{B_r}v_\lambda<t_b$, and we put $\delta_0:=t_b-\sup_{B_r}v_\lambda>0$.

Let $x\in B_r$.
For $0\leq\tau<\delta_0$, both $v_\lambda(x)$ and $v_\lambda(x)+\tau$ lie in $[\nu,t_b)$, where $\rho=g^{-\beta}g^\p$ is Lipschitz with some constant $L_r$, so that $|\tilde f(x,\tau)|\leq L_r\tau$.
For $\tau\geq\delta_0$, $|\tilde f(x,\tau)|\leq\rho(v_\lambda)\leq\rho(\nu)\leq\rho(\nu)\tau/\delta_0$.

Hence $|\tilde f(x,\tau)|\leq C_r\tau$ and $|\widetilde F(x,s)|\leq C_rs^2/2$ on $B_r$. Since $w$ vanishes outside $B_r$, \eqref{sd51} gives
\begin{align}\label{sd53}
-\lambda\int_\Omega\widetilde F(x,w)\,dx\leq\lambda C_rt^2\int_\Omega M_n^2\,dx\leq\frac{Y_2\,t^2}{\log n}\qquad(t\geq0,\ n\geq2).
\end{align}
Write $\Psi(x,V)=\Psi_0(x)+\Psi_1(x)+\mathcal{R}(x)$, where
\begin{align*}
\Psi_0:=\Psi(x,v_\lambda),\qquad \Psi_1:=\Psi_t(x,v_\lambda)w,\qquad \mathcal{R}(x):=\int_{v_\lambda}^{v_\lambda+w}\bigl(\Psi_t(x,s)-\Psi_t(x,v_\lambda)\bigr)ds .
\end{align*}
All three functions are nonnegative. For $\mathcal{R}$, this follows from $w\geq0$ and the monotonicity of $\Psi_t(x,\cdot)$, which is a consequence of $(f_4)$. With $k(x,y):=|x-y|^{-\mu}$ and the symmetry of $k$,
\begin{align*}
D(V)-D(v_\lambda)={}&2\int_\Omega\int_\Omega \Psi_0(y)\Psi_1(x)k\,dxdy+2\int_\Omega\int_\Omega \Psi_0(y)\mathcal{R}(x)k\,dxdy\\
&+\int_\Omega\int_\Omega \Psi_1(x)\Psi_1(y)k\,dxdy+2\int_\Omega\int_\Omega \Psi_1(x)\mathcal{R}(y)k\,dxdy\\
&+\int_\Omega\int_\Omega \mathcal{R}(x)\mathcal{R}(y)k\,dxdy .
\end{align*}
The first term equals $2\int_\Omega\mathcal{K}[v_\lambda]\Psi_t(x,v_\lambda)w\,dx$ and cancels against the linear term of $\mathcal{J}_\lambda$. The next three terms are nonnegative, and we drop them. In the last term we restrict the integral to $B_{r/n}\times B_{r/n}$, where $\mathcal{R}\geq0$ and $|x-y|\leq2r/n$. Hence
\begin{align}\label{sd54}
\begin{split}
-\frac\lambda2\bigl[D(V)-D(v_\lambda)\bigr]+\lambda\int_\Omega\mathcal{K}[v_\lambda]\Psi_t(x,v_\lambda)w\,dx
&\leq-\frac\lambda2\int_\Omega\int_\Omega\frac{\mathcal{R}(x)\mathcal{R}(y)}{|x-y|^\mu}dxdy\\
&\leq-\frac{\lambda}{2}\Bigl(\frac{n}{2r}\Bigr)^{\mu}\Bigl(\int_{B_{r/n}}\mathcal{R}(x)\,dx\Bigr)^2\leq0 .
\end{split}
\end{align}
In particular, by \eqref{sd52} and \eqref{sd53},
\begin{align}\label{sd55}
\mathcal{J}_\lambda(w_{n,t})\leq c^*(1+\xi)+\frac{Y_2t^2}{\log n}\qquad\text{for all }t\geq0,\ n\geq2 .
\end{align}

By $(f_5)$, there are $\theta\in[0,1)$, $c_5>0$ and $s_5>0$ such that
\begin{align}\label{sd56}
\Psi_*(s):=\inf_{x\in\bar\Omega}\Psi(x,s)\geq c_5\,e^{2\alpha s^2-s^\theta}\qquad\text{for }s\geq s_5 .
\end{align}
On $B_{r/n}$, we have $w=tm_n$ and $V\geq T_n:=\nu+tm_n$. Fix $T_1:=\sqrt{c^*}$. We claim that there is $n_0$, independent of $t$, such that
\begin{align}\label{sd57}
\mathcal{R}(x)\geq\tfrac12\Psi_*(T_n)\quad\text{on }B_{r/n},\ \text{for all }t\geq T_1\text{ and }n\geq n_0 .
\end{align}
Indeed, $\Psi(\cdot,v_\lambda)$ and $\Psi_t(\cdot,v_\lambda)$ are bounded on $B_r$ by a constant $C$.
Since $\Psi(x,\cdot)$ is nondecreasing and $V\geq T_n$ on $B_{r/n}$,
\[
\mathcal{R}(x)\geq\Psi_*(T_n)-C(1+tm_n)\geq\Psi_*(T_n)-C(1+T_n)\qquad\text{on }B_{r/n},
\]
because $tm_n=T_n-\nu\leq T_n$. By \eqref{sd56} and $\theta<1$, $\Psi_*(T)/(1+T)\to\infty$
as $T\to\infty$. Hence there is $T_*\geq s_5$, depending only on $C$, $c_5$, $s_5$, $\alpha$ and $\theta$, such that $\Psi_*(T)\geq2C(1+T)$ for $T\geq T_*$. For $t\geq T_1$, we have $T_n\geq tm_n\geq T_1m_n$.
Choose $n_0$ with $T_1m_{n_0}\geq T_*$. Then $T_n\geq T_*$ for all $t\geq T_1$ and $n\geq n_0$,
and \eqref{sd57} follows.

By \eqref{sd54}, \eqref{sd57} and $|B_{r/n}|=\pi r^2/n^2$, for $t\geq T_1$ and $n\geq n_0$,
\begin{align*}
-\frac\lambda2\bigl[D(V)-D(v_\lambda)\bigr]+\lambda\int_\Omega\mathcal{K}[v_\lambda]\Psi_t(x,v_\lambda)w\,dx\leq-\bar c\,r^{4-\mu}\,n^{\mu-4}\,\Psi_*(T_n)^2,\quad \bar c:=\frac{\lambda\pi^2}{2^{3+\mu}} .
\end{align*}
By \eqref{sd52}, $4\alpha T_n^2\geq4\alpha t^2m_n^2+8\alpha\nu tm_n=(4-\mu)(1+\xi)\log n+8\alpha\nu tm_n$. Hence
\begin{align*}
n^{\mu-4}\Psi_*(T_n)^2=n^{\mu-4}e^{4\alpha T_n^2}\bigl[\Psi_*(T_n)e^{-2\alpha T_n^2}\bigr]^2\geq n^{(4-\mu)\xi}e^{8\alpha\nu tm_n}\bigl[\Psi_*(T_n)e^{-2\alpha T_n^2}\bigr]^2 .
\end{align*}
Together with \eqref{sd52} and \eqref{sd53}, we obtain
\begin{align}\label{sd61}
\mathcal{J}_\lambda(w_{n,t})\leq c^*(1+\xi)+\frac{Y_2t^2}{\log n}-A_n(t)\,n^{(4-\mu)\xi}\qquad(t\geq T_1,\ n\geq n_0),
\end{align}
where
\begin{align}\label{sd60}
A_n(t):=\bar c\,r^{4-\mu}\,e^{8\alpha\nu tm_n}\bigl[\Psi_*(T_n)e^{-2\alpha T_n^2}\bigr]^2.
\end{align}
Let $t\geq T_1$ and $n\geq n_0$. Then $tm_n=T_n-\nu$ and $T_n\geq T_*\geq s_5$, so \eqref{sd56} gives
\begin{align*}
A_n(t)\geq \bar cr^{4-\mu}c_5^2e^{-8\alpha\nu^2}\exp\bigl(8\alpha\nu T_n-2T_n^\theta\bigr)\geq c_2\,a_n,
\end{align*}
where
\begin{align*}
a_n:=\inf_{T\geq\nu+T_1m_n}\exp\bigl(8\alpha\nu T-2T^\theta\bigr),\quad c_2:=\bar cr^{4-\mu}c_5^2e^{-8\alpha\nu^2}.
\end{align*}
Since $\nu>0$ and $\theta<1$, $8\alpha\nu T-2T^\theta\to\infty$ as $T\to\infty$. As $T_1m_n\to\infty$, we conclude that $a_n\to\infty$ as $n\to\infty$. The bounds \eqref{sd57} and $A_n(t)\geq c_2a_n$ are the only places where $(f_5)$ is used, and the factor $e^{8\alpha\nu tm_n}$, which comes from the translation by $v_\lambda$, is what makes a threshold in $(f_5)$ unnecessary.

Let $T_0:=2\sqrt{c^*}$. By \eqref{sd52}, $\xi\leq-\frac12$ for $t\leq T_1$ and $\xi\geq1$ for $t\geq T_0$. Put
\begin{align*}
\bar\epsilon_n:=\frac{Y_2T_0^2}{\log n},\qquad Y_3:=\exp\Bigl(\frac{2(4-\mu)Y_2T_0^2}{c^*}\Bigr).
\end{align*}
Fix $n\geq n_0$ so large that
\begin{align}\label{sd62}
\frac{Y_2T_1^2}{\log n}<\frac{c^*}{2},\quad c_2a_n(4-\mu)\log n\geq\max\Bigl\{c^*Y_3,\ 2c^*+\frac{4Y_2c^*}{\log n}\Bigr\};
\end{align}
this is possible because $a_n\to\infty$. We show that $\mathcal{J}_\lambda(w_{n,t})<c^*$ for all $t\geq0$, and that $\mathcal{J}_\lambda(w_{n,t})\to-\infty$ as $t\to\infty$. We distinguish three cases.

(a) Let $0\leq t\leq T_1$. Then, by \eqref{sd55} and \eqref{sd62}, $\mathcal{J}_\lambda(w_{n,t})\leq\frac{c^*}{2}+\frac{Y_2T_1^2}{\log n}<c^*$.

(b) Let $T_1\leq t\leq T_0$. Here $Y_2t^2/\log n\leq\bar\epsilon_n$. If $\xi\leq-2\bar\epsilon_n/c^*$, then \eqref{sd55} gives $\mathcal{J}_\lambda(w_{n,t})\leq c^*-2\bar\epsilon_n+\bar\epsilon_n<c^*$. If $\xi\geq-2\bar\epsilon_n/c^*$, put $\xi^\p:=\xi+2\bar\epsilon_n/c^*\geq0$, and note that $n^{-2(4-\mu)\bar\epsilon_n/c^*}=Y_3^{-1}$. Then \eqref{sd61} and the bound $A_n(t)\geq c_2a_n$ give
\begin{align*}
\mathcal{J}_\lambda(w_{n,t})\leq c^*-\bar\epsilon_n+\Bigl(c^*\xi^\p-\frac{c_2a_n}{Y_3}\,n^{(4-\mu)\xi^\p}\Bigr)\leq c^*-\bar\epsilon_n-\frac{c_2a_n}{Y_3}<c^* .
\end{align*}
Here we used that $\xi^\p\mapsto c^*\xi^\p-\frac{c_2a_n}{Y_3}n^{(4-\mu)\xi^\p}$ is nonincreasing on $[0,\infty)$: its derivative is at most $c^*-\frac{c_2a_n}{Y_3}(4-\mu)\log n\leq0$, by \eqref{sd62}.

(c) Finally, let $t\geq T_0$. Then $\xi\geq1$, $n^{(4-\mu)\xi}\geq1+(4-\mu)\xi\log n$, and
\begin{align*}
\frac{Y_2t^2}{\log n}=\frac{2Y_2c^*(1+\xi)}{\log n}\leq\frac{4Y_2c^*\xi}{\log n}.
\end{align*}
Hence \eqref{sd61}, the bound $A_n(t)\geq c_2a_n$ and \eqref{sd62} give
\begin{align*}
\mathcal{J}_\lambda(w_{n,t})\leq c^*-c_2a_n-\xi\Bigl(c_2a_n(4-\mu)\log n-c^*-\frac{4Y_2c^*}{\log n}\Bigr)\leq c^*-c_2a_n-c^*\xi<c^* .
\end{align*}
The right-hand side tends to $-\infty$ as $t\to\infty$, since $\xi=t^2/(2c^*)-1$.

With this $n$ fixed, choose $T\geq\max\{T_0,r_0+1\}$ with $\mathcal{J}_\lambda(w_{n,T})<0$, and set $\gamma(\tau):=w_{n,\tau T}=\tau TM_n$ for $\tau\in[0,1]$. Then $\gamma$ is continuous, $\gamma(0)=0$, $\gamma([0,1])\subset X^+$, $\|\nabla\gamma(1)\|_2=T>r_0$ and $\mathcal{J}_\lambda(\gamma(1))<0$. By (a)--(c), $\max_{[0,1]}\mathcal{J}_\lambda\circ\gamma<c^*$. \qed
\end{proof}
The next lemma gives a Palais--Smale sequence at the mountain-pass level. We cannot quote the mountain pass theorem of \cite{szulkin1986} directly, because it assumes the Palais--Smale condition, while our Palais--Smale sequences are only known to have weak limits (Lemma~\ref{lemsd16}). What we need is the sequence itself, in the form \eqref{sd63}.

\begin{lemma}\label{lemsd19}
Let $e\in X^+$, let $\mathcal{C}:=\{\zeta\in C([0,1],X^+):\zeta(0)=0,\ \zeta(1)=e\}$ and $c:=\inf_{\zeta\in\mathcal{C}}\max_{[0,1]}\mathcal{J}_\lambda\circ\zeta$, and assume $c>\max\{\mathcal{J}_\lambda(0),\mathcal{J}_\lambda(e)\}$. Then there is $\{w_k\}\subset X^+$ with $\mathcal{J}_\lambda(w_k)\to c$ and
\begin{align}\label{sd63}
\langle\Lambda^\p(w_k),\ell-w_k\rangle+\psi(\ell)-\psi(w_k)\geq-\tfrac1k\max\{1,\|\nabla(\ell-w_k)\|_2\}\quad\text{for all }\ell\in X^+ .
\end{align}
\end{lemma}
The proof of Lemma~\ref{lemsd19} is given in Appendix~\ref{secsd3}.

\noindent\textbf{Proof of Theorem~\ref{thmsd1}(iii).} By Lemma~\ref{lemsd14}, it suffices to find a critical point $w\not\equiv0$ of $\mathcal{J}_\lambda$ on $X^+$, or a sequence as in Lemma~\ref{lemsd15} at a level $c\in(0,c^*)$, to which Lemma~\ref{lemsd16} applies. If the first alternative of Lemma~\ref{lemsd17} holds, which is the degenerate case discussed in Remark~\ref{remsd5}, we are done. Otherwise, fix $r_0\in(0,r_1)$, let $a$ be as in Lemma~\ref{lemsd17}, and let $\gamma$ be the path of Lemma~\ref{lemsd18} for this $r_0$, so that $\|\nabla\gamma(1)\|_2>r_0$. Take $e:=\gamma(1)$ in Lemma~\ref{lemsd19}. Every $\zeta\in\mathcal{C}$ crosses $\{\|\nabla v\|_2=r_0\}$, so $c\geq a>0=\mathcal{J}_\lambda(0)>\mathcal{J}_\lambda(e)$. Moreover $c\leq\max_{[0,1]}\mathcal{J}_\lambda\circ\gamma<c^*$ by Lemma~\ref{lemsd18}. Lemma~\ref{lemsd19} gives a sequence $\{w_k\}\subset X^+$ with $\mathcal{J}_\lambda(w_k)\to c$ satisfying \eqref{sd63}, which is \eqref{sd45} with $\epsilon_k=\frac1k$. Lemma~\ref{lemsd16} then yields, along a subsequence, $w_k\rightharpoonup w\not\equiv0$, and $v_\lambda+w$ is a second weak solution of \eqref{sd6}. By Lemma~\ref{lemsd14}, $u=g(v_\lambda+w)$ is a second weak solution of \eqref{sd1}, and $v_\lambda+w\geq v_\lambda$. \qed

For $\alpha>1$, the defect $\Xi$ of Lemma~\ref{lemsd22} grows only linearly. Therefore, in the proof of Lemma~\ref{lemsd18}, the factor $e^{8\alpha\nu tm_n}$ can compensate the loss $e^{-2\Xi(T_n)}$ when $v_\lambda$ is large enough. This gives a second solution also when $f$ is a power of $u$. Recall that $\delta_\infty$ is defined in \eqref{sd83}.

\begin{prop}\label{propsd2}
Let $\alpha>1$ and $\lambda\in(0,\lambda^*_b)$. Assume that $f$ satisfies $(f_1)$--$(f_4)$ and
\begin{align}\label{sd87}
\liminf_{u\to\infty}\,\inf_{x\in\bar\Omega}\,f(x,u)e^{\epsilon u^{2\alpha}}>0\qquad\text{for every }\epsilon>0 .
\end{align}
Let $v_\lambda$ be the solution given by Lemma~\ref{lemsd10}, and assume that
\begin{align}\label{sd88}
\|v_\lambda\|_{L^\infty(\Omega)}>\frac{\delta_\infty}{(2\alpha)^{1/2}} .
\end{align}
Then \eqref{sd1} has a second weak solution $u=g(v)$ with $v\geq v_\lambda$ and $v\not\equiv v_\lambda$.
\end{prop}
\begin{proof}
In the proof of Theorem~\ref{thmsd1}(iii), $(f_5)$ is used only in Lemma~\ref{lemsd18}. In the proof of Lemma~\ref{lemsd18}, it is used only for \eqref{sd57} and for the bound $A_n(t)\geq c_2a_n$ with $a_n\to\infty$. So it suffices to choose $x_0$ and $r$ suitably and to prove these two facts under \eqref{sd87} and \eqref{sd88}. The rest of the proof of Lemma~\ref{lemsd18}, and the proof of Theorem~\ref{thmsd1}(iii), do not change.

\textbf{Step 1: choice of $x_0$ and $r$.} Put $\nu_0:=\delta_\infty/(2\alpha)^{1/2}$. The function $v_\lambda$ is continuous in $\Omega$, so the set
\[
U:=\{x\in\Omega:\ v_\lambda(x)>\nu_0\}
\]
is open. It is nonempty by \eqref{sd88}. The set $\{v_\lambda=t_b\}$ has measure zero by Lemma~\ref{lemsd10}(i), so it has empty interior. Hence there is $x_0\in U$ with $v_\lambda(x_0)\neq t_b$. By continuity, we can choose $r>0$ such that $\bar B_r(x_0)\subset U$ and either $\sup_{B_r}v_\lambda<t_b$ or $\inf_{B_r}v_\lambda>t_b$. Then
\[
\nu:=\inf_{B_r}v_\lambda>\nu_0 ,
\]
and the proof of Lemma~\ref{lemsd18} holds without change up to \eqref{sd56}.

\textbf{Step 2: a lower bound for $\Psi_*$.} Since $\nu>\nu_0$, we can fix $\epsilon>0$ such that
\[
\kappa_0:=8\alpha\nu-2(2\alpha)^{1/2}(2\delta_\infty+\epsilon)>0 .
\]
By \eqref{sd87}, there are $c_\epsilon>0$ and $u_\epsilon\geq2$ such that
\[
f(x,\tau)\geq c_\epsilon e^{-\epsilon\tau^{2\alpha}}\qquad\text{for }x\in\bar\Omega\text{ and }\tau\geq u_\epsilon-1 .
\]
Let $u\geq u_\epsilon$. In the integral that defines $H(x,u)$, we keep only the interval $[u-u^{1-4\alpha},u]$. On this interval, $\tau^{4\alpha}\geq u^{4\alpha}-4\alpha$, as in the proof of Remark~\ref{remsd7}(b). Hence
\[
H(x,u)\geq c_\epsilon e^{-4\alpha}u^{1-4\alpha}e^{u^{4\alpha}-\epsilon u^{2\alpha}} .
\]
Now let $T>0$ with $u=g(T)\geq u_\epsilon$, and use the notation of the proof of Lemma~\ref{lemsd22}. Since $k_\alpha>0$, we have $\delta(T)\leq\delta_\infty$. Moreover, $u^{2\alpha}\leq(2\alpha)^{1/2}T$ by $(g_3)$. Therefore
\begin{align*}
\Xi(T)&=\delta(T)\bigl(2u^{2\alpha}+\delta(T)\bigr)\leq2(2\alpha)^{1/2}\delta_\infty T+\delta_\infty^2,\\
u^{1-4\alpha}&\geq\bigl((2\alpha)^{1/2}T\bigr)^{-\frac{4\alpha-1}{2\alpha}} .
\end{align*}
We insert these bounds in the lower bound for $H$ and use $u^{4\alpha}=2\alpha T^2-\Xi(T)$. With $T_\epsilon:=g^{-1}(u_\epsilon)$, this gives a constant $c_6>0$ such that
\begin{align}\label{sd89}
\Psi_*(T)e^{-2\alpha T^2}\geq c_6\,T^{-\frac{4\alpha-1}{2\alpha}}\exp\bigl(-(2\alpha)^{1/2}(2\delta_\infty+\epsilon)T\bigr)\qquad\text{for }T\geq T_\epsilon .
\end{align}

\textbf{Step 3: the two facts.} By \eqref{sd89}, $\Psi_*(T)/(1+T)\to\infty$ as $T\to\infty$. This is the only property of $\Psi_*$ that the proof of \eqref{sd57} uses, so \eqref{sd57} holds if we choose $T_*\geq T_\epsilon$ there. Next, let $t\geq T_1$ and $n\geq n_0$. Then $T_n\geq T_*$ and $tm_n=T_n-\nu$, so \eqref{sd60} and \eqref{sd89} give
\[
A_n(t)\geq c_2\,T_n^{-\frac{4\alpha-1}{\alpha}}e^{\kappa_0T_n}\geq c_2a_n ,
\]
where
\[
c_2:=\bar c\,r^{4-\mu}c_6^2e^{-8\alpha\nu^2},\qquad a_n:=\inf\Bigl\{T^{-\frac{4\alpha-1}{\alpha}}e^{\kappa_0T}:\ T\geq\nu+T_1m_n\Bigr\}.
\]
Since $\kappa_0>0$ and $T_1m_n\to\infty$, we have $a_n\to\infty$. \qed
\end{proof}

We comment on the two hypotheses of Proposition~\ref{propsd2}. Since $v_\lambda\geq z_\lambda$ by Lemma~\ref{lemsd7}, condition \eqref{sd88} holds whenever $\|z_\lambda\|_{L^\infty(\Omega)}>\delta_\infty/(2\alpha)^{1/2}$. Condition \eqref{sd87} holds whenever $\inf_xf(x,u)\geq cu^{-k}$ for large $u$, with $c>0$ and $k\geq0$. In particular, $f(x,u)=u^{4\alpha-1}$ for $u>0$ satisfies \eqref{sd87}, and it satisfies $(f_1)$--$(f_4)$ by the argument of Example~\ref{exsd1}(b). Thus Proposition~\ref{propsd2} applies to $h(x,u)=u^{4\alpha-1}e^{u^{4\alpha}}$, which is not covered by Theorem~\ref{thmsd1}(iii).

\begin{rem}\label{remsd6}
We comment on hypothesis $(f_5)$.
\begin{enumerate}[label=\textnormal{(\alph*)}]
\item Condition $(f_5)$ is the analogue, for the transformed nonlinearity, of condition (H3) of \cite{dhanya2014}, which in the semilinear case reads $\lim_{t\to\infty}h(t)te^{t^r}=\infty$ for some $r\in[0,1)$. As there, no threshold is needed. In the proof of Lemma~\ref{lemsd18}, the factor $e^{8\alpha\nu tm_n}$ comes from the translation by the first solution $v_\lambda\geq\nu>0$, and it dominates the factor $e^{-2T_n^\theta}$ allowed by $(f_5)$. When the mountain-pass path starts at $0$, as in \cite{defigueiredo1995}, there is no such factor, and the argument there uses a lower bound with an explicit constant.  Condition (H3) is a lower bound for the nonlinearity, while $(f_5)$ is a lower bound for its primitive. We use $H$ because $H(x,g(\cdot))$ is the quantity that enters the nonlocal term in the proof of Lemma~\ref{lemsd18}, see \eqref{sd56} and \eqref{sd57}. By $(f_3)$, $(f_5)$ also gives a lower bound of the same form for $h(x,g(s))$.
\item By Remark~\ref{remsd7}(b), $(f_5)$ fails for every $f$ of polynomial growth. In the approach of Lemma~\ref{lemsd18}, this loss is intrinsic to the scale of the transformed problem. Indeed, assume that, for some $k\geq0$,
\[
\sup_{x\in\bar\Omega}f(x,u)\leq C(1+u^k)\qquad\text{for all }u>0 .
\]
With $u=g(T)$, we have $e^{u^{4\alpha}}=e^{2\alpha T^2-\Xi(T)}$, and $u\leq T$ by $(g_3)$. Hence $(f_3)$ gives
\[
\Psi_*(T)e^{-2\alpha T^2}\leq M\bigl(1+h(x,u)\bigr)e^{-2\alpha T^2}\leq C(1+T^k)e^{-\Xi(T)},
\]
and the quantity $A_n(t)$ in \eqref{sd60} is at most $C(1+T_n^{2k})e^{8\alpha\nu T_n-2\Xi(T_n)}$. For $\alpha\leq1$, $\Xi(T)/T\to\infty$ by Lemma~\ref{lemsd22} and $(g_5)$, so $A_n(t)\to0$ as $T_n\to\infty$. Thus the gain $e^{8\alpha\nu tm_n}$ from the translation by $v_\lambda$ cannot compensate the factor $e^{-2\Xi(T_n)}$, and \eqref{sd61} gives no level below $c^*$. Lower bounds on $h$ at the scale $e^{u^4}$ of the original variable, such as those of \cite{biswas2024} for $\alpha=1$, allow $f$ of polynomial growth, and therefore they do not suffice for our proof. For $\alpha\leq1$ and $f$ of polynomial growth, in particular for $h(x,u)=u^{4\alpha-1}e^{u^{4\alpha}}$, the existence of a second solution remains open. For $\alpha>1$, $\Xi$ grows only linearly, and Proposition~\ref{propsd2} applies to such $f$ if, in addition, \eqref{sd87} holds.
\end{enumerate}
\end{rem}
	\appendix
	\section{Proofs of auxiliary results}\label{secsd3}
	In this appendix we give the proofs of Lemma~\ref{lemsd1}, Lemma~\ref{lemsd2}, Lemma~\ref{lemsd3}(c), Lemma~\ref{lemsd20} and Lemma~\ref{lemsd19}. They are elementary or follow standard schemes.
	
	\subsection{Proof of Lemma~\ref{lemsd1}}
	\begin{proof}
		Let $\mathcal{T}(u):=\int_0^u(1+2\alpha|r|^{4\alpha-2})^{1/2}dr$, $u\in\R$. Since $\alpha>\frac12$, $\mathcal{T}$ is odd, of class $C^1$ on $\R$ and of class $C^\infty$ on $\R\setminus\{0\}$, and $\mathcal{T}^\p\geq1$. Hence $\mathcal{T}$ is strictly increasing and admits an odd inverse $g:=\mathcal{T}^{-1}$, which solves \eqref{sd82}. Conversely, if $g$ solves \eqref{sd82}, then $(\mathcal{T}\circ g)^\p=1$ on $(0,\infty)$. Together with $g(0)=0$, this forces $g=\mathcal{T}^{-1}$, which proves $(g_1)$. Property $(g_2)$ follows directly from $g^\p=1/\mathcal{T}^\p(g)$. For $(g_3)$--$(g_9)$, oddness allows us to take $t>0$. We put $u:=g(t)>0$, so that $t=\mathcal{T}(u)$ and hence $tg^\p(t)=\mathcal{T}(u)/\mathcal{T}^\p(u)$.
		
		$(g_3)$ By $(g_2)$, $g(0)=0$ and $0<g^\prime\leq1$, so the fundamental theorem of calculus gives $g(t)\leq t$. Also $\mathcal{T}(u)\geq\int_0^u(2\alpha)^{1/2}r^{2\alpha-1}dr=(2\alpha)^{-1/2}u^{2\alpha}$, that is, $u^{2\alpha}\leq(2\alpha)^{1/2}t$.
		
		$(g_4)$ We show that $u/(2\alpha)\leq\mathcal{T}(u)/\mathcal{T}^\p(u)\leq u$. The upper bound holds because $\mathcal{T}^\p$ is nondecreasing on $(0,\infty)$. For the lower bound, the function $2\alpha\mathcal{T}(u)-u\mathcal{T}^\p(u)$ vanishes at $u=0$. Since $u\mathcal{T}^{\p\p}(u)=2\alpha(2\alpha-1)u^{4\alpha-2}/\mathcal{T}^\p(u)$, its derivative is
		\begin{align*}
			(2\alpha-1)\mathcal{T}^\p(u)-u\mathcal{T}^{\p\p}(u)=\frac{2\alpha-1}{\mathcal{T}^\p(u)}>0 .
		\end{align*}
        Thus $2\alpha\mathcal{T}(u)-u\mathcal{T}^\p(u)>0$ for $u>0$, which is the lower bound.
		
		$(g_5)$ By l'H\^opital's rule, $\mathcal{T}(u)/u^{2\alpha}\to(2\alpha)^{-1/2}$ as $u\to\infty$. Thus $g(t)^{2\alpha}/t\to(2\alpha)^{1/2}$ as $t\to\infty$, which is equivalent to $(g_5)$.
		
		$(g_6)$ The derivative $g^\p=(1+2\alpha g^{4\alpha-2})^{-1/2}$ is nonincreasing on $(0,\infty)$, so $g$ is concave there. As $g(0)=0$, this gives $g(t)\geq g(1)t$ for $0<t\leq1$. By $(g_4)$, $\bigl(t^{-1/(2\alpha)}g(t)\bigr)^\p=t^{-1/(2\alpha)-1}\bigl(tg^\p(t)-g(t)/(2\alpha)\bigr)\geq0$, so $g(t)\geq g(1)t^{1/(2\alpha)}$ for $t\geq1$.
		
		$(g_7)$ By \eqref{sd82}, $g^\p\leq(2\alpha g^{4\alpha-2})^{-1/2}=(2\alpha)^{-1/2}g^{1-2\alpha}$.
		
		$(g_8)$ Let $\gamma\geq4\alpha-1$. In terms of $u$, $g(t)^\gamma g^\p(t)/t=u^\gamma/(\mathcal{T}(u)\mathcal{T}^\p(u))$, and $u$ is increasing in $t$. The logarithmic derivative of this function of $u$ is $\gamma/u-\mathcal{T}^{\p\p}/\mathcal{T}^\p-\mathcal{T}^\p/\mathcal{T}$. By the lower bound in $(g_4)$, $\mathcal{T}^\p/\mathcal{T}\leq2\alpha/u$. Moreover,
		\begin{align*}
			\frac{\mathcal{T}^{\p\p}(u)}{\mathcal{T}^\p(u)}=\frac{2\alpha(2\alpha-1)u^{4\alpha-3}}{1+2\alpha u^{4\alpha-2}}<\frac{2\alpha-1}{u}.
		\end{align*}
		Hence the logarithmic derivative exceeds $(\gamma-4\alpha+1)/u\geq0$, and the function is strictly increasing in $u$, and so in $t$.
		
		$(g_9)$ On $(0,\infty)$, $g^{-\gamma}$ is decreasing and $g^\p$ is positive and nonincreasing, hence $t\mapsto g(t)^{-\gamma} g^\prime(t)$ is decreasing. \qed
	\end{proof}
	
	\subsection{Proof of Lemma~\ref{lemsd2}}
	\begin{proof}
		Let $F(x):=\lambda\bigl(\rho(v)+\mathcal{K}[v]\,h(x,g(v))g^\p(v)\bigr)$. On a compact $K\subset\Omega$, $\rho(v)\leq g(\operatorname*{ess\,inf}_K v)^{-\beta}$ is bounded. By \eqref{sd85}, $h(\cdot,g(v))\in L^p(\Omega)$ for every $p<\infty$ and $\mathcal{K}[v]\in L^\infty(\Omega)$. Thus $F\in L^p_{loc}(\Omega)$, and since $-\Delta v=F$ in $\mathcal{D}^\p(\Omega)$, elliptic regularity \cite[Chapter 9]{gilbarg2001} yields $v\in W^{2,p}_{loc}(\Omega)\subset C^{1}_{loc}(\Omega)$. In particular, $v\in L^\infty_{loc}(\Omega)$. By mollification, the identity \eqref{sd7} then holds for every $w\in H^1_0(\Omega)\cap L^\infty(\Omega)$ with compact support, because all terms are bounded on the support of $w$ and $h(\cdot,g(v))\in L^p(\Omega)$.
		
		Set $u=g(v)$. By $(g_2)$, $(g_3)$, and $(g_7)$, $|u|\leq|v|$, $|\nabla u|=g^\p(v)|\nabla v|\leq|\nabla v|$, $u^{2\alpha}\leq(2\alpha)^{1/2}|v|$ and $|\nabla(u^{2\alpha})|=2\alpha u^{2\alpha-1}g^\p(v)|\nabla v|\leq(2\alpha)^{1/2}|\nabla v|$. Hence $u,u^{2\alpha}\in H^1_0(\Omega)$, $u\in L^\infty_{loc}(\Omega)$, and $u$ inherits the positivity of $v$. Let $\phi\in H^1_0(\Omega)\cap L^\infty(\Omega)$ have compact support $K$ and put $w:=\phi\,(1+2\alpha u^{4\alpha-2})^{1/2}=\phi/g^\p(v)$. Since $u\in C^1(K)$ is bounded below on $K$, $w\in H^1_0(\Omega)\cap L^\infty(\Omega)$ has compact support, so \eqref{sd7} holds with this $w$. Using $\nabla v=(1+2\alpha u^{4\alpha-2})^{1/2}\nabla u$ and
		\begin{equation*}
			\nabla w=(1+2\alpha u^{4\alpha-2})^{1/2}\nabla\phi+\alpha(4\alpha-2)\,u^{4\alpha-3}(1+2\alpha u^{4\alpha-2})^{-1/2}\phi\,\nabla u ,
		\end{equation*}
		we obtain \[\nabla v\cdot\nabla w=(1+2\alpha u^{4\alpha-2})\nabla u\cdot\nabla\phi+2\alpha(2\alpha-1)u^{4\alpha-3}|\nabla u|^2\phi,\] while the right-hand side of \eqref{sd7} equals the right-hand side of \eqref{sd8}, because $g^\p(v)w=\phi$ and $H(y,g(v(y)))=H(y,u(y))$. This is \eqref{sd8}.
		
		Conversely, let $u$ be a weak solution of \eqref{sd1} and set $v:=g^{-1}(u)$, where $(g^{-1})^\p(s)=(1+2\alpha s^{4\alpha-2})^{1/2}\leq1+(2\alpha)^{1/2}s^{2\alpha-1}$ for $s>0$. Then $0\leq v\leq u+(2\alpha)^{-1/2}u^{2\alpha}$ and $|\nabla v|\leq|\nabla u|+(2\alpha)^{-1/2}|\nabla(u^{2\alpha})|$, so $v\in H^1_0(\Omega)$, and $v$ has the positivity required in Definition~\ref{defsd2}. For $w\in C_0^\infty(\Omega)$, put $\phi:=g^\p(v)w=w(1+2\alpha u^{4\alpha-2})^{-1/2}$. On $\operatorname{supp} w$, the function $u$ is bounded and bounded away from zero, so $\phi\in H^1_0(\Omega)\cap L^\infty(\Omega)$ has compact support and is admissible in \eqref{sd8}. Now
		\begin{equation*}
			\nabla\phi=(1+2\alpha u^{4\alpha-2})^{-1/2}\nabla w-\alpha(4\alpha-2)\,u^{4\alpha-3}(1+2\alpha u^{4\alpha-2})^{-3/2}w\,\nabla u.
		\end{equation*}
        When this expression is substituted into the left-hand side of \eqref{sd8}, the two terms involving $u^{4\alpha-3}|\nabla u|^2$ cancel, and we are left with
        \[(1+2\alpha u^{4\alpha-2})\nabla u\cdot\nabla\phi+2\alpha(2\alpha-1)u^{4\alpha-3}|\nabla u|^2\phi=(1+2\alpha u^{4\alpha-2})^{1/2}\nabla u\cdot\nabla w=\nabla v\cdot\nabla w.\] The right-hand side of \eqref{sd8} equals $\lambda\int_\Omega\bigl(\rho(v)+\mathcal{K}[v]\,h(x,g(v))g^\p(v)\bigr)w\,dx$, the right-hand side of \eqref{sd7}, and \eqref{sd7} follows. \qed
	\end{proof}
	
	\subsection{Proof of Lemma~\ref{lemsd3}(c)}
	\begin{proof}
		 For $\epsilon\in(0,1]$, define
		\begin{align*}
			P_\epsilon(t):=
			\begin{cases}
				\int_0^t\rho(s+\epsilon)ds, & t\geq0,\\
				\rho(\epsilon)t, & t<0,
			\end{cases}
			\qquad
			E_\epsilon(v):=\frac12\int_\Omega|\nabla v|^2dx-\lambda\int_\Omega P_\epsilon(v)dx .
		\end{align*}
		Since $\rho$ is continuous, nonnegative and nonincreasing, $P_\epsilon$ is of class $C^1$ and concave, with
		\begin{align*}
			P_\epsilon^\p(t)=\rho(t+\epsilon)\ \ (t\geq0),\quad P_\epsilon^\p(t)=\rho(\epsilon)\ \ (t<0),\quad 0\leq P_\epsilon^\p\leq\rho(\epsilon),\quad |P_\epsilon(t)|\leq\rho(\epsilon)|t| .
		\end{align*}
		Hence $E_\epsilon$ is strictly convex, coercive, and weakly lower semicontinuous on $H^1_0(\Omega)$, and it has a unique minimizer $z_\epsilon$. The minimizer solves $-\Delta z_\epsilon=\lambda P_\epsilon^\p(z_\epsilon)$ weakly. Since $P_\epsilon^\p\geq0$ is bounded, elliptic regularity \cite[Chapter 9]{gilbarg2001} yields $z_\epsilon\in W^{2,p}(\Omega)\cap C^1(\bar\Omega)$ for all $p<\infty$, and the maximum principle gives $z_\epsilon\geq0$. Consequently, $-\Delta z_\epsilon=\lambda\rho(z_\epsilon+\epsilon)$. The comparison argument of (b), now applied to bounded functions, shows that $z_\epsilon$ is nondecreasing as $\epsilon$ decreases.
		
		Let $\zeta^{(1)}$ be the solution of
		\begin{align*}
			-\Delta\zeta=\lambda C_\rho\zeta^{-\beta}\ \text{ in }\Omega,\quad \zeta>0\ \text{ in }\Omega,\quad \zeta=0\ \text{ on }\partial\Omega,
		\end{align*}
		given by \cite{crandall1977dirichlet}. By \cite{lazer1991}, $\zeta^{(1)}\in H^1_0(\Omega)\cap C(\bar\Omega)$ precisely because $\beta<3$, and (a), applied with $\rho(t)=C_\rho t^{-\beta}$, gives $\zeta^{(1)}\geq\kappa_1\varphi_\beta$ for some $\kappa_1>0$. Let $z_0$ solve $-\Delta z_0=1$ in $\Omega$, $z_0=0$ on $\partial\Omega$, and put $\zeta_\lambda:=\zeta^{(1)}+\lambda C_\rho z_0$. Then
		\begin{align*}
			-\Delta\zeta_\lambda=\lambda C_\rho\bigl((\zeta^{(1)})^{-\beta}+1\bigr)\geq\lambda C_\rho(\zeta_\lambda^{-\beta}+1)\geq\lambda\rho(\zeta_\lambda+\epsilon),
		\end{align*}
		so $\zeta_\lambda$ is a supersolution of the problem solved by $z_\epsilon$. Testing with $(z_\epsilon-\zeta_\lambda)^+$, which is admissible by the estimate in (b), we get $z_\epsilon\leq\zeta_\lambda$.
		
       By \eqref{sd15}, $w:=\kappa\varphi_\beta-\epsilon$ satisfies $-\Delta w\leq\lambda\rho(w+\epsilon)$ in $\Omega$, and $(w-z_\epsilon)^+$ is supported in $\{\kappa\varphi_\beta\geq\epsilon\}\Subset\Omega$. The argument of (a) gives $z_\epsilon+\epsilon\geq\kappa\varphi_\beta$. Consequently $\rho(z_\epsilon+\epsilon)\leq C_\rho((\kappa\varphi_\beta)^{-\beta}+1)$, and testing the equation with $z_\epsilon$,
		\begin{align*}
			\int_\Omega|\nabla z_\epsilon|^2dx=\lambda\int_\Omega\rho(z_\epsilon+\epsilon)z_\epsilon\,dx\leq\lambda C_\rho\int_\Omega\bigl((\kappa\varphi_\beta)^{-\beta}+1\bigr)z_\epsilon\,dx\leq C\|\nabla z_\epsilon\|_2
		\end{align*}
		by the estimate of Remark~\ref{remsd1}. Hence $\{z_\epsilon\}$ is bounded in $H^1_0(\Omega)$.
		
	    As $\epsilon\downarrow0$, $z_\epsilon\uparrow z$ a.e.\ for some $z$ with $\kappa\varphi_\beta\leq z\leq\zeta_\lambda$, and $z_\epsilon\rightharpoonup z$ in $H^1_0(\Omega)$. Let $\phi\in C^\infty_0(\Omega)$. Since $\rho$ is continuous, $\rho(z_\epsilon+\epsilon)\to\rho(z)$ a.e., and these functions are bounded by $C_\rho((\kappa\varphi_\beta)^{-\beta}+1)\in L^\infty(\operatorname{supp}\phi)$. By dominated convergence,
		\begin{align*}
			\int_\Omega\nabla z\cdot\nabla\phi\,dx=\lambda\int_\Omega\rho(z)\phi\,dx .
		\end{align*}
		Since $\rho(z)\in L^\infty_{loc}(\Omega)$, elliptic regularity gives $z\in W^{2,p}_{loc}(\Omega)$, hence $z\in C(\Omega)$. Finally, $0\leq z\leq\zeta_\lambda\in C(\bar\Omega)$ and $\zeta_\lambda=0$ on $\partial\Omega$ give $z\in C(\bar\Omega)$. \qed
	\end{proof}
	
	\subsection{Proof of Lemma~\ref{lemsd20}}
	\begin{proof}
		Put $t^*:=\|\bar v\|_{L^\infty(\Omega)}$. Recall from Subsection~\ref{secsd1} that $\Psi_t(x,t)=h(x,g(t))g^\p(t)$ is continuous and nondecreasing in $t\geq0$, so it is bounded on $\bar\Omega\times[0,t^*]$.
		By \eqref{sd72}, the singular term is controlled on $\mathcal{N}_*$. Moreover, with $C_\mu$ as in Subsection~\ref{secsd1},
		\begin{align*}
			0\leq\mathcal{K}[v]\Psi_t(x,v)\leq C_\mu\max_{\bar\Omega\times[0,t^*]}\Psi\,\max_{\bar\Omega\times[0,t^*]}\Psi_t\qquad(v\in\mathcal{N}_*).
		\end{align*}
		
		Since $\rho_*\geq0$, $P_*$ is nondecreasing, so
		\begin{align*}
			P_*(\underline v)\leq P_*(v)\leq P_*(\bar v)\leq P_*(t^*)\quad\text{on }\mathcal{N}_* .
		\end{align*}
		We claim that $P_*(\underline v)\in L^1(\Omega)$. For $0<t\leq1$,
		\begin{align*}
			P_*(t)=-\int_t^1\rho_*(s)\,ds\geq
			\begin{cases}
				-C_\rho\bigl(\frac{t^{1-\beta}}{\beta-1}+1\bigr), & \beta>1,\\[1mm]
				-C_\rho\bigl(\log\frac1t+1\bigr), & \beta=1,\\[1mm]
				-C_\rho\bigl(\frac{1}{1-\beta}+1\bigr), & \beta<1.
			\end{cases}
		\end{align*}
		For $\beta>1$, we have $\underline v\geq c\varphi_\beta\geq c_1d(x)^{2/(\beta+1)}$, so $|P_*(\underline v)|\leq C(1+d(x)^{-2(\beta-1)/(\beta+1)})$. This function is integrable because $\frac{2(\beta-1)}{\beta+1}<1$, that is, $\beta<3$. For $\beta=1$, we have $|P_*(\underline v)|\leq C(1+\log\frac1{d(x)})\in L^1(\Omega)$, and for $\beta<1$, the function $P_*(\underline v)$ is bounded. This proves the claim. Finally, $0\leq\Psi(\cdot,v)\leq\Psi(\cdot,\bar v)\leq C$ on $\mathcal{N}_*$. Hence, for $v\in\mathcal{N}_*$,
		\begin{align}\label{sd73}
			\mathcal{E}(v)\geq\frac12\|\nabla v\|_2^2-\lambda\int_\Omega P_*(\bar v)\,dx-\frac\lambda2\int_\Omega\int_\Omega\frac{\Psi(y,\bar v)\Psi(x,\bar v)}{|x-y|^\mu}dxdy\geq\frac12\|\nabla v\|_2^2-C_* ,
		\end{align}
		and $\mathcal{E}$ is finite on $\mathcal{N}_*$.
		The set $\mathcal{N}_*$ is closed and convex, hence weakly closed. Let $v_k\rightharpoonup v$ in $H^1_0(\Omega)$ with $v_k\in\mathcal{N}_*$. Along a subsequence, $v_k\to v$ a.e. The two nonlinear terms of $\mathcal{E}$ converge by dominated convergence, with the majorants $|P_*(\underline v)|+|P_*(\bar v)|$ and $\Psi(y,\bar v)\Psi(x,\bar v)|x-y|^{-\mu}$, and $\liminf\|\nabla v_k\|_2^2\geq\|\nabla v\|_2^2$. Thus $\mathcal{E}$ is weakly lower semicontinuous on $\mathcal{N}_*$. It is coercive by \eqref{sd73}, and therefore it attains its minimum at some $u\in\mathcal{N}_*$.
		
		Let $u$ be a minimizer, $\phi\in C^\infty_0(\Omega)$ and $\epsilon>0$. Put $K:=\operatorname{supp}\phi$ and $c_K:=\min_Kc\varphi_\beta>0$, and
		\begin{align*}
			\overline\vartheta_\epsilon:=(u+\epsilon\phi-\bar v)^+,\qquad \underline\vartheta_\epsilon:=(u+\epsilon\phi-\underline v)^-,\\
			\ell_\epsilon:=\min\{\bar v,\max\{\underline v,u+\epsilon\phi\}\}=u+\epsilon\phi-\overline\vartheta_\epsilon+\underline\vartheta_\epsilon .
		\end{align*}
		Then $\ell_\epsilon\in\mathcal{N}_*$ and $0\leq\overline\vartheta_\epsilon,\underline\vartheta_\epsilon\leq\epsilon|\phi|$. Since $\underline v\leq u\leq\bar v$, the function $\overline\vartheta_\epsilon$ is supported in $\Omega^\epsilon:=\{u+\epsilon\phi>\bar v\}\subset K\cap\{\phi>0\}$, and $\underline\vartheta_\epsilon$ is supported in $\Omega_\epsilon:=\{u+\epsilon\phi<\underline v\}\subset K\cap\{\phi<0\}$.
		
		Since $\mathcal{N}_*$ is convex, $u+s(\ell_\epsilon-u)\in\mathcal{N}_*$ for $s\in[0,1]$, and therefore $\mathcal{E}(u+s(\ell_\epsilon-u))\geq \mathcal{E}(u)$. Dividing by $s$ and letting $s\downarrow0$, we obtain
		\begin{align}\label{sd74}
			0\leq\langle \mathcal{E}^\p(u),\ell_\epsilon-u\rangle:=\int_\Omega\nabla u\cdot\nabla(\ell_\epsilon-u)\,dx-\lambda\int_\Omega\bigl(\rho_*(u)+\mathcal{K}[u]\Psi_t(x,u)\bigr)(\ell_\epsilon-u)\,dx .
		\end{align}
		To justify the limit, put $w:=\ell_\epsilon-u$. As $s\downarrow0$, $s^{-1}\bigl(P_*(u+sw)-P_*(u)\bigr)\to\rho_*(u)w$ a.e. Since $u+sw\in\mathcal{N}_*$ and $\rho_*$ is nonincreasing, these quotients are bounded by $\rho_*(\underline v)|w|$, which is integrable by \eqref{sd72}. Dominated convergence also shows that the double integral is differentiable along the segment: $\Psi_t$ is continuous and bounded on $\bar\Omega\times[0,t^*]$ and the kernel is symmetric.
		
		By linearity,
		\begin{align*}
			\langle \mathcal{E}^\p(u),\ell_\epsilon-u\rangle=\epsilon\langle \mathcal{E}^\p(u),\phi\rangle-\langle \mathcal{E}^\p(u),\overline\vartheta_\epsilon\rangle+\langle \mathcal{E}^\p(u),\underline\vartheta_\epsilon\rangle .
		\end{align*}
		We claim that
		\begin{align}\label{sd75}
			\langle \mathcal{E}^\p(u),\overline\vartheta_\epsilon\rangle\geq-o(\epsilon),\qquad \langle \mathcal{E}^\p(u),\underline\vartheta_\epsilon\rangle\leq o(\epsilon)\qquad\text{as }\epsilon\downarrow0 .
		\end{align}
		For the first claim, we use the second inequality in \eqref{sd70} with the test function $\overline\vartheta_\epsilon\geq0$:
		\begin{align*}
			\langle \mathcal{E}^\p(u),\overline\vartheta_\epsilon\rangle&\geq\int_\Omega\nabla(u-\bar v)\cdot\nabla\overline\vartheta_\epsilon\,dx-\lambda\int_\Omega\bigl(\rho_*(u)-\rho_*(\bar v)\bigr)\overline\vartheta_\epsilon\,dx\\
			&\qquad-\lambda\int_\Omega\bigl(\mathcal{K}[u]\Psi_t(x,u)-\mathcal{K}[\bar v]\Psi_t(x,\bar v)\bigr)\overline\vartheta_\epsilon\,dx .
		\end{align*}
		We now estimate each of the three terms.
		 The last integral is nonpositive, because $u\leq\bar v$ and $\mathcal{K}[\cdot]$, $\Psi_t(x,\cdot)$ are nondecreasing.
			 For the first term, note that $\nabla\overline\vartheta_\epsilon=\nabla(u-\bar v)+\epsilon\nabla\phi$ on $\Omega^\epsilon$, so
			\begin{align*}
				\int_\Omega\nabla(u-\bar v)\cdot\nabla\overline\vartheta_\epsilon\,dx\geq-\epsilon\int_{\Omega^\epsilon}|\nabla(u-\bar v)||\nabla\phi|\,dx=o(\epsilon),
			\end{align*}
			because $\chi_{\Omega^\epsilon}\to\chi_{\{u=\bar v\}\cap\{\phi>0\}}$ a.e.\ as $\epsilon\downarrow0$ and $\nabla(u-\bar v)=0$ a.e.\ on $\{u=\bar v\}$.
			For the second term, on $\Omega^\epsilon$ we have $0\leq\bar v-u<\epsilon\|\phi\|_\infty$ and $c_K\leq u\leq\bar v\leq t^*$. Let $\omega_K$ be a modulus of continuity of $\rho_*$ on $[c_K,t^*]$. Since $\rho_*(u)\geq\rho_*(\bar v)$,
			\begin{align*}
				0\leq\lambda\int_\Omega\bigl(\rho_*(u)-\rho_*(\bar v)\bigr)\overline\vartheta_\epsilon\,dx\leq\lambda\,\omega_K(\epsilon\|\phi\|_\infty)\,\epsilon\|\phi\|_{L^1(\Omega)}=o(\epsilon).
			\end{align*}
		
		This proves the first claim in \eqref{sd75}. The second is proved in the same way, with the first inequality in \eqref{sd70} and the test function $\underline\vartheta_\epsilon\geq0$. Here $\nabla\underline\vartheta_\epsilon=\nabla(\underline v-u)-\epsilon\nabla\phi$ on $\Omega_\epsilon$, the difference $\mathcal{K}[\underline v]\Psi_t(x,\underline v)-\mathcal{K}[u]\Psi_t(x,u)$ is nonpositive, and on $\Omega_\epsilon$, we have $0\leq u-\underline v<\epsilon\|\phi\|_\infty$, with $u,\underline v\in[c_K,t^*]$.
		
		Inserting \eqref{sd75} into \eqref{sd74}, dividing by $\epsilon$ and letting $\epsilon\downarrow0$, we get $\langle \mathcal{E}^\p(u),\phi\rangle\geq0$. Replacing $\phi$ by $-\phi$ gives $\langle \mathcal{E}^\p(u),\phi\rangle=0$ for all $\phi\in C^\infty_0(\Omega)$. By \eqref{sd72} and the boundedness of $\mathcal{K}[u]\Psi_t(\cdot,u)$, the map $\phi\mapsto\langle \mathcal{E}^\p(u),\phi\rangle$ is continuous on $H^1_0(\Omega)$, so the identity holds for all $\phi\in H^1_0(\Omega)$. Finally, $\operatorname*{ess\,inf}_Ku\geq\min_Kc\varphi_\beta>0$ for every compact $K\subset\Omega$. \qed
	\end{proof}
	
	\subsection{Proof of Lemma~\ref{lemsd19}}
	\begin{proof}
		With $d_{\mathcal{C}}(\zeta,\zeta^\p):=\max_{[0,1]}\|\nabla(\zeta-\zeta^\p)\|_2$, $(\mathcal{C},d_{\mathcal{C}})$ is a complete metric space, since $X^+$ is closed. The functional $\Upsilon(\zeta):=\max_{[0,1]}\mathcal{J}_\lambda\circ\zeta$ is continuous on it, because $\mathcal{J}_\lambda=\Lambda+\psi$ is locally Lipschitz and $\zeta([0,1])$ is compact. Ekeland's variational principle \cite{ekeland1974} gives $\zeta_k\in\mathcal{C}$ with
		\begin{align}\label{sd64}
			\Upsilon(\zeta_k)\leq c+\tfrac1k,\qquad \Upsilon(\zeta)\geq\Upsilon(\zeta_k)-\tfrac1k\,d_{\mathcal{C}}(\zeta,\zeta_k)\quad\text{for all }\zeta\in\mathcal{C}.
		\end{align}
		The set $\Gamma_k:=\{t\in[0,1]:\mathcal{J}_\lambda(\zeta_k(t))=\Upsilon(\zeta_k)\}$ is compact. Since $\Upsilon(\zeta_k)\geq c>\max\{\mathcal{J}_\lambda(0),\mathcal{J}_\lambda(e)\}$, it does not contain $0$ or $1$. It suffices to find $t_k\in\Gamma_k$ such that $w_k:=\zeta_k(t_k)$ satisfies \eqref{sd63}; then $\mathcal{J}_\lambda(w_k)=\Upsilon(\zeta_k)\in[c,c+\frac1k]$.
		
	    Suppose no such $t_k$ exists. For each $t\in\Gamma_k$, pick $\ell_t\in X^+$ violating \eqref{sd63} at $w=\zeta_k(t)$, put $\mathfrak{m}_t:=\max\{1,\|\nabla(\ell_t-\zeta_k(t))\|_2\}$ and, for $\tau\in[0,1]$,
		\begin{align*}
			e_t(\tau):=\frac{\ell_t-\zeta_k(\tau)}{\mathfrak{m}_t},\qquad q_t(\tau):=\frac{\langle\Lambda^\p(\zeta_k(\tau)),\ell_t-\zeta_k(\tau)\rangle+\psi(\ell_t)-\psi(\zeta_k(\tau))}{\mathfrak{m}_t} .
		\end{align*}
		Then $\delta_t:=-\frac1k-q_t(t)>0$ and $\|\nabla e_t(t)\|_2\leq1$. The maps $\tau\mapsto q_t(\tau)$ and $\tau\mapsto e_t(\tau)$ are continuous, because $\Lambda^\p$, $\psi$ and $\zeta_k$ are continuous. Hence, there is an open interval $I_t\ni t$ such that
		\begin{align*}
			q_t(\tau)<-\tfrac1k-\tfrac{\delta_t}2,\qquad \|\nabla e_t(\tau)\|_2<1+\tfrac{k\delta_t}8\qquad\text{for }\tau\in I_t\cap[0,1].
		\end{align*}
		
		Cover $\Gamma_k$ by finitely many $I_{t_1},\dots,I_{t_m}$, and write $e_j:=e_{t_j}$, $q_j:=q_{t_j}$, $\ell_j:=\ell_{t_j}$, $\mathfrak{m}_j:=\mathfrak{m}_{t_j}$, $\delta_j:=\delta_{t_j}$, $\delta_*:=\min_j\delta_j$, $\delta^*:=\max_j\delta_j$. First, choose closed neighborhoods $N^\p\subset N$ of $\Gamma_k$ in $(0,1)$ with $N^\p\subset\operatorname{int}N$ and $N\subset\bigcup_jI_{t_j}$. Next, choose continuous functions $\pi_j\geq0$ with $\operatorname{supp}\pi_j\subset I_{t_j}$ and $\sum_j\pi_j=1$ on $N$. Finally, take a continuous function $\chi:[0,1]\to[0,1]$ with $\chi=1$ on $N^\p$ and $\operatorname{supp}\chi\subset N$. 
		Put
		\begin{align*}
       W(\tau):=\sum_{j=1}^m\pi_j(\tau)e_j(\tau),\quad \widehat W(\tau):=\frac{W(\tau)}{\max\{1,\|\nabla W(\tau)\|_2\}} .
		\end{align*}
		For $\tau\in N$, with $S:=\sum_j\pi_j(\tau)\delta_j\in[\delta_*,\delta^*]$,
		\begin{align}\label{sd65}
			\sum_{j=1}^m\pi_j(\tau)q_j(\tau)\leq-\tfrac1k-\tfrac S2,\quad \|\nabla W(\tau)\|_2\leq1+\tfrac{kS}8 ,\quad \|\nabla\widehat W(\tau)\|_2\leq1 .
		\end{align}
		For $s\in[0,1]$, let $\zeta^s(\tau):=\zeta_k(\tau)+s\chi(\tau)\widehat W(\tau)$. On $\operatorname{supp}\chi$, put $\varsigma(\tau):=s\chi(\tau)/\max\{1,\|\nabla W(\tau)\|_2\}\in[0,1]$. Then
		\begin{align}\label{sd66}
			\zeta^s(\tau)=\sum_{j=1}^m\pi_j(\tau)\Bigl(\zeta_k(\tau)+\frac{\varsigma(\tau)}{\mathfrak{m}_{j}}\bigl(\ell_{j}-\zeta_k(\tau)\bigr)\Bigr)
		\end{align}
		is a convex combination of elements of $X^+$, because $\varsigma/\mathfrak{m}_{j}\in[0,1]$. Hence $\zeta^s\in\mathcal{C}$ and $d_{\mathcal{C}}(\zeta^s,\zeta_k)\leq s$, and \eqref{sd64} gives
		\begin{align*}
			\Upsilon(\zeta^s)\geq\Upsilon(\zeta_k)-\frac sk .
		\end{align*}
		
		Let $\tau\in N$. We apply the convexity of $\psi$ to \eqref{sd66}, first to the convex combination and then to each term, and obtain
		\begin{align*}
			\psi(\zeta^s(\tau))-\psi(\zeta_k(\tau))\leq\sum_{j=1}^m\pi_j(\tau)\frac{\varsigma(\tau)}{\mathfrak{m}_j}\bigl(\psi(\ell_j)-\psi(\zeta_k(\tau))\bigr).
		\end{align*}
		By the mean value theorem for the $C^1$ functional $\Lambda$, there is $y_\tau^s$ on the segment from $\zeta_k(\tau)$ to $\zeta^s(\tau)$ such that
		\begin{align*}
			\Lambda(\zeta^s(\tau))-\Lambda(\zeta_k(\tau))=s\chi(\tau)\langle\Lambda^\p(y^s_\tau),\widehat W(\tau)\rangle
			\leq s\chi(\tau)\langle\Lambda^\p(\zeta_k(\tau)),\widehat W(\tau)\rangle+s\,\varpi(s).
		\end{align*}
		Here $\varpi(s):=\sup\{\|\Lambda^\p(y)-\Lambda^\p(z)\|:z\in\zeta_k([0,1]),\ \|\nabla(y-z)\|_2\leq s\}$, and we used $\|\nabla\widehat W\|_2\leq1$. Since $\Lambda^\p$ is continuous and $\zeta_k([0,1])$ is compact, $\varpi(s)\to0$ as $s\downarrow0$. Moreover,
		\begin{align*}
			s\chi(\tau)\widehat W(\tau)=\varsigma(\tau)\sum_j\pi_j(\tau)\frac{\ell_j-\zeta_k(\tau)}{\mathfrak{m}_j}.
		\end{align*}
		Adding the two estimates and using \eqref{sd65}, we get
		\begin{align*}
			\mathcal{J}_\lambda(\zeta^s(\tau))-\mathcal{J}_\lambda(\zeta_k(\tau))\leq\varsigma(\tau)\sum_{j=1}^m\pi_j(\tau)q_j(\tau)+s\varpi(s)\leq-\varsigma(\tau)\Bigl(\frac1k+\frac S2\Bigr)+s\varpi(s)\qquad(\tau\in N).
		\end{align*}
		For $\tau\in N^\p$, we have $\chi(\tau)=1$, hence $\varsigma(\tau)\geq s/(1+\frac{kS}8)$ by \eqref{sd65}. Moreover,
		\begin{align*}
			\frac{\frac1k+\frac S2}{1+\frac{kS}8}=\frac1k+\frac{3S/8}{1+kS/8}\geq\frac1k+\hat\delta,\quad \hat\delta:=\frac{3\delta_*/8}{1+k\delta^*/8}>0 .
		\end{align*}
		Hence
		\begin{align}\label{sd67}
			\mathcal{J}_\lambda(\zeta^s(\tau))\leq\mathcal{J}_\lambda(\zeta_k(\tau))-s\Bigl(\frac1k+\hat\delta\Bigr)+s\varpi(s)\qquad(\tau\in N^\p,\ s\in(0,1]).
		\end{align}
		
		Let $\tau_s$ be a maximum point of $\mathcal{J}_\lambda\circ\zeta^s$. As $s\downarrow0$, $\zeta^s\to\zeta_k$ uniformly, so every limit point of $\tau_s$ lies in $\Gamma_k$. Take $s_i\downarrow0$ with $\tau_{s_i}\to\tau_*\in\Gamma_k$. Then $\tau_{s_i}\in N^\p$ and $\varpi(s_i)\leq\hat\delta/2$ for $i$ large. By \eqref{sd67},
		\begin{align*}
			\Upsilon(\zeta^{s_i})=\mathcal{J}_\lambda(\zeta^{s_i}(\tau_{s_i}))\leq\mathcal{J}_\lambda(\zeta_k(\tau_{s_i}))-s_i\Bigl(\tfrac1k+\tfrac{\hat\delta}2\Bigr)\leq\Upsilon(\zeta_k)-s_i\Bigl(\tfrac1k+\tfrac{\hat\delta}2\Bigr),
		\end{align*}
		which contradicts $\Upsilon(\zeta^{s_i})\geq\Upsilon(\zeta_k)-\frac{s_i}k$. \qed
	\end{proof}
	\noindent \textbf{Funding} \\
The authors declare that no funds, grants, or other support were received during the preparation of this manuscript.

\noindent \textbf{Competing Interests} \\
The authors have no relevant financial or non-financial interests to disclose.

\noindent \textbf{Ethics Approval and Informed Consent} \\
Not applicable.

	\bibliographystyle{siam}
	\bibliography{ref}
\end{document}